\documentclass{amsart}
\usepackage{amscd,amsmath,amssymb,amsfonts, dsfont}
\usepackage[cmtip, all]{xy}
\usepackage{pstricks}
\usepackage{color}

\begin{document}
\newtheorem{theorem}{Theorem}[subsection]
\newtheorem{proposition}[theorem]{Proposition}
\newtheorem{lemma}[theorem]{Lemma}
\newtheorem{corollary}[theorem]{Corollary}
\newtheorem{conjecture}[theorem]{Conjecture}
\newtheorem{IntroThm}{Theorem}
\newtheorem{IntroConj}[IntroThm]{Conjecture}
\newtheorem{claim}{Claim}
\renewcommand{\theclaim}{\kern-4pt}
\theoremstyle{definition}
\newtheorem{definition}[theorem]{Definition}
\theoremstyle{remark}
\newtheorem{remark}[theorem]{Remark}
\newtheorem{SecRemark}{Remark}
\renewcommand{\theSecRemark}{\kern-4.5pt}
\newtheorem{remarks}[theorem]{Remarks}
\newtheorem{example}[theorem]{Example}
\newtheorem{examples}[theorem]{Examples}
\newtheorem{nota}[theorem]{Notation}
\numberwithin{equation}{subsection}
\newcommand{\eq}[2]{\begin{equation}\label{#1}#2 \end{equation}}
\newcommand{\ml}[2]{\begin{multline}\label{#1}#2 \end{multline}}
\newcommand{\ga}[2]{\begin{gather}\label{#1}#2 \end{gather}}
\newcommand{\sA}{{\mathcal A}}
\newcommand{\sB}{{\mathcal B}}
\newcommand{\sC}{{\mathcal C}}
\newcommand{\sD}{{\mathcal D}}
\newcommand{\sE}{{\mathcal E}}
\newcommand{\sF}{{\mathcal F}}
\newcommand{\sG}{{\mathcal G}}
\newcommand{\sH}{{\mathcal H}}
\newcommand{\sI}{{\mathcal I}}
\newcommand{\sJ}{{\mathcal J}}
\newcommand{\sK}{{\mathcal K}}
\newcommand{\sL}{{\mathcal L}}
\newcommand{\sM}{{\mathcal M}}
\newcommand{\sN}{{\mathcal N}}
\newcommand{\sO}{{\mathcal O}}
\newcommand{\sP}{{\mathcal P}}
\newcommand{\sQ}{{\mathcal Q}}
\newcommand{\sR}{{\mathcal R}}
\newcommand{\sS}{{\mathcal S}}
\newcommand{\sT}{{\mathcal T}}
\newcommand{\sU}{{\mathcal U}}
\newcommand{\sV}{{\mathcal V}}
\newcommand{\sW}{{\mathcal W}}
\newcommand{\sX}{{\mathcal X}}
\newcommand{\sY}{{\mathcal Y}}
\newcommand{\sZ}{{\mathcal Z}}
\newcommand{\A}{{\mathbb A}}
\newcommand{\B}{{\mathbb B}}
\newcommand{\C}{{\mathbb C}}
\newcommand{\F}{{\mathbb F}}
\newcommand{\G}{{\mathbb G}}
\renewcommand{\H}{{\mathbb H}}
\newcommand{\J}{{\mathbb J}}
\newcommand{\M}{{\mathbb M}}
\newcommand{\N}{{\mathbb N}}
\renewcommand{\P}{{\mathbb P}}
\newcommand{\Q}{{\mathbb Q}}
\newcommand{\R}{{\mathbb R}}
\newcommand{\T}{{\mathbb T}}
\newcommand{\U}{{\mathbb U}}
\newcommand{\V}{{\mathbb V}}
\newcommand{\W}{{\mathbb W}}
\newcommand{\X}{{\mathbb X}}
\newcommand{\Y}{{\mathbb Y}}
\newcommand{\Z}{{\mathbb Z}}
\renewcommand{\1}{{\mathds{1}}}

\newcommand{\fa}{{\mathfrak{a}}}
\newcommand{\fb}{{\mathfrak{b}}}
\newcommand{\fc}{{\mathfrak{c}}}
\newcommand{\fd}{{\mathfrak{d}}}
\newcommand{\fe}{{\mathfrak{e}}}
\newcommand{\ff}{{\mathfrak{f}}}
\newcommand{\fg}{{\mathfrak{g}}}
\newcommand{\fh}{{\mathfrak{h}}}
\newcommand{\fj}{{\mathfrak{j}}}
\newcommand{\fk}{{\mathfrak{k}}}
\newcommand{\fl}{{\mathfrak{l}}}
\newcommand{\fm}{{\mathfrak{m}}}
\newcommand{\fn}{{\mathfrak{n}}}
\newcommand{\fo}{{\mathfrak{o}}}
\newcommand{\fp}{{\mathfrak{p}}}
\newcommand{\fq}{{\mathfrak{q}}}
\newcommand{\fr}{{\mathfrak{r}}}
\newcommand{\fs}{{\mathfrak{s}}}
\newcommand{\ft}{{\mathfrak{t}}}
\newcommand{\fu}{{\mathfrak{u}}}
\newcommand{\fv}{{\mathfrak{v}}}
\newcommand{\fw}{{\mathfrak{w}}}
\newcommand{\fx}{{\mathfrak{x}}}
\newcommand{\fy}{{\mathfrak{y}}}
\newcommand{\fz}{{\mathfrak{z}}}

\newcommand{\an}{{\rm an}}
\newcommand{\alg}{{\rm alg}}
\newcommand{\cl}{{\rm cl}}
\newcommand{\Alb}{{\rm Alb}}
\newcommand{\CH}{{\rm CH}}
\newcommand{\mc}{\mathcal}
\newcommand{\mb}{\mathbb}
\newcommand{\surj}{\twoheadrightarrow}
\newcommand{\codim}{{\rm codim}}
\newcommand{\rank}{{\rm rank}}
\newcommand{\Pic}{{\rm Pic}}
\newcommand{\Div}{{\rm Div}}
\newcommand{\Hom}{{\rm Hom}}
\newcommand{\im}{{\rm im}}
\newcommand{\Spec}{{\rm Spec \,}}
\newcommand{\sing}{{\rm sing}}
\newcommand{\Char}{{\rm char}}
\newcommand{\Tr}{{\rm Tr}}
\newcommand{\Gal}{{\rm Gal}}
\newcommand{\Min}{{\rm Min \ }}
\newcommand{\Max}{{\rm Max \ }}
\newcommand{\supp}{{\rm supp}\,}
\newcommand{\0}{\emptyset}
\newcommand{\sHom}{{\mathcal{H}{om}}}

\newcommand{\Nm}{{\operatorname{Nm}}}
\newcommand{\NS}{{\operatorname{NS}}}
\newcommand{\id}{{\operatorname{id}}}
\newcommand{\Zar}{{\text{\rm Zar}}} 
\newcommand{\Ord}{{\mathbf{Ord}}}
\newcommand{\FSimp}{{{\mathcal FS}}}
\newcommand{\inj}{{\text{\rm inj}}}  
\newcommand{\Sch}{{\operatorname{\mathbf{Sch}}}} 
\newcommand{\cosk}{{\operatorname{\rm cosk}}} 
\newcommand{\sk}{{\operatorname{\rm sk}}} 
\newcommand{\subv}{{\operatorname{\rm sub}}}
\newcommand{\bary}{{\operatorname{\rm bary}}}
\newcommand{\Comp}{{\mathbf{SC}}}
\newcommand{\IComp}{{\mathbf{sSC}}}
\newcommand{\Top}{{\mathbf{Top}}}
\newcommand{\holim}{\mathop{{\rm holim}}}
\newcommand{\op}{{\text{\rm op}}}
\newcommand{\<}{\mathopen<}
\renewcommand{\>}{\mathclose>}
\newcommand{\Sets}{{\mathbf{Sets}}}
\newcommand{\del}{\partial}
\newcommand{\fib}{{\operatorname{\rm fib}}}
\renewcommand{\max}{{\operatorname{\rm max}}}
\newcommand{\bad}{{\operatorname{\rm bad}}}
\newcommand{\Spt}{{\mathbf{Spt}}}
\newcommand{\Spc}{{\mathbf{Spc}}}
\newcommand{\Sm}{{\mathbf{Sm}}}
\newcommand{\cofib}{{\operatorname{\rm cofib}}}
\newcommand{\hocolim}{\mathop{{\rm hocolim}}}
\newcommand{\Glu}{{\mathbf{Glu}}}
\newcommand{\can}{{\operatorname{\rm can}}}
\newcommand{\Ho}{{\mathbf{Ho}}}
\newcommand{\GL}{{\operatorname{\rm GL}}}
\newcommand{\sq}{\square}
 \newcommand{\Ab}{{\mathbf{Ab}}}
\newcommand{\Tot}{{\operatorname{\rm Tot}}}
\newcommand{\loc}{{\operatorname{\rm s.l.}}}
\newcommand{\HZ}{{\operatorname{\sH \Z}}}
\newcommand{\Cyc}{{\operatorname{\rm Cyc}}}
\newcommand{\cyc}{{\operatorname{\rm cyc}}}
 \newcommand{\RCyc}{{\operatorname{\widetilde{\rm Cyc}}}}
 \newcommand{\Rcyc}{{\operatorname{\widetilde{\rm cyc}}}}
\newcommand{\Sym}{{\operatorname{\rm Sym}}}
\newcommand{\fin}{{\operatorname{\rm fin}}}
\newcommand{\qfin}{{\operatorname{\rm q.fin}}}
\newcommand{\SH}{{\operatorname{\sS\sH}}}
\newcommand{\Anna}{{$\mathcal{ANNA\ SARAH\ LEVINE}$}}
\newcommand{\Wedge}{{\Lambda}}
\newcommand{\eff}{{\operatorname{\rm eff}}}
\newcommand{\rcyc}{{\operatorname{\rm rev}}}
\newcommand{\DM}{DM}
\newcommand{\GW}{{\operatorname{\rm{GW}}}}
\newcommand{\sSets}{{\mathbf{sSets}}}
\newcommand{\Nis}{{\operatorname{Nis}}}
\newcommand{\Cat}{{\mathbf{Cat}}}
\newcommand{\ds}{{/\kern-3pt/}}
\newcommand{\bD}{{\mathbf{D}}}
\newcommand{\res}{{\operatorname{res}}}

\newcommand{\fil}{\phi}
\newcommand{\barfil}{\sigma}

\newcommand{\Fac}{{\mathop{\rm{Fac}}}}
\newcommand{\Fun}{{\mathbf{Func}}}
\newcommand{\lci}{\text{l.c.i.}}
\newcommand{\lng}{\operatorname{lng}}
\newcommand{\Tor}{{\operatorname{Tor}}}
\newcommand{\Proj}{{\operatorname{Proj}}}
\newcommand{\BM}{{\operatorname{B.M.}}}
\newcommand{\qf}{{\operatorname{q.fin.}}}
\newcommand{\sgn}{{\operatorname{sgn}}}
\newcommand{\Cone}{{\operatorname{Cone}}}
\newcommand{\cone}{{\operatorname{Cone}}}
\newcommand{\BGL}{{\operatorname{BGL}}}
\newcommand{\ch}{{\operatorname{ch}}}
\newcommand{\Gr}{{\operatorname{Gr}}}
\newcommand{\Sus}{{\operatorname{Sus}}}
\newcommand{\FS}{{\operatorname{FS}}}
\newcommand{\Bl}{{\operatorname{Bl}}}
\newcommand{\rat}{{\operatorname{rat}}}
\newcommand{\num}{{\operatorname{num}}}
\renewcommand{\hom}{{\operatorname{hom}}}
\newcommand{\ord}{{\operatorname{ord}}}
\newcommand{\rnk}{{\operatorname{rnk}}}
\newcommand{\et}{{\text{\'et}}}
\newcommand{\dis}{\displaystyle}
\newcommand{\tor}{{\operatorname{tor}}}
\newcommand{\SL}{{\operatorname{SL}}}
\newcommand{\Jac}{{\operatorname{Jac}}}
\newcommand{\Ext}{{\operatorname{Ext}}}
\newcommand{\HS}{{\operatorname{HS}}}
\newcommand{\MHS}{{\operatorname{MHS}}}
\newcommand{\gr}{{\operatorname{gr}}}
\newcommand{\cb}{{\operatorname{cb}}}
\newcommand{\Alt}{{\operatorname{Alt}}}
\newcommand{\HCor}{{\operatorname{HCor}}}
\newcommand{\SmProj}{{\mathbf{SmProj}}}

\newcommand{\hZ}{\Z}

\newcommand{\Cor}{{\operatorname{Cor}}}

\newcommand{\PST}{{\operatorname{PST}}}
\newcommand{\tr}{{\operatorname{tr}}}
\newcommand{\DTM}{{\operatorname{DMT}}}
\newcommand{\DMT}{{\operatorname{DMT}}}
\newcommand{\TM}{{\operatorname{TM}}}
\newcommand{\Vect}{{\operatorname{Vec}}}

\newcommand{\BTM}{{\operatorname{BTM}}}
\newcommand{\BKTM}{{\operatorname{BKTM}}}
\newcommand{\ind}{{\operatorname{ind}}}
\newcommand{\DCM}{{\operatorname{DCM}}}
\newcommand{\DT}{{\operatorname{DT}}}
\newcommand{\coRep}{{\operatorname{co-rep}}}
\newcommand{\sym}{{\operatorname{sym}}}
\newcommand{\odd}{{\operatorname{odd}}}
\newcommand{\ev}{{\operatorname{ev}}}
\newcommand{\Mod}{{\operatorname{mod}}}
\newcommand{\End}{{\operatorname{End}}}
\newcommand{\ECM}{{\operatorname{ECM}}}
\newcommand{\EHM}{{\operatorname{EHM}}}
\newcommand{\Rep}{{\operatorname{Rep}}}
\newcommand{\Period}{{\operatorname{Period}}}
\newcommand{\dash}{{\text{-}}}
\newcommand{\Aff}{{\mathbf{Aff}}}
\newcommand{\NM}{{\operatorname{NMM}}}
\newcommand{\MT}{{\operatorname{MT}}}
\newcommand{\mR}{{\operatorname{MR}}}
\newcommand{\dR}{{\operatorname{DR}}}
\newcommand{\Aut}{{\operatorname{Aut}}}
\newcommand{\JMM}{{\operatorname{JMM}}}
\newcommand{\JM}{{\operatorname{M}}}
\newcommand{\MM}{{\sM\sM}}
\newcommand{\coker}{{\operatorname{coker}}}
\newcommand{\MotGal}{{\operatorname{MotGal}}}
\newcommand{\MotU}{{\operatorname{\widehat{\sU}}}}
\newcommand{\ssim}{{\operatorname{ss}}}
\newcommand{\Hdg}{{\operatorname{Hdg}}}

\newcommand{\ess}{{\rm ess}}
\newcommand{\PSh}{{\mathop{PSh}}}
\newcommand{\Sh}{{\mathop{Sh}}}
\newcommand{\EM}{{\mathop{EM}}}
\newcommand{\SmCor}{{\mathop{SmCor}}}
\newcommand{\mot}{{\mathop{mot}}}
\newcommand{\Mot}{{\mathop{Mot}}}
\newcommand{\gm}{{\mathop{gm}}}
\newcommand{\Sing}{{\operatorname{Sing}}}
\newcommand{\equi}{{\mathop{equi}}}
\newcommand{\CM}{{\operatorname{\mathcal{CM}}}}
\newcommand{\KCM}{{\operatorname{\mathcal{KCM}}}}
\newcommand{\Conn}{{\mathop{Conn}}}
\newcommand{\PreSh}{{\mathop{PreSh}}}
\newcommand{\nh}{{\mathfrak{h}}}
\newcommand{\calg}{{\mathfrak{ca}}}

\newcommand{\colim}{\operatornamewithlimits{\varinjlim}}

\newcommand{\MAT}{{\operatorname{MAT}}}
\newcommand{\un}{\underline}
\newcommand{\nai}{{pd}}
\newcommand{\cycmod}{{\mathfrak{cm}}}
\newcommand{\efc}{{}}

\title{Tate motives and the fundamental group}
\date{August 26, 2007}
\author{H\'el\`ene Esnault}
\address{ Universit\"at Duisburg-Essen\\
Mathematik\\
45117 Essen\\ Germany
}
\email{esnault@uni-due.de}

\author{Marc Levine}

\address{ 
Department of Mathematics\\
Northeastern University\\
Boston, MA 02115\\
USA}
\email{marc@neu.edu}

\keywords{Motivic $\pi_1$, Tate motive, bar construction}

\subjclass{Primary 14C25, 19E15; Secondary 19E08 14F42, 55P42}
\thanks{The first-named author is supported by the Leibniz Preis der DFG. The second-named author gratefully acknowledges the support of the Humboldt Foundation  and support of the NSF via grant DMS-0457195}

\renewcommand{\abstractname}{Abstract}
\begin{abstract}  Let $k$ be a number field,   and let  $S\subset \P^1(k)$  be a finite set of rational points. We relate the Deligne-Goncharov contruction of the motivic  fundamental group  of $X:=\P^1\setminus S$ to the Tannaka group scheme over $\Q$ of the category of mixed Tate motives over $X$.
\end{abstract}

\maketitle
\tableofcontents
\section*{Introduction}  In \cite{Del}, P. Deligne defined  the motivic fundamental group of $X=\P^1\setminus\{0,1,\infty\}$ over a number field $k$ as an object in the  category of systems of realizations. This is a  Tannakian category over $\Q$, which he constructed as  tuples (Betti, de Rham, $\ell$-adic, crystalline), with compatibilities between them, a definition   close to the one given by U. Jannsen \cite{J}.  The Betti-de Rham  component is  the mixed Hodge structure, defined by J. Morgan \cite{Morgan}, on the nilpotent completion $\varprojlim_N \Q[\pi^{{\rm top}}(X,a)]/I^N$ of the topological fundamental group $\pi_1^{{\rm top}}(X(\C),a)$,  for all complex embeddings $k\subset \C$, where  the base-point $a$ is  either a point in $X(k)$ or a non-trivial tangent vector at $\bar{a}\in (\P^1\setminus X)(k)$.

A. Beilinson (\cite[proposition~3.4]{DelGon}) showed that for any smooth complex variety $X$,  and for base-point  $a\in X(\C)$, the ind-system 
\[
\varinjlim_N {\rm Hom}_{\Q}(\Q[\pi^{{\rm top}}(X,a)]/I^N, \Q),
\]
 which is a Hopf algebra over $\Q$,  arises from  the cohomology of a cosimplicial scheme   $\sP_a(X)$.  As pointed out by Z. Wojtkoviak \cite{Woj}, the Hopf algebra structure on   $\varinjlim_N (\Q[\pi^{{\rm top}}(X,a)]/I^N)^\vee$ similarly arises from operations on $\sP_a(X)$. 
These  key  results have many consequences.  For instance, one can use $\sP_a(X)$ to define the mixed Hodge structure on $\varprojlim_N \Q[\pi^{{\rm top}}(X,a)]/I^N$,  {\it cf.} \cite{Hain}.  Even more,  the cosimplicial scheme   $\sP_a(X)$, regardless of the geometry of $X$, defines an ind- Hopf algebra object   $h_k(\sP_a(X))$ in Voevodsky's triangulated category of motives $\DM_{gm}(k)$ \cite[chapter~V]{FSV}; here 
\[ 
h_k:\Sm/k^\op\to \DM_{gm}(k)
\]  
is the ``cohomological motive" functor, dual to the canonical functor $M_\gm:\Sm/k\to \DM_{gm}(k)$. If in addition $X$   is the complement in $\P^1_k$ of a finite set of $k$-rational points, then   $h_k(\sP_a(X))$ lies in the full triangulated subcategory $\DMT(k)$ of $\DM(k)_\Q$ spanned by the Tate objects $\Q(n)$.

As explained in \cite{LevineTate}, if $k$ satisfies the Beilinson-Soul\'e vanishing conjecture, that is,  if the motivic cohomology $H^p(k,\Z(q))$ vanishes for $p\le 0$, $q>0$, there is a $t$-structure defined on $\DMT(k)$, the heart of which is the abelian category $\MT(k)$ of mixed Tate motives over $k$.  $\MT(k)$ is a $\Q$-linear, abelian rigid tensor category with the structure of a functorial exact weight filtration $W_*$. Taking the associated graded object with respect to $W_*$ defines a neutral fiber functor $\gr^W_*$,   endowing $\MT(k)$ 
with the structure of a Tannakian category over $\Q$. 

By the work of Borel \cite{Borel}, we know that if $k$ is a number field,  then $k$ does satisfy the Beilinson-Soul\'e conjecture. Thus Beilinson's construction allows one to define  the  ind-Hopf algebra object   $H^0(h_k(\sP_a(X))$ in  $\MT(k)$,   if $k$ is a number field.  In \cite[th\'eor\`eme~4.4]{DelGon} P. Deligne and A. Goncharov show that the  dual $\pi^{{\rm mot}}_1(X,a)$ of  $H^0(h_k(\sP_a(X))$, which is a pro-group scheme object in $\MT(k)$, yields  Deligne's original motivic fundamental group upon applying the appropriate realization functors, in case $a\in X(k)$ and $X\subset \P^1_k$ is the complement of a finite set of $k$-points of $\P^1$.  In addition, they show that, even for  a tangential base-point $a$, there is a pro-group scheme object $\pi^{{\rm mot}}_1(X,a)$ in $\MT(k)$ which maps to Deligne's motivic fundamental group under realization, without, however, making an explicit construction of $\pi^{{\rm mot}}_1(X,a)$ in this case. Using this construction as starting point, they go on to construct a motivic fundamental group for any unirational variety over the number field $k$, as a pro-group scheme over the larger Tannakian category of Artin-Tate motives $\MAT(k)$ (see \cite{DelGon} for details).

Using recent work of Cisinski-D\'eglise \cite{CisinskiDeglise}, one now has available a reasonable candidate for the category of motives over a base $X$, at least if $X$ is a smooth variety over a perfect field $k$. In any case, the resulting triangulated category $\DM(X)$ has Tate objects $\Z_X(n)$ which properly compute the motivic cohomology of $X$ (defined using Voevodsky's category $\DM_{gm}(k)$). 
In addition, if $X\subset \P^1_k$ is an open defined over a number field $k$, then 
the observation made in  \cite{LevineTate} carries over to the full triangulated subcategory $\DMT(X)$ of  the category $\DM(X)_\Q$  generated by the Tate objects $\Q_X(n)$. Indeed, the localization sequence and homotopy invariance for motivic cohomology
allow one to reduce Beilinson-Soul\'e vanishing conjectures for the motivic cohomology $H^p(X,\Z(q))$ of $X$ to   the conjectures  for finite extensions of $k$. Thus, assuming $k$ is a number field,  there is a heart  $\MT(X) \subset \DMT(X)$ which is a $\Q$-linear abelian rigid tensor category, and which receives $\MT(k)$ by pull back via the structure morphism $p:X \to \Spec k$.  

By Tannaka duality, we therefore have the Tannaka group schemes over $\Q$,  $G(\MT(X),\gr^W_*)$  and $G(\MT(k),\gr^W_*)$, and $p^*:\MT(k)\to \MT(X)$ gives a canonical short exact sequence
\[
1\to K\to G(\MT(X),\gr^W_*)\xrightarrow{p_*} G(\MT(k),\gr^W_*)\to 1.
\]
$K$ is defined as the kernel of $p_*$; the surjectivity of $p_*$ follows from the fact that each  $a\in X(k)$ 
defines a splitting $s_{a*}:G(\MT(k),\gr^W_*)\to G(\MT(X),\gr^W_*)$ to $p_*$. In fact, the splitting $s_{a*}$ also defines an action of $G(\MT(k),\gr^W_*)$ on $K$, which lifts the $\Q$ group-scheme $K$ to a group-scheme object $K_a$ in $\MT(k)$.

  In \cite[section~4.19]{DelGon}, Deligne and Goncharov  use the group-scheme $\pi_1^{{\rm mot}}(X,a)$ over $\MT(k)$ to {\em define}  $\MT(X)$ as the category of $\MT(k)$ representations of $\pi_1^{{\rm mot}}(X,a)$. In \cite[section~4.22]{DelGon} they ask about the relationship between 
$\MT(k(\P^1))$, defined as above as a subcategory of Voevodsky's category $\DM(k(\P^1))_\Q$,  and $\varinjlim_{X\subset \P^1} \MT(X)$  (this is the formulation for $k=\bar\Q$, in general, one needs to use the Artin-Tate motives $\MAT$).
The purpose of this article is to give an answer to this question in the following shape: the intrinsic
definition of $\MT(X)$ we outlined above is equivalent to the category of $K_a$-representations in $\MT(k)$, assuming $\P^1\setminus X$ consists of $k$-rational points.

We now  describe our main result (corollary~\ref{cor:main}; see also theorem~\ref{thm:Main} for a more general statement). 

\begin{IntroThm}
 Let $k$ be a number field, $S\subset X(k)$ a finite set of $k$-points of $\P^1$, let $X:=\P^1\setminus S$  and take $a\in X(k)$. Then  the pro-group scheme objects $K_a$ and $\pi_1^{{\rm mot}}(X,a)$ are isomorphic   as group-schemes in $\MT(k)$. 
\end{IntroThm}
 The  equivalence of $\MT(X)$ with the category of $K_a$-representations in $\MT(k)$ follows directly from this.

We  now explain the  ideas that go into the proof.  In \cite{BlochKriz} S. Bloch and I. Kriz construct a group-scheme $G_{BK}(k)$ over $\Q$, by applying the bar construction to the {\em cycle algebra}  $\sN_k:=\Q\oplus \oplus_{r\ge1}\sN_k(r)$. The $r$th component $\sN_k(r)$ of $\sN_k$ is a shifted, alternating version of Bloch cycle complex,
\[
\sN^m_k(r)=z^r(k, 2r-m)^{{\rm Alt}};
\]
the alternation makes the product on $\sN_k$ strictly graded-commutative. The additional grading $r$ is the {\em Adams grading}. The reduced bar construction gives us the Adams graded Hopf algebra  $H^0( \bar{B}(\sN_k))$ and $G_{BK}(k)$ is the pro group scheme $\Spec H^0( \bar{B}(\sN_k))$. Bloch-Kriz define the category of ``Bloch-Kriz" mixed Tate motives over $k$, $\MT_{BK}(k)$, as the finite dimensional graded representations of $G_{BK}(k)$  in  $\Q$-vector spaces.

 In  \cite{KrizMay}, I. Kriz and P. May consider, for an Adams graded commutative differential graded algebra (cdga) $\sA=\Q\cdot\id\oplus\oplus_{r\ge 1}\sA(r)$ over $\Q$, the ``bounded" derived category $\sD^f_\sA$ of Adams graded dg $\sA$ modules. $\sD^f_\sA$ admits a  functorial exact weight filtration, arising from the Adams grading; in case $\sA$ is cohomologically connected, $\sD^f_\sA$ has a $t$-structure, defined by pulling back the usual $t$-structure on $\sD^f_\Q\cong \oplus_n D^b(\Q)$ via the functor $M\mapsto M\otimes^L_\sA\Q$ from $\sD^f_\sA$ to $\sD^f_\Q$. In particular, they define the heart $\sH^f_\sA$. Next,  assuming $\sA$ cohomologically connected,  they construct an exact functor
 \[
 \rho:D^b\big(\coRep_\Q^f(H^0( \bar{B}(\sA)))\big)\to \sD^f_\sA
 \]
 where $\coRep_\Q^f(H^0( \bar{B}(\sA)))$ is the category of graded  co-representations of $H^0( \bar{B}(\sA))$ in finite-dimension $\Q$-vector spaces, and show that $\rho$ identifies the categories $\sH^f_\sA$ and
 $\coRep_\Q^f(H^0( \bar{B}(\sA)))$ (although $\rho$ is not in general an equivalence).  For those who prefer group-schemes to Hopf algebras, let $G_\sA:=\Spec H^0( \bar{B}(\sA))$. Then $G_\sA$ is a pro-affine group scheme over $\Q$ with $\G_m$ action, and $\coRep_\Q^f(H^0( \bar{B}(\sA)))$ is equivalent to the category of graded representions of $G_\sA$ in finite dimensional $\Q$-vector spaces.

 Taking $\sA=\sN_k$, and noting that the Beilinson-Soul\'e vanishing conjectures for $k$ are equivalent to the cohomological connectedness of $\sA$, this gives an equivalence of the heart $\sH^f_{\sN_k}$ with the Bloch-Kriz mixed Tate motives $\MT_{BK}(k)$.

 M. Spitzweck \cite{Spitzweck} (see \cite[section~5]{LevineHandbook} for a detailed account)  defines an equivalence
 \[
\theta_k: \sD^f_{\sN_k}\to \DMT(k) \subset \DM_{gm}(k)_\Q
\]
for $k$ an arbitrary field. In addition, under the assumption that $k$ satisfies the Beilinson-Soul\'e conjectures, or equivalently, that $\sN_k$ is cohomologically connected, $\theta_k$ restricts to an equivalence 
\[
\theta_k:\sH^f_{\sN_k}\to \MT(k).
\]
From the discussion above, this gives an equivalence of $\coRep_\Q^f(H^0( \bar{B}(\sN_k)))$ with $\MT(k)$, and in fact identifies $G_{BK}(k)\ltimes\G_m$ as the Tannaka group of $(\MT(k),\gr^*_W)$.

Our first task is  to extend this picture from $k$ to $X$. To this aim, one defines the cycle algebra $\sN_X$ by replacing $k$ with $X$ in the definition of $\sN_k$ and modifying the construction further by using  complexes of cycles which are equi-dimensional over $X$. This yields  an Adams graded cdga over $\Q$ together with a map of Adams graded cdgas $p^*:\sN_k\to \sN_X$ arising from the structure morphism $p:X\to \Spec k$. Essentially the same construction as for $k$ gives an equivalence
\begin{equation}\label{eqn:Spitzweck}
\theta_X:\sD^f_{\sN_X}\to \DMT(X)\subset \DM(X)_\Q \tag{$*$}
\end{equation}
and if $X$ satisfies the Beilinson-Soul\'e vanishing conjectures, $\theta_X$ restricts to an equivalence $\sH^f_{\sN_X}\sim \MT(X)$. Defining the $\Q$ pro-group scheme $G_{BK}(X)$ as above, 
\[
G_{BK}(X):= G_{\sN_X}=\Spec(H^0( \bar{B}(\sN_X))),
\]
we also have the equivalence of $\MT(X)$ with the graded representations of $G_{BK}(X)$ in finite dimensional $\Q$-vector spaces, giving the identification of $G_{BK}(X)\ltimes\G_m$ with the Tannaka group of $(\MT(X),\gr^W_*)$, and identifying $p_*:G(\MT(X),\gr^W_*)\to G(\MT(k),\gr^W_*)$ with the map
\[
\tilde{p}\times\id:G_{BK}(X)\ltimes\G_m\to G_{BK}(k)\ltimes\G_m,
\]
with $\tilde{p}$ induced from $p^*:\sN_k\to \sN_X$. 

A $k$-point $a$ of $X$ gives an augmentation $\epsilon_a:\sN_X\to \sN_k$.  We discuss the general theory of augmented cdgas in section~\ref{sec:RelDga},  leading to  the  {\em relative bar construction}  $H^0_\sN(\bar{B}_\sN(\sA,\epsilon))$ of a cdg $\sN$ algebra $\sA$ with augmentation $\epsilon:\sA\to \sN$, as an ind-Hopf algebra in $\sH^f_\sN$. Let $G_{\sA/\sN}(\epsilon)=\Spec H^0_\sN(\bar{B}_\sN(\sA,\epsilon))$
and let $G_{\sA/\sN}(\epsilon)_\Q$ be the pro-group scheme over $\Q$ gotten from $G_{\sA/\sN}(\epsilon)$ by applying the fiber functor $\gr^W_*:\sH^f_\sN\to \Vect_\Q$. Note that Tannaka duality gives a canonical action of $G_\sN$ on $G_{\sA/\sN}(\epsilon)_\Q$.

Of course, in order to make a reasonable relative bar construction, one needs to use a good model for $\sA$ as an $\sN$-algebra. This is provided by using the {\em relative minimal model} $\sA\{\infty\}_\sN$ of $\sA$ over $\sN$, for which the derived tensor product is just the usual tensor product. 

In section~\ref{subsec:SemiDirect}, especially  theorem~\ref{thm:Kernel}, we show that
\begin{enumerate}
\item $G_{\sA/\sN}(\epsilon)_\Q=\Spec H^0(\bar{B}(\sA\{\infty\}_\sN\otimes_\sN\Q))$.
\item There is an exact sequence of pro-group schemes over $\Q$:
\[
1\to G_{\sA/\sN}(\epsilon)_\Q\to G_\sA\xrightarrow{p_*} G_\sN\to 1
\]
The splitting $\epsilon^*$ to $p^*$ defines a splitting $\epsilon_*:G_\sN\to G_\sA$ to $p_*$.
\item The conjugation action of $G_\sN$ on  $G_{\sA/\sN}(\epsilon)_\Q$ given by the splitting $\epsilon_*$ is the same as the canonical action.
\end{enumerate}

To do this, we use an alternate description of dg modules over an Adams graded cdga $\sN$, that of {\em flat dg connections}.  Kriz and May describe dg modules $M$ over $\sN$ as $\sN^{+}:=\bigoplus_{r>0} \sN(r)$-valued  connections over $M\otimes_\sN \Q$ (for the canonical augmentation $\sN\to \Q$). Writing $\sA\{\infty\}_\sN^{+}$ as $\sN^{+} \oplus \sI$,  with this decomposition coming from the augmentation $\sA\{\infty\}_\sN\to \sN$, the absolute (i.e. $\sA\{\infty\}_\sN^{+}$-valued) 
connection on $H^0(\bar{B}(\sA))=H^0(\bar{B}(\sA\{\infty\}_\sN))$ induces a $\sN^{+}$-valued  connection on
$H^0(\bar{B}(\sA\{\infty\}_\sN\otimes_{\sN}\Q))$. Similarly, the structure of $H^0_\sN(\bar{B}_\sN(\sA,\epsilon))$ as an ind-Hopf algebra in $\sH^f_\sN$ gives an $\sN^{+}$-valued  connection on 
\[
H^0_\sN(\bar{B}_\sN(\sA,\epsilon))\otimes_\sN\Q=H^0(\bar{B}(\sA\{\infty\}_\sN\otimes_{\sN}\Q)).
\]
Using this description, it is easy to make the identifications necessary for proving (1)-(3) above.
In a  metaphorical way, we could say that  $G_{\sA/\sN}(\epsilon)$ is the Gau{\ss}-Manin connection of $G_{\sA}$ associated to $\sA/\sN$.

Applying this theory to the splitting $\epsilon_a:\sN_X\to \sN_k$, the $\Q$ pro-group scheme $K$, and lifting $K_a$ to a $\MT(k)$ pro-group scheme, gives us the isomorphism of pro-group schemes 
\[
K\cong \Spec H^0(\bar{B}(\sN_X\{\infty\}_{\sN_k}\otimes_{\sN_k}\Q))
\]
 and the isomorphism of  pro-group scheme objects in $\sH^f_{\sN_k}$
\begin{equation}\label{eqn:KernIdent}
K_a\cong H^0_\sN(\bar{B}_{\sN_k}(\sN_X,\epsilon_a)).\tag{$**$}
\end{equation}

One can make  the dg $\sN_k$-module $H^0_\sN(\bar{B}_{\sN_k}(\sN_X,\epsilon_a))$ explicit  as an object in $\MT(k)$ via Spitzweck's theorem.  This relies on a crucial property of the transformation from dg  $\sN_k$ modules to motives:\\
\\
Take $X\in \Sm/k$. If the motive $h_k(X)$ is in $\DTM(k)$ and $X$ satisfies the Beilinson-Soul\'e vanishing conjectures, then the motive
of the cycle module  $\sN_X\{\infty\}_{\sN_k}$ is canonically isomorphic to $h_k(X)_\Q$. \\
\\
The explicit decription of the Beilinson simplicial scheme  underlying the Deligne-Goncharov construction, together with this essential fact, allows one to conclude that the Gau{\ss}-Manin connection is precisely  $\pi_1^{\rm \mot}(X,a)$, when $a$ comes from a rational point $a\in X(k)$ (see sections~\ref{subsec:CompThm} and \ref{subsec:FundSeq}). In other words, we have the isomorphism of pro-group schemes over $\MT(k)$:
\[
\pi_1^{\rm \mot}(X,a)\cong \Spec  H^0_\sN(\bar{B}_{\sN_k}(\sN_X,\epsilon_a)).
\]
Combining this with our identification \eqref{eqn:KernIdent}   proves the main theorem.

 The generalization of Spitzweck's theorem, i.e.,   the identification \eqref{eqn:Spitzweck} of the triangulated Tate subcategory  $\DMT(X)$ of  $\DM(X)_\Q$  with the bounded derived category of dg modules over $\sN_X$,  is entirely due to the second author. This result together with a sketch of the proof is summarized in section~\ref{sec:NModMot}. Along with the necessary constructions on motives over a base and cycle algebras (sketched here in sections~\ref{sec:Motives} and \ref{sec:CycAlg}) this material will be developed to a greater extent in a forthcoming article \cite{LevineTateMotive} by the second author.

We now discuss a few questions. In this article, 
we do not consider the case  of the base-point $a$ being  a non-trivial tangent vector at some point $\bar{a}\in \P^1\setminus X$. As mentioned above,  Deligne-Goncharov \cite{DelGon} show in this case as well that  the motivic $\pi_1$, defined by Deligne  \cite{Del} as a system of realizations,  comes from $\MT(k)$. This defines $\pi_1^{{\rm mot}}(X,a)$ as an object in $\MT(k)$, but there is no direct construction in $\MT(k)$. 
However, the results of \cite{LevineTubular} give a section $\epsilon_a$ to $p_*:\sN_k\to \sN_X$ (in the homotopy category of cdgas) for tangential base-points $a$ as well as for $k$-points, so we do have a relative bar construction available even for tangential base-points.  In order to extend our main theorem~\ref{thm:Main}  to this case, one should define realization functors on the categories of Tate motives, described as  dg modules over the cycle algebra, and check that the realization of $\Spec H^0_\sN(\bar{B}_{\sN_k}(\sN_X,\epsilon_a))$ agrees with Deligne's motivic $\pi_1$.

At any rate this raises another question.  Given {\em any} section $s: G_{\sN_k}\to G_{\sN_X}$, the resulting action of $G(\MT(k), {\rm gr}^W_*)$ on $K$ defines a pro-group scheme object $K_s$   over  $\MT(k)$. We have as well a descripton of the action in terms of the Gau{\ss}-Manin connection.  We could ask for a description of  all the sections and a condition which separates out the ones which are geometric, i.e., arising from a rational point or tangential base-point.

Recall that sections  $s^{\acute{e}t}$ of the natural homomorphism $\pi_1^{\acute{e}t}(X) \to {\rm Gal}(\bar k/k)$ of Grothendieck's pro-finite fundamental group to the Galois  group of $k$ 
are predicted by Grothendieck's section conjecture to be geometric (see \cite{GroFa}), if $|(\P^1\setminus X)(\bar k)|\ge 3$, (or more generally if $X$ is a hyperbolic curve),  and $k$ is a number field. Geometric in this setting means that there is a rational point $a\in \P^1$ such that the restriction of the section to $\pi^{\acute{e}t}_1(X\cup \{a\})$ comes from $a\in \P^1(k)$. This implies in particular  that each section
$s^{\acute{e}t}$ lifts to a section of   $\pi_1^{\acute{e}t}(U)\to  {\rm Gal}(\bar k/k)$ for all $U\subset X$ open. If one knew this, 
then the conjecture would reduce to the case $X=\P^1\setminus \{0,1,\infty\}$ (see \cite[proposition~7.9]{EH}). In addition, it is shown in \cite[theorem~6.5]{EH} that the geometric sections are described precisely by $\varprojlim_n k^\times/(k^\times)^n$ at each rational point at infinity. This group contains $k^\times$, this subgroup corresponding to Deligne's tangential base points (\cite[section~6]{EH}). On the other hand, sections $s^{\acute{e}t}$ should define 
sections $s_\ell$ of $G(\MT(X), \gr^W_*)\otimes \Q_\ell\to G(\MT(k),\gr^W_*)\otimes \Q_\ell$ for all $\ell$.

It is thus natural to ask about the sections of 
\[
p_*:G(\MT(X),\gr^W_*)\to G(\MT(k),\gr^W_*),
\]
 i..e. graded sections of $G_{BK}(X)\xrightarrow{\tilde{p}} G_{BK}(k)$. It is rather easy to see that the Lie algebra of $G_{BK}(k)$ is a free pro-nilpotent Lie algebra over $\Q$, so there is essentially no restriction on the graded sections to $\tilde{p}$. Thus, we need to look deeper for an analog to Grothendieck's conjecture.  
 
 If we fix a $k$-point $a\in X(k)$, we have  the conjugation representation $\rho_a$ of $G_{BK}(k)$ on the kernel $K$.  If we want to pick out the geometric sections, we can first ask:   suppose we have a  section $s$ to $p_*$ with the property that the resulting conjugation representation $\rho_s$ is isomorphic to $\rho_a$. Is $s$ then geometric, in fact, is $s=s_a$? 
 
 If the answer is yes, we can then ask what properties make the representations $\rho_a$ special. An analog of the Deligne-Ihara conjecture suggests the following: Note that the Lie algebra ${\rm Lie}(G_{BK}(k))$ is a  graded pro-nilpotent Lie algebra over $\Q$, concentrated in negative degrees. Let $\sI\subset {\rm Lie}(G_{BK}(k))$ be the ideal generated by the degree $-1$ homogeneous elements. We note that the space of degree $-1$ elements in  ${\rm Lie}(G_{BK}(k))$ is the pro-vector space dual to $H^1(k,\Q(1))\cong k^\times_\Q$.
 
 \begin{IntroConj} Let $a$ be the tangential base point $(0,\partial/\partial t_{|0})$ for $X:=\P^1_\Q\setminus\{0,1,\infty\}$. Then the map
 \[
 d\rho_a:{\rm Lie}(G_{BK}( \Q))\to \End({\rm Lie}(K))
 \]
  has kernel equal to  $\sI$.
  \end{IntroConj}
It is not difficult to show that the kernel of $d\rho_a$ contains $\sI$.
\\ \ \\
{\it Acknowledgements:} We gave a seminar in the winter 2006-7 at the university of Duisburg-Essen on \cite{DelGon},  to try to understand the constructions and results of Deligne-Goncharov, as well as the various constructions of mixed Tate motives and the relationships between them, as developed in the works of Bloch, Bloch-Kriz, Kriz-May and Spitzweck,  and summarized in \cite{LevineHandbook};  our paper  is to a large extent a product of that seminar. We thank all the seminar participants  for their willingness to give talks, in particular we thank Ph\`ung H\^o Hai for various discussions on Tannakian categories.

\section{Differential graded algebras}\label{sec:Dga} We fix notation  and recall some basic facts on commutative differential graded algebras (cdgas) over $\Q$. This material is taken from  \cite{KrizMay}. 

In what follows a cdga will always mean a cdga over $\Q$.

\subsection{Adams graded cdgas}
\begin{definition} \index{cdga} (1)  A {\em cdga} $(A^*, d,\cdot)$ (over $\Q$) consists of a unital, graded-commutative $\Q$-algebra $(A^*:=\oplus_{n\in \Z}A^n,\cdot)$ together with a graded homomorphism $d=\oplus_nd^n$, $d^n:A^n\to A^{n+1}$, such that
\begin{enumerate}
\item $d^{n+1}\circ d^n=0$.
\item $d^{n+m}(a\cdot b)=d^n a\cdot b+(-1)^na\cdot d^mb$; $a\in A^n$, $b\in A^m$.
\end{enumerate}
$A^*$ is called {\em connected} if $A^n=0$ for $n<0$ and $A^0=\Q\cdot 1$, {\em cohomologically connected} if $H^n(A^*)=0$ for $n<0$ and $H^0(A^*)=\Q\cdot 1$.

\medskip
\noindent
(2) An {\em Adams graded} cdga  is a cdga $A$ together with a direct sum decomposition into subcomplexes $A^*:=\oplus_{r\ge0}A^*(r)$ such that $A^*(r)\cdot A^*(s)\subset A^*(r+s)$. In addition, we require that $A^*(0)=\Q\cdot \id$.
An Adams graded cdga is  said to be (cohomologically) connected if the underlying cdga is (cohomologically) connected. 

For $x\in A^n(r)$, we call $n$ the {\em cohomological degree} of $x$, $n:=\deg x$, and $r$ the {\em Adams degree} of $x$, $r:=|x|$.
\end{definition}
Note that an Adams graded cdga $A$ has a canonical augmentation  $A\to\Q$ with augmentation ideal $A^+:=\oplus_{r>0}A^*(r)$.

\subsection{The bar construction}\label{subsection:bar}  We let $\Ord$ denote the category with objects the   sets $[n]:=\{0,\ldots, n\}$, $n=0, 1,\ldots$, and morphisms the non-decreasing maps of sets. The morphisms in $\Ord$ are generated by the {\em coface maps} $\delta^n_i:[n]\to [n+1]$ and the {\em codegeneracy maps} $\sigma^n_i:[n]\to [n-1]$, where $\delta^n_i$ is the strictly increasing map omitting $i$ from its image and $\sigma^n_i$ is the non-decreasing surjective map sending $i$ and $i+1$ to $i$. For a category $\sC$, we have the categories of {\em cosimplicial objects} in $\sC$ and  {\em simplicial objects} in $\sC$, namely, the categories of functors $\Ord\to \sC$ and $\Ord^\op\to \sC$, respectively. For a cosimplicial object $X:\Ord\to \sC$, we often write $\delta^n_i$ and $\sigma^n_i$ for the coface maps $X(\delta^n_i)$ and $X(\sigma^n_i)$, and for a simplicial object $S:\Ord^\op\to \sC$, we often write $d^n_i$ and $s^n_i$ for the {\em face} and {\em degeneracy} maps $S(\delta^n_i)$ and $S(\sigma^n_i)$.

  Let $A$ be a cdga. We begin by defining the simplicial cdga   $B_\bullet(A)$ as follows: Tensor product (over $\Q$) is the coproduct in the category of cdgas, so for a finite set $S$, we have $A^{\otimes S}$, giving the functor $A^{\otimes ?}$ from finite sets to cdgas. Concretely, 
for 
\[
\varphi: S:=\{i_1,\ldots, i_s\}  \to T:= \{j_1,\ldots, j_t\}
\]
 a  map of finite sets, the induced map   $A^{\otimes \varphi}:A^{\otimes S}\to A^{\otimes T}$ of cdgas is defined by 
 \begin{align*}
 &a_{i_1}\otimes \ldots \otimes a_{i_s} \mapsto b_{j_1}\otimes \ldots \otimes b_{j_t} \\
& b_j= \prod_{\varphi(i)=j}a_i;\ 
b_j=1\text{ if  }\varphi^{-1}(j)=\0.
\end{align*} 
Thus, if we have a simplicial set $S$ such that $S[n]$ is a finite set for all $n$, we may form the simplicial cdga $A^{\otimes S}$, $n\mapsto A^{\otimes S[n]}$.  We have the representable simplicial sets $\Delta[n]:=\Hom_\Ord(-,[n])$; setting $[0,1]:=\Delta[1]$ gives us the simplicial cdga
\[
B_\bullet(A):=A^{\otimes [0,1]}.
\]

The two inclusion $[0]\to [1]$ define the maps $i_0, i_1:\Delta[0]\to \Delta[1]$. Letting $\{0,1\}$ denote the constant simplicial set with two elements, the maps $i_0, i_1$ give rise to  the map of simplicial sets $i_0\amalg i_1:\{0,1\}\to [0,1]$, which makes $B_\bullet(A)$ into a simplicial $A\otimes A=A^{\otimes\{0,1\}}$ algebra. 

Suppose we have augmentations $\epsilon_1,\epsilon_2:A\to\Q$. Define $\bar{B}_\bullet(A,\epsilon_1, \epsilon_2)$ by
\[
\bar{B}_\bullet(A,\epsilon_1, \epsilon_2):=B_\bullet(A)_{A\otimes A}\Q
\]
using $\epsilon_1\otimes\epsilon_2:A\otimes A\to \Q$ as structure map. Since $\bar{B}_n(A,\epsilon_1, \epsilon_2)$ is a complex for each $n$, we can form a double complex by using the usual alternating sum of the face maps $d^n_i:\bar{B}_{n+1}(A,\epsilon_1, \epsilon_2)\to \bar{B}_n(A,\epsilon_1, \epsilon_2)$ as the second differential, and let $\bar{B}(A,\epsilon_1, \epsilon_2)$ denote the total complex of this double complex. We use cohomological grading throughout, so $\bar{B}_n(A,\epsilon_1, \epsilon_2)^m$ has total degree $m-n$.  For $\epsilon_1=\epsilon_2=\epsilon$, we write $\bar{B}(A,\epsilon)$ or simply $\bar{B}(A)$; this is the {\em reduced bar construction} for $(A,\epsilon)$.  As is usual, we denote a decomposable element $x_1\otimes\ldots\otimes x_n$ of $\bar{B}(A)$ by $[x_1|,\ldots|x_n]$. Note that
\[
\deg([x_1|\ldots|x_m])=-m+\sum_i\deg(x_i).
\]

The bar construction $\bar B:=\bar B(A)$ has the following structures: a differential $d:\bar{B}\to \bar{B}$ of degree +1  coming from the differential in $A$, a  product (the shuffle product)
\begin{align*}
&\cup:\bar{B}\otimes\bar{B}\to\bar{B}\\
&[x_1|\ldots|x_p]\cup[x_{p+1}|\ldots|x_{p+q}]=\sum_\sigma \sgn(\sigma)[x_{\sigma(1)}|\ldots|x_{\sigma(p+q)}]
\end{align*}
where the sum is over all $(p,q)$ shuffles $\sigma\in S_{p+q}$,  a co-product 
\begin{align*}
&\delta:\bar{B}\to\bar{B}\otimes\bar{B}\\
&\delta([x_1|\ldots|x_n]):=\sum_{i=0}^n[x_1|\ldots x_i]\otimes[x_{i+1}|\ldots|x_n]
\end{align*}
and an involution 
\begin{align*}
&\iota:\bar B\to \bar B, \\
&\iota([x_1|x_2|\ldots|x_{n-1}|x_n] ):=(-1)^n[x_n|x_{n-1}|\ldots|x_2|x_1]
\end{align*}
making   $(\bar{B}(A),d,\cup,\delta, \iota)$ a differential graded  Hopf algebra over $\Q$, which is graded-commutative with respect to the product $\cup$. The cohomology $H^*(\bar{B}(A))$ is thus a graded Hopf algebra over  $\Q$, in particular, $H^0(\bar{B}(A))$ is a  commutative Hopf algebra over $\Q$. 

Let $\sI(A)$ be the kernel of the augmentation $H^0(\bar{B}(A))\to \Q$  induced by $\epsilon$. The coproduct $\delta$ on $H^0(\bar{B}(A))$ induces the structure of a co-Lie algebra on $\gamma_A:=\sI(A)/\sI(A)^2$.

 From the formula for the coproduct, we see that,  modulo tensors of  degree $<m$, we have 
\[
\delta([x_1|\ldots|x_m])=1\otimes [x_1|\ldots|x_m]+[x_1|\ldots|x_m]\otimes 1.
\]
This implies that the pro-affine $\Q$-algebraic group $G:=\Spec H^0(\bar{B}(A))$ is pro-unipotent. In addition, it is known that, as an algebra over $\Q$, $H^0(\bar{B}(A))$ is a polynomial algebra, with indecomposables  $\gamma_A$.

Suppose $A=\oplus_{r\ge0}A^*(r)$ is an Adams graded cdga, with canonical augmentation $\epsilon:A\to \Q$. The Adams grading on $A$ induces an Adams grading on $B_\bullet(A)$ and thus on $\bar{B}(A)$; explicitly $\bar{B}(A)$ has the Adams grading $\bar{B}(A)=\oplus_{r\ge0}\bar{B}(A)(r)$ where the Adams degree of $[x_1|\ldots|x_m]$ is
\[
|[x_1|\ldots|x_m]|:=\sum_j|x_j|.
\]
Thus $H^0(\bar{B}(A))=\oplus_{r\ge0}H^0(\bar{B}(A)(r))$ becomes an Adams graded Hopf algebra over $\Q$, commutative as a $\Q$-algebra. We also have the Adams graded co-Lie algebra $\gamma_A=\oplus_{r>0}\gamma_A(r)$. 

\begin{remark} Let $A$ be an Adams graded cdga. The Adams grading  equips the  pro-unipotent affine $\Q$ group scheme  $G:=\Spec H^0(\bar{B}(A))$   with a grading, or, equivalently, with a $\G_m$-action. 
Thus $\gamma_A$ is a graded  nilpotent  co-Lie algebra, and there is an equivalence of categories between the graded co-representations of $H^0(\bar{B}(A))$ in finite dimensional graded $\Q$-vector spaces, $\coRep^f_\Q(H^0(\bar{B}(A)))$, and the graded co-representations of $\gamma_A$ in  finite dimensional graded $\Q$-vector spaces,
$\coRep^f_\Q(\gamma_A)$. 
\end{remark}

\subsection{The category of cell modules} Kriz and May \cite{KrizMay}  define a triangulated category directly from an Adams graded cdga $A$ without passing to the bar construction or using a co-Lie algebra. We recall some of their work here.

Let $A^*$ be a graded algebra over  $\Q$.   We let $A[n]$ be the left $A^*$-module which is $A^{m+n}$ in degree $m$, with the $A^*$-action given by left multiplication. If  $A^*(*)=\oplus_{n\ge 0,r\ge 0 }A^n(r)$ is a bi-graded  $\Q$-algebra, we let $A\<r\>[n]$ be the left $A^*(*)$-module which is $A^{m+n}(r+s)$ in bi-degree $(m,s)$, with action given by left multiplication.

\begin{definition}\label{def:DgMod} Let $A$ be a cdga.\\
\\
(1) A  {\em dg $A$-module} $(M^*,d)$ consists of a complex $M^*=\oplus_nM^n$  of $\Q$-vector spaces with differential $d$, together with a graded, degree zero map  $A^*\otimes_\Q M^*\to M^*, \ a\otimes m \mapsto a\cdot m$, which makes $M^*$ into a graded $A^*$-module, and satisfies the Leibniz rule
\[
d(a\cdot m)=da\cdot m+(-1)^{\deg a}a\cdot dm;\ a\in A^*, m\in M^*.
\]
\\
(2) If $A=\oplus_{r\ge0}A^*(r)$ is an Adams graded cdga, an {\em Adams graded dg $A$-module} is a dg $A$-module $M^*$ together with a decomposition into subcomplexes $M^*=\oplus_sM^*(s)$ such that $A^*(r)\cdot M^*(s)\subset M^*(r+s)$.   We say $x\in M^*$ has {\em Adams degree $s$} if $x\in M^*(s)$, and write this as $|x|=s$. \\
\\
(3) An Adams graded dg $A$-module $M$ is a {\em cell module} if 
\begin{enumerate}
\item[(a)] $M$ is free as a bi-graded $A$-module, where we forget the differential structure. That is, there is a set $J$ and   elements $b_j\in M^{n_j}(r_j)$, $j\in J$, such that the maps $a\mapsto a\cdot b_j$ induces an isomorphism of bi-graded $A$-modules
\[
\oplus_{j\in J}A\<-r_j\>[-n_j]\to M.
\]
\item[(b)] There is a filtration on the index set $J$:
\[
J_{-1}=\0\subset J_0\subset J_1\subset \ldots J_n\subset \ldots  \subset J
\]
such that $J=\cup_{n=0}^\infty J_n$ and for $j\in J_n$,  
\[
db_j=\sum_{i\in J_{n-1}}a_{ij}b_i.
\]
\end{enumerate}
A {\em finite cell module} is a cell module  with index set $J$ finite.\\
\\
We denote the category of dg $A$-modules by $\sM_A$, the $A$-cell modules by $\CM_A$ and the finite cell modules by $\CM_A^f$.
\end{definition}

\subsection{The derived category}  Let $A$ be an Adams graded cdga and let $M$ and $N$ be Adams graded dg $A$-modules. Let $\sHom_A(M,N)$ be the Adams graded dg $A$-module with
$\sHom(M,N)^n(r)$ the $A$-module   consisting of maps $f:M\to N$ with $f(M^a(s))\subset N^{a+n}(s+r)$, and differential $d$ defined by $df(m)=d(f(m))+(-1)^{n+1}f(dm)$ for $f\in \sHom(M,N)^n(r)$. Similarly, let $M\otimes_A N$ be the Adams graded dg $A$-module with underlying module $M\otimes_AN$ and 
with differential $d(m\otimes n)=dm\otimes n+(-1)^{\deg m}m\otimes dn$. 

For $f:M\to N$ a morphism of Adams graded dg $A$-modules, we let $\Cone(f)$ be the Adams graded dg $A$-module with
\[
\Cone(f)^n(r):=N^n(r)\oplus M^{n+1}(r)
\]
and differential $d(n,m)=(dn+f(m),-dm)$. Let $M[1]$ be the Adams graded dg $A$-module with $M[1]^n(r):=M^{n+1}(r)$ and differential $-d$, where $d$ is the differential of $M$. A sequence of the form
\[
M\xrightarrow{f}N\xrightarrow{i}\Cone(f)\xrightarrow{j}M[1]
\]
where $i$ and $j$ are the evident inclusion and projection, is called a {\em cone sequence}.

\begin{definition} \label{defn:Homoto_Cat} Let $A$ be an Adams graded cdga over $\Q$. We let $\sM_A$ denote the category of Adams graded dg $A$-modules, $\sK_A$ the {\em  homotopy category}, i.e.  the objects of $\sK_A$ are the objects of $\sM_A$ and
\[
\Hom_{\sK_A}(M,N)=H^0(\sHom_A(M,N)(0)).
\]
\end{definition}

The category $\sK_A$ is a triangulated category, with distinguished triangles as usual those triangles which are isomorphic   in $\sK_A$ to a cone sequence.

\begin{definition} The  {\em derived category}  $\sD_A$  of dg $A$-modules  is the localization of $\sK_A$ with respect to morphisms $M\to N$ which are quasi-isomorphisms on the underlying complexes of $\Q$-vector spaces.  For $M$ in $\sD_A$, we denote the $n$th cohomology of $M$, as a complex of $\Q$-vector spaces, by $H^n(M)$.
\end{definition}

We define the homotopy category of $A$-cell modules, resp. finite cell modules,   as the full subcategory of $\sK_A$ with objects in $\sC\sM_A$, resp. in $\sC\sM_A^f$,
\[
\KCM^f_A\subset \KCM_A\subset \sK_A.
\]
 
Note that for $A=\Q$, $\sM_\Q$ is just the category of complexes of graded $\Q$-vector spaces, and $\sD_\Q$ is the unbounded derived category of graded $\Q$-vector spaces.

\begin{proposition}[\hbox{\cite[construction~2.7]{KrizMay}}] \label{prop:Whitehead} The evident functor
\[
\KCM_A\to  \sD_A
\]
is an  equivalence of triangulated   categories.
\end{proposition}

We let $\sD_A^f\subset\sD_A$ be the full subcategory  with objects those $M$ isomorphic in $\sD_A$ to a finite cell module. As an immediate consequence of proposition~\ref{prop:Whitehead}, we have

\begin{proposition}\label{prop:WhiteheadFin} $\KCM^f_A\to \sD_A^f$ is an equivalence of triangulated categories.
\end{proposition}

\begin{example}[Tate objects]   For $n\in\Z$, let $\Q(n)$ be the object of $\CM^f_A$ which is the free rank one $A$-module with generator $b_n$ having Adams degree $-n$, cohomological degree 0 and $db_n=0$,  i.e., $\Q(n)=A\<n\>$. We sometimes write $\Q_A(n)$ for $\Q(n)$; $\Q(n)$ is called a {\em Tate object}.   
\end{example}

\subsection{Weight filtration}  Let $M$ be an Adams graded dg $A$-module which is free as a bi-graded $A$-module. Choose a basis $\sB:=\{b_j\ |\ j\in J\}$, $M=\oplus_jA\cdot b_j$. Write
\[
db_j=\sum_ia_{ij}b_i;\ a_{ij}\in A.
\]
Since $|a_{ij}|\ge0$ and $d$ has Adams degree 0, it follows that
\[
|b_i|\le |b_j|\text{ if } a_{ij}\neq0.
\]
Thus, we have the subcomplex  
\[
W^\sB_nM=\oplus_{\substack{\{j, \ |b_j|\le n}\}}A\cdot b_j
\]
of $M$.

The subcomplex $W^\sB_nM$ is independent of the choice of basis: if $\sB'=\{b'_j\}$ is another basis and if $|b'_j|=n$, then as $b'_j=\sum_ie_{ij}b_i$ and $|e_{ij}|\ge0$, it follows that $b'_j\in W_n^\sB M$ and hence $W_n^{\sB'}M\subset W_n^\sB M$. By symmetry,  $W_n^\sB M\subset W_n^{\sB'}M$. We may thus write $W_nM$ for $W_n^\sB M$.

This gives us the increasing filtration as an Adams graded dg $A$-module  
\[
W_*M: \ \ \ldots\subset W_nM\subset W_{n+1}M\subset \ldots \subset M
\]
with $M=\cup_nW_nM$.

Similarly,  for $n\ge n'$, define $W_{n/n'}M$ as the cokernel of the inclusion
$W_{n'}M\to W_nM$, i.e., $W_{n/n'}M$ is the  Adams graded dg $A$-module with basis the $b_j$ having $n'<|b_j|\le n$ and with differential induced by the differential in $W_nM$. We write $\gr^W_n$ for $W_{n/n-1}$ and $W^{>n}$ for $W_{\infty/n}$.

It is not hard to see that $W_nM$ is functorial in $M$. In particular, if $f:M\to M'$ is a homotopy equivalence of cell modules with homotopy inverse $g:M'\to M$, then $f$ and $g$ restricted to $W_nM$ and $W_nM'$ give inverse homotopy equivalences $W_nf:W_nM\to W_nM'$, $W_ng:W_nM'\to W_nM$. Thus the $W$ filtration  in $\sC\sM_A$ defines a functorial tower of 
endo-functors on $\sK\sC\sM_A$:
\begin{equation}\label{eqn:WeightFilt}
\ldots \to W_n\to W_{n+1}\to \ldots\to \id
\end{equation}

\begin{lemma} 1. The endo-functor $W_n$ is exact for each $n$.\\
\\
2. For $n'\le n\le \infty$, the sequence of endo-functors $W_{n'}\to W_n\to W_{n/n'}$ canonically extends to a distinguished triangle of endo-functors.
\end{lemma}

\begin{proof} For (1), it follows directly from the definition that $W_n$ transforms a cone sequence into a cone sequence. For (2), take $M\in \sC\sM_A$. The sequence 
\[
0\to W_{n'}M\to W_nM\to W_{n/n'}M\to 0
\]
is  split exact as a sequence of bi-graded $A$-modules. Thus (2) follows from the general fact that a sequence in $\sC\sM_A$
\[
0\to N'\xrightarrow{i} N\xrightarrow{p} N''\to 0
\]
that is  split exact as a sequence of bi-graded $A$-modules extends canonically to a distinguished triangle in $\sK\sC\sM_A$. To see this,  choose a splitting $s$ to $p$ (as bi-graded $A$-modules), and define $t:N''\to N'[1]$ by $i\circ t=s\circ d_{N''}-d_N\circ s$. It is then easy to check that $t$ is a map of complexes and $(s,t):N''\to N\oplus N'[1]$ defines the map of complexes
\[
(s,t):N'\to \cone(i)
\]
making the diagram
\[
\xymatrix{
 N'\ar[r]^{i}\ar@{=}[d]& N\ar[r]^{p}\ar@{=}[d]& N''\ar[r]^t\ar[d]^{(s,t)}&N'[1]\ar@{=}[d]\\
  N'\ar[r]^{i}& N\ar[r]&\cone(i)\ar[r] &N'[1]}
  \]
  commute. In particular, $(s,t)$ is an isomorphism in $\sK\sC\sM_A$. One sees similarly that another choice $s'$ of splitting   leads to a homotopic map $(s',t')$.
\end{proof}

Note that it is not necessary for $M$ to be a cell module to define $W_nM$; being free as a bi-graded $A$-module suffices. It is not however clear that $W_nM$ is a quasi-isomorphism invariant in general. To side-step this issue, we use instead

\begin{definition} \label{defn:weight} Define the tower of exact endo-functors on $\sD_A$
\[
\ldots \to W_n\to W_{n+1}\to \ldots\to \id
\]
using \eqref{eqn:WeightFilt}  and the equivalence of categories in proposition~\ref{prop:Whitehead}. Concretely, choose for each $M$ in $\sD_A$ an object $P_M$ in   $\sK\sC\sM_A$ and an isomorphism $\psi_M:P_M\to M$ in $\sD_A$. Define
$W_nM$ to be the image of $W_nP_M$ in $\sD_A$, and the map $W_nM\to M$ by composing $W_nP_M\to P$ with the chosen isomorphism $\psi_M$. $W_n$ is defined on morphisms via the isomorphism
\[
\Hom_{\sK\sC\sM_A}(P_M,P_N)\cong
\Hom_{\sD_A}(P_M,P_N)\xrightarrow{(\psi_M^*)^{-1}\circ\psi_{N*}} \Hom_{\sD_A}(M,N)
\]
We define $W_{n/n'}$, $\gr^W_n$  and $W^{>n}$ on $\sD_A$ similarly.
\end{definition}

\begin{remark}\label{rem:weight} 
Since $\sK\sC\sM_A\to \sD_A$ is an equivalence of triangulated categories, the natural distinguished triangles
\[
 W_{n'}\to W_n\to W_{n/n'}\to W_{n'}[1]
\]
in $\sK\sC\sM_A$ give us  natural distinguished triangles
\[
 W_{n'}\to W_n\to W_{n/n'}\to W_{n'}[1]
\]
in $\sD_A$. 
\end{remark}

One uses the weight filtration for inductive arguments, for example:

\begin{lemma}\label{lem:pseudoab} Let $M$ be a finite $A$-cell module. Suppose $N$ is a summand of $M$ in $\sD_A$. Then there is a finite $A$-cell module $M'$ with $N\cong M'$ in $\sD_A$.
\end{lemma}

\begin{proof} 
By proposition~\ref{prop:Whitehead} there  is an isomorphism   $N'\cong N$ in $\sD_A$ with 
with $N'$ an object in $\sC\sM_A$. Thus we may assume that $N$ is a cell module. Since $\KCM_A\to \sD_A$ is an equivalence, $N$ is a summand of $M$ in $\KCM_A$. Write $M=N\oplus N'$ in $\KCM_A$ and let $p:M\to M$ be the projection $M\to N$ followed by the inclusion $N\to M$.

Since $M$ is finite, there is a minimal $n$ with $W_nM\neq 0$. Thus $W_{n-1}N$ is homotopy equivalent to zero and $N\cong W_{\infty/n-1}N$ in $\KCM_A$. Hence, we may assume that $W_{n-1}N=0$  in $\sC\sM_A$. Similarly, we may assume that $M=W_{n+r}M$ and $N=W_{n+r}N$  in $\sC\sM_A$ for some $r\ge0$. We proceed by induction on $r$.

As $A^*(0)=\Q\cdot\id$, it follows that $W_nM=A\otimes_\Q M_0$ for a finite complex of finite dimensional  graded $\Q$-vector spaces $M_0$.
 Indeed, choose a finite bi-graded $A$-basis $\{b_j\}$ for $W_nM$ and let $M_0$ be the finite dimensional $\Q$-vector space spanned by the $b_j$. Since   $W_{n-1}M=0$, all the $b_j$ have Adams degree $n$. Writing $db_j=\sum_i a_{ij}b_i$ and noting that the differential has Adams degree 0, it follows that $|a_{ij}|=0$ for all $i,j$, i.e., $a_{ij}\in \Q\cdot\id$. Consequently $M_0$ is a subcomplex of $M$ and $W_nM=A\otimes_{\Q} M_0$ as an Adams graded dg module.

 But such an $M_0$ is homotopy equivalent to the direct sum of its cohomologies; replacing $M_0$ with $\oplus_nH^n(M_0)[-n]$ and changing notation, we may assume that $d_{M_0}=0$. Thus  $W_nM=A\otimes_\Q M_0$ for $M_0$ a finite dimensional bi-graded $\Q$-vector space; using again the fact that $A(r)=0$ for $r<0$ and $A(0)=\Q\cdot\id$, we see that $W_np=\id \otimes q$ with $q:M_0\to M_0$ an idempotent endomorphism of the bi-graded $\Q$-vector space $M_0$.  Thus $W_nN\cong A\otimes \text{im}(q)$, hence $W_nN$ is homotopy equivalent to a finite $A$-cell module. This also takes care of the case $r=0$.

Using the distinguished triangle 
\[
W_nN\to N\to W_{n+r/n}N\to W_nN[1]
\]
we may replace $N$ with the shifted cone of the map $W_{n+r/n}N\to A\otimes \text{im}(q)[1]$. Since 
$W_{n+r/n}N$ is a summand of $W_{n+r/n}M$, it follows by induction on $r$ that $W_{n+r/n}N$ is homotopy equivalent to a finite cell module, hence the cone of $W_{n+r/n}N\to A\otimes \text{im}(q)$ is homotopy equivalent to a finite cell module as well.
\end{proof}

\begin{definition} Let $\sD_A^{+w}\subset \sD_A$ be the full subcategory of $\sD_A$ with objects $M$ such that $W_nM\cong 0$ for some $n$. 
Similarly, let $\CM_A^{+w}\subset \CM_A$ be the full subcategory with objects $M$ such that $W_nM=0$ for some $n$ and let $\KCM^{+w}_A$ be the homotopy category of  $\CM_A^{+w}$.
\end{definition}

\begin{lemma}\label{lem:Whitehead+} 1. The natural map $\KCM^{+w}_A\to \KCM_A$ is an equivalence of $\KCM^{+w}_A$ with the full subcategory of $\KCM_A$ with objects the $M$ 
such that $W_nM\cong 0$ in $\sK\sC\sM_A$ for $n<<0$. \\
\\
2. The equivalence $\KCM_A\to \sD_A$ induces an equivalence $\KCM_A^{+w}\to \sD_A^{+w}$.
\end{lemma}

\begin{proof} 
Since $\KCM^{+w}_A$ is the homotopy category of the full subcategory $\CM^{+w}_A$ of $\CM_A$, the functor $\KCM^{+w}_A\to \KCM_A$ is a full embedding. Suppose that $W_nM\cong 0$ in $\KCM_A$. We have the cell module $W^{>n}M$ and the 
 distinguished triangle
\[
W_nM\to M\to W^{>n}M\to W_nM[1]
\]
in $\KCM_A$. Thus the map $M\to W^{>n}M$ is an isomorphism in $\KCM_A$; since $W^{>n}M$ is in $\CM^{+w}_A$, the essential image of $\KCM^{+w}_A$ in $\KCM_A$ is as described.

For (2), following definition~\ref{defn:weight}, $W_nM$ is defined by choosing an isomorphism $P\to M$ in $\sD_A$ with $P\in\CM_A$ and taking $W_nM:=W_nP$. Since $W_nP=W_nM\cong0$  in $\sD_A$, it follows that $W_nP\cong0$ in $\KCM_A$, so $P$ is isomorphic to an object in $\KCM_A^{+w}$. Thus $\sD_A^{+w}$ is the essential image of 
$\KCM_A^{+w}$ in $\sD_A$. Since $\KCM_A\to \sD_A$ is an equivalence, this proves (2).
\end{proof}

\begin{remark} Take $M\in \sD_A^{+w}$. Then there is an $n_0$ such that $W_nM\cong0$ for all $n\le n_0$. Indeed, by definition, $W_{n_0}M\cong0$ for some $n_0$. Thus $M\to W^{>n_0}M$ is an isomorphism  in $\sD_A$. If $n<n_0$, then $W_nM\to W_nW^{>n_0}M\cong 0$ is an isomorphism in $\sD_A$.
\end{remark}

Another result using induction on the weight filtration is
 
\begin{lemma} \label{lem:FinGen} Let $M$ be an Adams graded dg $A$-module. \\
\\
1.  $M$ is a finite $A$-cell module if and only if $M$ is free and finitely generated as a bi-graded $A$-module.  \\
\\
2. $M$ is in $\CM^{+w}_A$ if and only if $M$ is free  as a bi-graded $A$-module and there is an integer $r_0$ such that $|m|\ge r_0$ for all $m\in M$.
\end{lemma}

\begin{proof} We first prove (1). Clearly a finite $A$-cell module is free and finitely generated as a bi-graded $A$-module. Conversely, suppose $M$ is  free and finitely generated over $A$; choose a basis $\sB$ for $M$. 

Clearly $W^\sB_nM=0$ for $n<<0$; let $N$ be the minimum integer $n$ such that $W^\sB_nM\neq0$ and let $\sB_N$ be the set of basis elements $b$ of Adams degree $N$. Since $A(0)=\Q\cdot\id$, it follows that $\sB_N$ forms a $\Q$ basis for $W_NM$ and the differential on $\sB_N$ is given by
\[
de_\alpha=\sum_\beta a_{\alpha\beta}e_\beta
\]
with $a_{\alpha\beta}\in\Q$ and $e_\beta\in \sB_N$. Changing the $\Q$ basis for $W^\sB_NM$, we may assume that the subset $\sB_N^0$ of $\sB_N$ of $e_\alpha$ such that $de_\alpha=0$ forms an $\Q$ basis for the kernel of $d$ on the $\Q$-span of $\sB_N$. Since $d^2=0$, the two-step filtration
\[
\sB_N^0\subset \sB_N
\]
exhibits $W_NM$ as a finite cell module.

The result follows by induction on the length of the weight filtration: By induction $W^{>N}_\sB M:=M/W^\sB_nM$ is a finite cell module with basis say $\{b'_j\ |\ j\in J\}$ for some filtration on $J$. Since $M=W_N^\sB M\oplus W_\sB^{>N}M$ as an $A$-module, we just take the union of the two bases, and the concatenation of the filtrations, to present $M$ as a finite cell module. 

The proof of (2) is similar. In fact,  the same proof as for (1) shows that the sub-dg $A$-module $W^\sB_nM$ of $M$ is in $\CM^{+w}_A$ for all $n$ and that we may find an $A$ basis $\sB_n$ for $W^\sB_nM$ and a filtration
\[
\0=\sB_n^{r_0-1}\subset \sB_n^{r_0}\subset\ldots\subset \sB^{2n-1}_n\subset \sB_n^{2n}=\sB_n
\]
that exhibits $W_n^\sB M$ as a cell module. In addition, we may assume that $\sB_i$ with its filtration is just $\sB_n^{2i}$ with the induced filtration, for all $i\le n$. Thus, taking the union of the $\sB_n$ gives the desired basis for $M$, showing that $M$ is in $\CM^{+w}_A$.
\end{proof}

\subsection{Bounded below modules}
\begin{definition}
Let $\sD^+_A\subset \sD_A$ be the full subcategory with objects the  Adams graded dg $A$-modules $M$ having $H^n(M)=0$ for $n<<0$, as usual, we call such an $M$ {\em bounded below}.
\end{definition}
\begin{lemma}\label{lem:BoundedBelow}
Suppose that $A$ is cohomologically connected, and $M$ is  an Adams graded dg $A$-module with  $H^n(M)=0$ for $n<n_0$. Then there is a quasi-isomorphism $P\to M$ with $P$ an $A$-cell module having basis $\{e_\alpha\}$ with $\deg(e_\alpha)\ge n_0$ for all $\alpha$.  If there is an $r_0$ such that $H^n(M)(r)=0$ for $r<r_0$ and all $n$, we may find $P\to M$ as above satisfying the additional condition $|e_\alpha|\ge r_0$ for all $\alpha$.
\end{lemma}

\begin{proof} We first note the following elementary facts: Let $V=\oplus_{n,r}V^n(r)$ be a bi-graded $\Q$-vector space, which we consider as a complex with zero differential. Then the complex $A\otimes_\Q V$ is a cell-module, as a bi-graded $\Q$ basis for $V$ gives a bi-graded $A$ basis with 0 differential. In addition, the map $v\mapsto 1\otimes v$ gives a map
\[
V^n:=\oplus_rV^n(r)\to H^n(A\otimes V).
\]
Finally, suppose there is an $n_0$ such that $V^{n_0}\neq0$ but $V^n=0$ for all $n<n_0$. Then as $H^n(A)=0$ for $n<0$ and $H^0(A)=\Q$, the map
\[
V^{n}\to H^{n}(A\otimes_\Q V)
\]
is an isomorphism for all $n\le n_0$.
 
We begin the construction of $P\to M$ by taking $V$ to be a   bi-graded $\Q$ subspace of  $\oplus_{n\ge n_0}M^n$ representing $\oplus_nH^n(M)$, giving  the map of Adams graded dg $A$  modules 
\[
\phi_{n_0}: P_0:=\oplus_{n\ge n_0 } A\otimes H^n(M)[-n] \to  M.
\]
From the discussion above, we see that  $\phi_{n_0}$ is an isomorphism on $H^n$ for   $n \le n_0$ and a surjection on $H^n$ for   $n> n_0$. If in addition there is an $r_0$ such that $H^n(M)(r)=0$ for $r<r_0$ and all $n$, then $P_0$ has a bi-graded $A$-basis $\sS_0$ with $|v|\ge r_0$ for each $v\in\sS_0$.

Suppose by induction we have constructed a sequence of inclusions of  $A$-cell modules
\[
P_0\to P_1\to \ldots\to  P_{r-1}
\]
and maps of Adams graded dg $A$-modules
\[
\phi_{n_0+i}:P_i\to M
\]
with the following properties:
\begin{enumerate}
\item The $P_i$ have $A$-bases $\sS(i):=\sS_0\cup\ldots\cup \sS_i$. In addition, for all $i\ge1$,   all the elements in $\sS_i$ are of cohomological degree $n_0+i-1$, and for $v\in \sS_i$, $dv$ is in $P_{i-1}$.
\item The map $P_i\to P_{i+1}$ is the one induced by the inclusion $\sS(i)\subset \sS(i+1)$.
\item $\phi_{n_0+i}:P_i\to M$ induces an isomorphism on $H^n$ for $n\le n_0+i$.
\item If  $H^n(M)(r)=0$ for $r<r_0$ and all $n$, then $v\in \sS(i)$ has Adams degree $|v|\ge r_0$.
\end{enumerate}
We now show how to continue the induction. For this,   let $n_r=n_0+r$ and let  $V\subset P_0^{n_r}$ be a  bi-graded $\Q$-subspace of representatives for the kernel of $H^{n_r}(P_{r-1})\to H^{n_r}(M)$. Let 
\[
P_r:=P_{r-1}\oplus  A\otimes_\Q V 
\]
as bi-graded $A$-module, where the differential is given by using the differential on $P_{r-1}$, setting
\[
d((0,1\otimes v))=(v,0)\in P_{r-1}^{n_r}
\]
for $v\in V$ and extending by the Leibniz rule.  Note that, for $v\in V$, there is an $m_v\in M^{n_r-1}$ with $dm_v=\phi_{r-1}(v)$; chosing a bi-graded $\Q$-basis $\sS_r$ for $V$ and extending the assignment $v\mapsto m_v$ from $\sS_r$ to all of $V$ by $\Q$-linearity, we have a $\Q$-linear map
\[
f:V\to M^{n_r-1}
\]
with $d(f(v))=\phi_{r-1}(v)$ for all $v\in V$. Thus, we may define the map of dg $A$-modules
\[
\phi_r: P_r\to M
\]
 by using $\phi_{r-1}$ on $P_r$, $f$ on $1\otimes V$ and extending by $A$-linearity. Clearly $P_r$ is an $A$-cell module with $A$-basis $\sS(r):=\sS(r-1)\cup \sS_r$.
 
 In case $H^n(M)(r)=0$ for $r<r_0$ and all $n$, clearly all bi-homogeneous elements of $V$ have Adams degree $\ge r_0$, so $|v|\ge r_0$ for all $v\in \sS_r$.
 
 We can compute the cohomology of $P_r$ by using the sequence of $A$-cell modules
 \[
 0\to P_{r-1}\to P_r\to  A\otimes_\Q V\to 0,
 \]
 where we consider $V$ as a complex with zero differential, which is split exact as a sequence of bi-graded $A$-modules. The resulting long exact cohomology sequence shows that $P_{r-1}\to P_r$ induces an isomorphism in cohomology $H^n$ for $n<n_r-1$ and we have the exact sequence
 \[
 0\to H^{n_r-1}(P_{r-1})\to  H^{n_r-1}(P_r)\to V
 \xrightarrow{\partial} H^{n_r}(P_{r-1})\to H^{n_r}(P_r)\to 0.
 \]
  In addition, one can compute the coboundary $\partial$ by lifting the element $1\otimes v\in (A\otimes_\Q V)^{n_r-1}$ to the element $(0, 1\otimes v)\in P_r^{n_r-1}$ and applying the differential $d_{P_r}$. From this, we see that the sequence
 \[
 0\to V\xrightarrow{\partial} H^{n_r}(P_{r-1})\to H^{n_r}(P_r)\to 0
 \]
 is exact, hence $ H^{n_r-1}(P_{r-1})\to  H^{n_r-1}(P_r)$ is an isomorphism. This also shows that $\phi_r:P_r\to M$ induces an isomorphism on $H^n$ for $n\le n_r$ and the induction continues.
 
 If we now take $P$ to be the direct limit of the $P_r$, it follows that $P$ is an $A$-cell module with basis elements all in cohomological degree $\ge n_0$, and that the map $\phi:P\to M$ induced from the $\phi_r$ is a quasi-isomorphism. If there is an $r_0$ such that $H^*(M)(r)=0$ for $r<r_0$, then by (4) above, the basis $\sS:=\cup_r\sS(r)$ clearly has  $|e|\ge r_0$ for all $e\in \sS$. This completes the proof.
 \end{proof}

\subsection{Tor and Ext}  The Hom functor $\sHom_A(M,N)$ and tensor product functor $M\otimes_AN$ define exact bi-functors
\begin{align*}
&\sHom_A:\KCM_A^\op\otimes\KCM_A\to \sD_A\\
&\otimes_A:\KCM_A \otimes\KCM_A\to \KCM_A.
\end{align*}
Via proposition~\ref{prop:Whitehead}, these give well-defined derived functors of $\sHom_A$ and $\otimes_A$:
\begin{align*}
&R\sHom_A:\sD_A^\op\otimes\sD_A\to \sD_A\\
&\otimes^L_A:\sD_A \otimes\sD_A\to \sD_A.
\end{align*}
Restricting to $\KCM^f_A$, we have the derived functors for the finite categories
\begin{align*}
&R\sHom_A:\sD_A^{f\op}\otimes\sD^f_A\to\sD^f_A\\
&\otimes^L_A:\sD^f_A\otimes\sD^f_A\to \sD^f_A.
\end{align*}

In both settings, these bi-functors are adjoint:
\[
R\sHom_A(M\otimes^LN,K)\cong R\sHom_A(M,R\sHom_A(N,K)).
\]
We have as well the restriction of $\otimes^L$ to $\sD^{+w}_A$ and $\sD^+_A$:
\begin{align*}
&\otimes^L_A:\sD^{+w}_A \otimes\sD^{+w}_A\to \sD^{+w}_A,
&\otimes^L_A:\sD^+_A \otimes\sD^+_A\to \sD^+_A.
\end{align*}

The derived tensor product makes $\sD_A$ into a triangulated tensor category with unit $\1:=A$ and $\sD^{+w}_A$, $\sD^{+}_A$ and  $\sD^f_A$ are triangulated tensor subcategories. By lemma~\ref{lem:pseudoab}, $\sD^f$ is closed under taking summands in $\sD_A$; this property is obvious for 
$\sD^{+w}_A$ and $\sD^{+}_A$.

Define $M^\vee:=R\sHom_A(M,A)$ and call $M$ {\em strongly dualizable} if the canonical map $M\to M^{\vee\vee}$ is an isomorphism in $\sD_A$. Note that $M$ is strongly dualizable if   $M$ is {\em rigid}, i.e., there exists an $N\in\sD_A$ and morphisms $\delta:A\to M\otimes^L_A N$ and $\epsilon:N\otimes^L_AM\to A$ such that
\begin{align*}
&(\id_M\otimes\epsilon)\circ (\delta\otimes \id_M)=\id_M\\
&(\id_N\otimes\delta)\circ(\epsilon\otimes\id_N)=\id_N
\end{align*}

We have

\begin{proposition}[\hbox{\cite[theorem~5.7]{KrizMay}}] \label{prop:rigid}$M\in \sD_A$ is rigid if and only if $M$ is in $\sD^f_A$, i.e., $M\cong   N$ in $\sD_A$ for some finite $A$-cell module $N$.
\end{proposition}
The precise statement found in \cite[theorem~5.7]{KrizMay} is that $M$ is rigid  if and only if $M$ is a {\em summand} in $\sD_A$ of some finite cell module, so the proposition follows from this and lemma~\ref{lem:pseudoab}; Kriz and May are working in a more general setting in which this lemma does not hold.

\begin{example} For $n\ge 0$,  $\Q(\pm n)\cong (\Q(\pm 1))^{\otimes n}$ and for all $n$, $\Q(n)^\vee\cong\Q(-n)$.
\end{example}

\subsection{Change of ring} If $\phi:A\to A'$ is a homomorphism of Adams graded cdgas,  we have the functor
\[
-\otimes_AA':\sM_A\to \sM_{A'}
\]
which induces a functor on   cell modules and the homotopy category
\[
\phi_*:\KCM_A\to \KCM_{A'}.
\]
Via proposition~\ref{prop:Whitehead}, we have the change of rings functor
\[
\phi_*:\sD_A\to \sD_{A'}
\]
 on the derived category. By proposition~\ref{prop:Whitehead} and  lemma~\ref{lem:Whitehead+}, the respective restrictions of $\phi_*$ define exact tensor functors
\begin{align*}
 &\phi_*:\sD^{+w}_A\to \sD^{+w}_{A'}\\
 &\phi_*:\sD^f_A\to \sD^f_{A'}  .\\
 \end{align*}
 
 \begin{lemma} Let $\phi:A\to A'$ be a homomorphism of  cohomologically connected cdgas. Then  the map $\phi_*$ restricts to
\[
\phi_*:\sD^+_{A}\to \sD^+_{A'} .
\]
\end{lemma}

\begin{proof} Take $M\in \sD^+_{A'}$. By lemma~\ref{lem:BoundedBelow}, there is an integer $N$, a cell $A'$ module $P$ with basis $\{e_\alpha\}$ such that $\deg e_\alpha\ge N$ and a quasi-isomorphism $P\to M$. In fact, looking at the proof of lemma~\ref{lem:BoundedBelow}, we can assume $P$ has an $A$-basis $\sS$  of the form
\[
\sS=\cup_{i=0}^\infty \sS_i
\]
with $\deg(e)\ge N$ for $e\in \sS_0$ and $\deg(e)= N+i-1$ for $e\in \sS_i$ for $i>0$, and such that $d\sS_{r+1}$ is contained in  the $A$-submodule $P_r$ of $P$ generated by $\cup_{i\le r}\sS_i$ for all $r\ge-1$ (where $\sS_{-1}=\0$). Thus we have the sequence of $A$-cell modules
\[
0\to P_{r-1}\to P_r \to P_r/P_{r-1}\to 0
\]
which is split exact as a sequence of bi-graded $A$-modules. Tensoring with $A'$ gives us the 
 sequence of $A'$-cell modules
\[
0\to P_{r-1}\otimes_{A}A'\to P_r\otimes_{A}A' \to P_r/P_{r-1}\otimes_{A}A'\to 0
\]
which is split exact as a sequence of bi-graded $A'$-modules. 

For all $r\ge0$, we have the isomorphism of $A'$-cell modules
\[
P_r/P_{r-1}\otimes_{A}A'\cong \oplus_jA'\<-r_j\>[-m_j]
\]
with $m_j\ge N$. Thus it follows by induction on $r$ and the fact that $A'$ is cohomologically connected that $H^n(P_r\otimes_{A}A' )=0$ for $n<N$. Taking the inductive limit, we see that $H^n(P\otimes_AA')=0$ for $n<N$.

Since  $\phi_*M\cong  P\otimes_{A}A'$ it follows that $\phi_*M$ is in $\sD^+_A$.
\end{proof}

\begin{theorem}[\hbox{\cite[proposition~4.2]{KrizMay}}] \label{thm:EquivCat} If $\phi$ is a quasi-isomorphism, then
\[
\phi_*:\sD_A\to \sD_{A'}
\]
is an equivalence of triangulated tensor categories.
\end{theorem}

Noting the $\phi_*$ is compatible with the weight filtrations, the theorem immediately yields
\begin{corollary}\label{cor:EquivWeightBoundedCat}  If $\phi$ is a quasi-isomorphism, then
\[
\phi_*:\sD_A^{+w}\to \sD_{A'}^{+w}
\]
is an equivalence of triangulated tensor categories.
\end{corollary}

In addition, we have

\begin{corollary} \label{cor:EquivCat} If $\phi$ is a quasi-isomorphism, then
\[
\phi_*:\sD^f_A\to \sD^f_{A'}
\]
is an equivalence of triangulated tensor categories.
\end{corollary}

\begin{proof} Since an equivalence of tensor triangulated categories induces an equivalence on the subcategories of rigid objects, the result follows from theorem~\ref{thm:EquivCat} and proposition~\ref{prop:rigid}.
\end{proof}

\begin{corollary}\label{cor:EquivBoundedCat} Let $\phi:A'\to A$ be a quasi-isomorphism of cohomologically connected cdgas. Then 
\[
\phi_*:\sD^+_{A'}\to \sD^+_A
\]
is an equivalence of tensor triangulated categories.
\end{corollary}

\begin{proof} For $M\in \sD_A$ and integer $n$, we have
\[
H^n(M)=\oplus_r\Hom_{\sD_A}(A,M\<r\>[n]).
\]
and similarly for $A'$.
Since $\phi_*:\sD_{A'}\to \sD_A$ is an equivalence, we have an isomorphism $H^n(M)\cong H^n(\phi_*M)$ for all $M\in \sD_{A'}$. This shows that $\phi_*$ restricts to an isomorphism of the isomorphism classes in $\sD^+_{A'}$ to those in  $\sD^+_A$, which proves the result.
\end{proof}

\begin{proposition} \label{prop:+conservative} Let $\phi:A\to B$ be a map of cdgas. Then $\phi_*:\sD^{+w}_A\to \sD^{+w}_B$ is conservative,  i.e., $\phi_*(M)\cong 0$ implies $M\cong 0$, or equivalently, if $\phi_*(f)$ is an isomorphism then $f$ is an isomorphism.  
\end{proposition}

\begin{proof}  Take $M\in \sD^{+w}$, and let   
\[
\sS:=\{n \ |\ M\cong W^{>n}M\}.
\]
Then $\sS\neq\0$; we claim that either $M\cong 0$ or $\sS$ has a maximal element. Indeed, if $\sS$ has no maximum then $W_nM\cong0$ for all $n$. But since
\[
\colim_nW_nM\to M
\]
is an isomorphism, this implies that $M$ is acyclic, hence $M\cong 0$ in $\sD_A$. 

Thus, we may find a cell module $P$ and quasi-isomorphism $P\to M$ such that $W_{n-1}P=0$, but $W_nP$ is not acyclic.   In particular $P$ has a basis $\{e_\alpha\}$ with $|e_\alpha|\ge n$ for all $\alpha$. If $|e_\alpha|=n$ then since there are no basis elements with Adams grading $<n$, we have
\[
de_\alpha =\sum_ja_{\alpha j}e_j
\]
with $|a_{\alpha j}|=0$, $|e_j|=n$, i.e., $ a_{\alpha j}\in\Q=A(0)$. Since $W_nP$ is not acyclic, it thus follows that $(W_nP)\otimes_A\Q$ is also not acyclic: if $(W_nP)\otimes_A\Q$ were acyclic, this complex would be zero in the homotopy category $\KCM_\Q$, which would make $W_nP$ 0 in $\KCM_A$. As
$W_n(P\otimes_AB)=(W_nP)\otimes_AB$ and
\[
(W_nP)\otimes_A\Q=(W_nP\otimes_AB)\otimes_B\Q
\]
it follows that $P\otimes_AB$ is not isomorphic to zero in $\KCM_B$, and thus $\phi_*(M)$ is non-zero in $\sD^{+w}_B$. 
\end{proof}

\begin{example} Each Adams graded cdga $A$ has a canonical  augmentation $\epsilon:A\to \Q$,  given by projection on $A^0(0)=\Q\cdot\id$. 

 In particular, we have the functor
\[
q:=\epsilon_*:\CM_A\to \sM_\Q,    \ qM:=M\otimes_A\Q
\]
and the exact tensor functors
\begin{align*}
&q:\sD_A\to \sD_\Q,\\
&q^{+w}:\sD^{+w}_A\to \sD^{+w}_\Q,\\
&q^f:\sD^f_A\to \sD^f_\Q.
\end{align*}
Explicitly, $q$ is given on the derived level by $qM:=M\otimes^L_A\Q$.  Assuming $A$ to be cohomologically connected, we have as well the exact tensor functor
\[
q^+:\sD^+_A\to \sD^+_\Q.
\]
\end{example}

\subsection{Finiteness conditions}

 $\sM_\Q$ is just the category of graded $\Q$-vector spaces, so $\sD_\Q$ is equivalent to the product of the unbounded derived categories
\[
\sD_\Q\cong \prod_{n\in\Z}D(\Q).
\]
Similarly
\[
\sD^f_\Q\cong \oplus_{n\in\Z}D^b(\Q),
\]
where $D^b(\Q)$ is the bounded derived category of finite dimensional $\Q$-vector spaces. Finally, 
\[
\sD^{+w}_\Q\cong \bigcup_N \prod_{n\ge N}D(\Q) \ \subset \ \prod_{n\in\Z}D(\Q).
\]
 and
\[
\sD^+_\Q\cong \bigcup_N\prod_{n\in\Z}D^{\ge N}(\Q)\subset \prod_{n\in\Z}D^+(\Q),
\]
where $D^{\ge N}(\Q)\subset D^+(\Q)$ is the full subcategory with objects those complexes $C$ having $H^n(C)=0$ for $n<N$.

\begin{remark}\label{rem:DiffDecomp} The inclusion $\Q\to A$ splits $\epsilon$, identifying $\sD_\Q$, $\sD^+_\Q$, etc., with full subcategories of $\sD_A$, $\sD^+_A$, etc. Under this identification, and the decomposition of $\sD_\Q$ into its Adams graded pieces described above,  the functor $q$ is identified with the functor $\gr^W_*:=\prod_{n\in\Z}\gr^W_n$. Indeed, if $P$ is an $A$-cell module with basis $\{e_\alpha\}$, then as $A(r)=0$ for $r<0$ and $A(0)=\Q\cdot \id$, the differential $d$ decomposes as $d=d^0+d^+$ with 
\[
d^0e_\alpha=\sum_\beta a^0_{\alpha\beta}e_\beta, \ d^+e_\alpha=\sum_\beta a^+_{\alpha\beta}e_\beta
\]
where $|a^0_{\alpha\beta}|=0$, $|a^+_{\alpha\beta}|>0$. Since $d$ has Adams degree 0, it follows that $|e_\beta|<|e_\alpha|$ if $a^+_{\alpha\beta}\neq0$, and $|e_\beta|=|e_\alpha|$ if $a^0_{\alpha\beta}\neq0$. Thus $\gr^W_*P$ is the complex of graded $\Q$-vector spaces with $\Q$ basis $\{e_\alpha\}$ and with  $d_{\gr^W_*P}e_\alpha=d^0e_\alpha$. As $qP$ has exactly the same description, we have the identification of $\gr^W_*$ and $q$ as described.
\end{remark}

\begin{lemma}\label{lem:finite}  Let $M$ be in $\sD_A^{+w}$. Then $M$ is in $\sD^f_A$ if and and only if 
\begin{enumerate}
\item $\gr^W_nM$ is in $D^b(\Q)\subset D(\Q)$ for all $n$.
\item $\gr^W_nM\cong0$ for all but finitely many $n$.
\end{enumerate}
\end{lemma}

\begin{proof} It is clear that $M\in  \sD^f_A$ satisfies the conditions (1) and (2). Conversely, suppose 
 $M\in  \sD^{+w}_A$ satisfies  (1) and (2).  If $M\cong0$, there is nothing to prove, so assume $M$ is not isomorphic to 0. By proposition~\ref{prop:+conservative}, $qM=\prod_n\gr^W_nM$ is not isomorphic to zero. Take $N$ minimal such that $\gr^W_NM$ is not isomorphic to zero. By (2), there is a maximal $N'$ such that $\gr^W_{N'}M$ is not isomorphic to zero.
 
 If $N=N'$, then $M\cong \gr^W_NM$ is in $D^b(\Q)$ by (1), hence $M\cong \oplus_{i=1}^sA\<-N\>[m_i]$, and thus $M$ is in $\sD^f_A$. In general, we apply remark~\ref{rem:weight}, giving
 the distinguished triangle
 \[
 \gr^W_NM\to M\to M^{>N}\to \gr^W_NM[1];
 \]
note that $\gr^W_nM^{>N}\cong0$ for $n>N'$. By induction on $N'-N$,  $M^{>N}$ is in $\sD^f_A$; since $\sD^f_A$ is a full triangulated subcategory of $\sD_A$, closed under isomorphism, it follows that $M$ is in $\sD^f_A$.
\end{proof}

\subsection{Minimal models}  A cdga $A$  is said to be {\em generalized nilpotent} if
\begin{enumerate}
\item As a graded $\Q$-algebra, $A=\Sym^*E$ for some $\Z$-graded $\Q$-vector space $E$, i.e.,
$A=\Lambda^*E_\odd\otimes\Sym^*E_\ev$. In addition,  $E_n=0$ for $n\le0$.
\item For $n\ge0$, let $A_{(n)}\subset A$ be the subalgebra generated by the elements of degree $\le n$. Set $A_{(n+1,0)}=A_{(n)}$ and for $q\ge0$  define $A_{(n+1,q+1)}$ inductively as the subalgebra generated by $A_{(n)}$ and 
\[
A_{(n+1,q+1)}^{n+1}:=\{x\in A^{n+1}_{(n+1)}\ | dx\in A_{(n+1,q)}\}.
\]
Then for all $n\ge0$,
\[
A_{(n+1)}=\cup_{q\ge0}A_{(n+1,q)}.
\]
\end{enumerate}
Note that a generalized nilpotent cdga is automatically connected.

\begin{definition}
Let $A$ be a   cdga.  An {\em $n$-minimal model} of $A$ is a map of cdgas
\[
s:A\{n\}\to A,
\]
with $A\{n\}$ generalized nilpotent and generated (as an algebra) in degree  $\le n$, such that $s$ induces an isomorphism on $H^m$ for $1\le m\le n$ and an injection on $H^{n+1}$. 
\end{definition}

\begin{remark} Let $s:A\{n\}\to A$ be an $n$-minimal model of $A$. Then $A\{n\}_{(n-1)}\subset A\{n\}$ is clearly generalized nilpotent and the inclusion in $A\{n\}$ is an isomorphism in degrees $\le n-1$. Thus
$H^p(A\{n\}_{(n-1)})\to H^p(A\{n\})$ is an isomorphism for $p\le n-1$ and injective for $p=n$, and hence $s:A\{n\}_{(n-1)}\to A$ is an $n-1$-minimal model.
\end{remark}

 Define the above notions for Adams graded cdgas by giving everything an Adams grading.
Let $\sH(cdga)$ be the localization of the category of Adams graded cdgas with respect to maps of cdgas that are quasi-isomorphisms on the underlying complexes. We recall that the category of cdgas has the structure of a model category (see \cite{BousfieldGuggenheim}; the model structure defined there easily passes to the Adams graded case), so that the relation of homotopy between  maps of cdgas is defined. Finally, a generalized nilpotent cdga is cofibrant, so, assuming $A$ to be cohomologically connected,  the minimal model $s:A\{\infty\}\to A$ is a cofibrant replacement ($s$ is a weak equivalence and $A\{\infty\}$ is cofibrant).

\begin{theorem}\label{thm:MinMod} Let $A$ be an Adams graded cdga. \\
1. For each $n=1, 2,\ldots, \infty$, there is an $n$-minimal model  of $A$: $A\{n\}\to A$.\\
\\
2. If $\psi:A\to B$ is a quasi-isomorphism of Adams graded cdgas, and $s:A\{n\}\to A$, $t:B\{n\}\to B$ are $n$-minimal models, then there is an isomorphism of Adams graded cdgas, $\phi:A\{n\}\to B\{n\}$ such that $\psi\circ s$ is homotopic to $t\circ \phi$.
\end{theorem}
See  \cite{BousfieldGuggenheim} or \cite{Quillen} for a proof.

\begin{corollary} \label{cor:ConnectQiso} If $A$ is cohomologically connected, there is a quasi-isomorphism of  Adams graded cdgas $A'\to A$ with $A'$ connected. Similarly, if $\phi:A\to B$ is a map of cohomologically connected  Adams graded cdgas, there is a diagram of  Adams graded cdgas
\[
\xymatrix{
A'\ar[r]\ar[d]&B'\ar[d]\\
A\ar[r]_\phi&B
}
\]
that commutes up to homotopy, with the vertical maps being quasi-isomorphisms, such that $A'$ and $B'$ are connected.
\end{corollary}

\begin{proof}
For the first assertion, just take $A'=A\{\infty\}$. For the second, let
 $B'=B\{\infty\}$. Since $\phi:A\{\infty\}\to A$ is a quasi-isomorphism of $A$-cell modules, $\phi$ is a homotopy equivalence of $A$-cell modules (proposition~\ref{prop:Whitehead}), so taking the tensor product yields a quasi-isomorphism
 \[
A\{\infty\}\otimes_AB\to B.
\]
Clearly $A\{\infty\}\otimes_AB$ is a generalized nilpotent cdga, so we need only  apply theorem~\ref{thm:MinMod}(2). 
\end{proof}

This result together with theorem~\ref{thm:EquivCat},  or corollaries~\ref{cor:EquivWeightBoundedCat}, \ref{cor:EquivCat}, or \ref{cor:EquivBoundedCat}, allows us to replace ``cohomologically connected" with ``connected" in statements involving $\sD_A$,  $\sD^{+w}_A$,  $\sD_A^f$ or $\sD^+_A$.

For example:

\begin{proposition} \label{prop:conservative} Let  $\phi:A\to B$ be a map of  cohomologically connected Adams graded cdgas. Then $\phi_*:\sD^+_A\to \sD^+_B$ is conservative.
\end{proposition}

\begin{proof}  Replacing $A$ and $B$ with $A\{\infty\}$ and $B\{\infty\}$, we may suppose that $A$ and $B$ are connected.  Take $M\in \sD^+_A$ and suppose $\phi_*(M)\cong 0$, i.e. $M\otimes^LB$ is acyclic.  By lemma~\ref{lem:BoundedBelow}, there is a quasi-isomorphism $P\to M$ with $P$ an $A$-cell module having basis $\{e_\alpha\}$ and with $\deg(e_\alpha)\ge n_0$ for some integer $n_0$. If $M$ is not acyclic, we may assume that $n_0$ is chosen so that $H^{n_0}(M)\neq0$. 

Since $A$ is cohomologically connected, this is equivalent to saying that we may take $P$ so that some $e_\alpha$ has $\deg(e_\alpha)=n_0$ and $de_\alpha=0$. As $B$ is connected, $(P\otimes_AB)^n=0$ for $n<n_0$; since $A$ is connected as well, the map $P^{n_0}\to (P\otimes_AB)^{n_0}$ is an isomorphism.
Thus the image of $e_\alpha$ in $P\otimes_AB$ represents a non-zero cohomology class, i.e., $H^{n_0}(P\otimes_AB)\neq 0$. As $M\otimes_A ^LB=P\otimes_AB$, this shows  that $\phi_*$ is conservative.
\end{proof}

\subsection{$t$-structure} \label{subsec:tStructure}To define a $t$-structure on $\sD^+_A$,  $\sD^{+w}_A$ or $\sD_A^f$, one needs to assume that $A$ is cohomologically connected; by   corollaries~\ref{cor:EquivWeightBoundedCat}, \ref{cor:EquivCat}, or \ref{cor:EquivBoundedCat}, we may assume that $A$ is connected.

Define full subcategories $\sD_A^{\le0}$, $\sD_A^{\ge0}$ and $\sH_A$ of $\sD^{+w}_A$ by
\begin{align*}
&\sD_A^{\le0}:=\{M\in \sD^{+w}_A\ |\ H^n(qM)=0\text{ for } n>0\}\\
&\sD_A^{\ge0}:=\{M\in \sD^{+w}_A\ |\ H^n(qM)=0\text{ for } n<0\}\\
&\sH_A:=\{M\in \sD^{+w}_A\ |\ H^n(qM)=0\text{ for } n\neq 0\}.
\end{align*}
The arguments of Kriz-May \cite{KrizMay} show that this defines a $t$-structure  $(\sD_A^{\le0},\sD_A^{\ge0})$  on $\sD^{+w}_A$ with heart $\sH_A$.  Since Kriz-May use $\sD^+_A$ instead of $\sD^{+w}_A$, we give a sketch of the argument here, with the necessary modifications.

\begin{remark} As we have identified the functor $q$ with $\prod_n\gr^W_n$ (remark~\ref{rem:DiffDecomp}) we can describe the category $\sD_A^{\le0}$ as the $M\in \sD^{+w}_A$ such that $H^m(\gr^W_nM)=0$ for all $m>0$ and all $n$. We have a similar description of $\sD_A^{\ge0}$ and
$\sH_A$.
\end{remark}

Recall that an {\em essentially full} subcategory $\sB$ of a category $\sA$ is a full subcategory such that,  if $b\to a$ is an isomorphism in $\sA$ with $b$ in $\sB$, then $a$ is in $\sB$.

\begin{definition}\label{def:tStructure} 
We recall that a {\em $t$-structure} on a triangulated category $\sD$ consists of essentially 
 full subcategories $(\sD^{\le 0}, \sD^{\ge 0})$ of $\sD$ 
 such that
\begin{enumerate}
\item $\sD^{\le0}[1]\subset \sD^{\le0}$, $\sD^{\ge0}[-1]\subset\sD^{\ge0}$
\item $\Hom_\sD(M,N[-1])=0$ for $M$ in $\sD^{\le0}$, $N$ in $\sD^{\ge0}$
\item Each $M\in \sD$ admits a distinguished triangle 
\[
M^{\le 0}\to M\to M^{>0}\to M^{\le 0}[1]
\]
with $M^{\le0}$ in $\sD^{\le0}$, $M^{>0}$ in $\sD^{\ge0}[-1]$.
\end{enumerate}
Write $\sD^{\le n}$ for $\sD^{\le0}[-n]$ and $\sD^{\ge n}$ for $\sD^{\ge0}[-n]$.

A $t$-structure $(\sD^{\le 0}, \sD^{\ge 0})$ is {\em non-degenerate} if 
$A\in \cap_{n\le0}\sD^{\le n}$, $B\in \cap_{n\ge0}\sD^{\ge n}$ imply $A\cong0\cong B$.
\end{definition}

\begin{lemma}\label{lem:truncation} Suppose that $A$ is connected.\\
\\
1. Take $M$ in $\sD_A^{\le0}$. Then there is an $A$-cell module $P\in \CM^{+w}_A$ with basis $\{e_\alpha\}$ such that $\deg(e_\alpha)\le 0$ for all $\alpha$, and a quasi-isomorphism $P\to M$.\\
\\
2.  For $N\in \sD_A^{\ge0}$,  there is an $A$-cell module $P\in \CM^{+w}_A$ with basis $\{e_\alpha\}$ such that $\deg(e_\alpha)\ge 0$ for all $\alpha$, and a quasi-isomorphism $P\to N$.
\end{lemma}

\begin{proof} For (1) choose a quasi-isomorphism $Q\to M$ with $Q$ in $\CM^{+w}_A$. Let $\{e_\alpha\}$ be a basis for $Q$. Decompose the differential $d_Q$ as $d_Q=d_Q^0+d_Q^+$ as in remark~\ref{rem:DiffDecomp}. Making a $\Q$-linear change of basis if necessary, we may assume that the collection $\sS_0$ of $e_\alpha$ with $\deg e_\alpha=0$ and $d^0_Qe_\alpha=0$ forms a basis of 
\[
\ker[d^0:\oplus_{\deg e_\alpha=0}\Q e_\alpha\to \oplus_{\deg e_\beta=1}\Q e_\beta].
\]
Let $\tau^{\le 0}Q$ be the $A$ submodule of $Q$ with basis $\{e_\alpha\ |\ \deg e_\alpha<0\}\cup\sS_0$. We claim that $\tau^{\le 0}Q$ is a subcomplex of $Q$. Indeed, we have
\begin{align*}
d_Qe_\alpha&=d_Q^0e_\alpha+d^+_Qe_\alpha\\
&=\sum_\beta a^0_{\alpha\beta}e_\beta+\sum_\beta a^+_{\alpha\beta}e_\beta
\end{align*}
with $|a^0_{\alpha\beta}|=0=\deg a^0_{\alpha\beta}$, $|a^+_{\alpha\beta}|>0$. Since $A$ is  connected, $\deg a^+_{\alpha\beta}\ge1$. As $d_Q$ has cohomological degree +1, it follows that $\deg e_\beta\le \deg e_\alpha$ if $a^+_{\alpha\beta}\neq0$. Similarly, $\deg e_\beta= \deg e_\alpha+1$ if $a^0_{\alpha\beta}\neq0$. 

Take $e_\alpha$ with $\deg e_\alpha=-1$. Since $d_Q^2=0$, it follows that $(d_Q^0)^2=0$, from which it follows that $e_\beta$ is in $\sS_0$ if $a^0_{\alpha\beta}\neq0$.
 Now take $e_\alpha\in \sS_0$. Write 
\[
de_\alpha=\sum_{\deg b^+_{\alpha\beta}=1} b^+_{\alpha\beta}f^0_\beta+\sum_{\deg b^+_{\alpha\beta}>1} b^+_{\alpha\beta}f_\beta
\]
with the $\{b^+_{\alpha\beta}\}$  being chosen $\Q$ independent in  $A^{*\ge 1}$, the $f_\beta$ in the $\Q$ span of the degree $\le -1$ part of the basis $\{e_\alpha\}$ and the $f^0_\beta$ in $\Q$ span of the degree $0$ part of  $\{e_\alpha\}$. We have
\[
0=d^2e_\alpha=\sum_{\deg b_{\alpha\beta}=1} b^+_{\alpha\beta}d^0(f^0_\beta)+\ldots
\]
with the $\ldots$ involving only the degree $\le0$ part of the basis (and coefficients from $A$).  Since the $b^+_{\alpha\beta}$ are $\Q$ independent, we have $d^0f^0_\beta=0$ for all $\beta$ in the first sum, hence the $f^0_\beta$ are in the $\Q$-span of $\sS_0$. Thus $\tau^{\le 0}Q$ is a subcomplex of $Q$, as claimed.

So far we have only  needed that $Q$ is a cell module. We will now use that $Q$  lies in $\sC\sM_A^{+w}$. We claim that $\tau^{\le 0}Q\to Q$ is a quasi-isomorphism. By   proposition~\ref{prop:+conservative} applied to the augmentation $A\to \Q$, 
\[
q:\sD^{+w}_A\to \sD^{+w}_\Q
\]
 is conservative, thus  it suffices to see that $q\tau^{\le 0}Q\to qQ$ is a quasi-isomorphism. Now, $qQ$ represents $qM\in \sD_\Q$, and by assumption $qM$ is in $\sD^{\le0}_\Q$, hence $qQ$ is in $\sD^{\le 0}_\Q$. But by construction 
$q\tau^{\le 0}Q\to qQ$  is an isomorphism on $H^n$ for all $n\le 0$. Since $H^n(q\tau^{\le 0}Q)=0$ for $n>0$, it follows that $q\tau^{\le 0}Q\to qQ$ is a quasi-isomorphism, as desired.

 For (2), we may assume that $N$ is an object in $\sC\sM_A^{+w}$ and thus  $W_{r_0-1}N=0$ for some $r_0$. The result then follows from lemma~\ref{lem:BoundedBelow}.
\end{proof}

\begin{lemma}\label{lem:Perp} Suppose that $A$ is connected. Then $\Hom_{\sD^{+w}_A}(M,N[-1])=0$ for $M$ in $\sD_A^{\le0}$, $N$ in $\sD_A^{\ge0}$.
\end{lemma}
\begin{proof} By lemma~\ref{lem:truncation} we may assume that $M$ and $N[-1]$ are $A$-cell modules with  bases $\{e_\alpha\}$ for $M$ and $\{f_\beta\}$ for  $N[-1]$ satisfying $\deg e_\alpha\le 0$ and $\deg f_\beta\ge1$ for all $\alpha, \beta$. By lemma~\ref{lem:Whitehead+}, we also have  
\[
\Hom_{\sD^{+w}_A}(M,N[-1])=\Hom_{\KCM^{+w}_A}(M,N[-1]).
\]
But if $\phi:M\to N[-1]$ is a map in $\KCM^{+w}_A$, then $\phi$ is given by a degree 0 map of complexes, so
\[
\phi(e_\alpha)=\sum_\beta a_{\alpha\beta}f_\beta
\]
for $a_{\alpha\beta}\in A$ with $\deg(a_{\alpha\beta})+\deg(f_\beta)=\deg(e_\alpha)$ Since $A^i=0$ for $i<0$, this is impossible.
\end{proof}

\begin{lemma}\label{lem:Truncation} Suppose that $A$ is connected. For $M\in \sD^{+w}_A$, there is a distinguished triangle
\[
M^{\le 0}\to M\to M^{>0}\to M^{\le 0}[1]
\]
with $M^{\le0}$ in in $\sD_A^{\le0}$, $M^{>0}$ in $\sD_A^{\ge1}$. 
\end{lemma}

\begin{proof} We may assume that $M$ is in $\CM^{+w}_A$. We perform exactly the same construction as in the proof of lemma~\ref{lem:truncation}, giving us a sub $A$-cell module $\tau^{\le 0}M$ of $M$ such that
\begin{enumerate}
\item[(a)] $\tau^{\le0}M$ has a basis $\{e_\alpha\}$ with $\deg e_\alpha\le 0$ for all $\alpha$
\item[(b)] The map $q\tau^{\le0}M\to qM$ induced by applying $q$ to $\tau^{\le 0}M\to M$ gives an isomorphism on $H^n$ for $n\le0$.
\end{enumerate}

Let $M^{\le 0}=\tau^{\le 0}M$ and let $M^{>0}$ be the cone of $\tau^{\le 0}M\to M$. This gives us the distinguished triangle 
\[
M^{\le 0}\to M\to M^{>0}\to M^{\le 0}[1]
\]
in $\sD^{+w}_A$. By construction, $M^{\le 0}$ is in $\sD^{+w}_A$. Applying $q$ to the distinguished triangle gives the distinguished triangle in $\sD^{+w}_\Q$
\[
qM^{\le 0}\to qM\to qM^{>0}\to qM^{\le 0}[1];
\]
by (b) and the fact that $H^1(qM^{\le 0})=0$, it follows that $H^n(qM^{>0})=0$ for $n\le0$. Thus $M^{>0}$ is in $\sD^{\ge1}_A$, as desired.
\end{proof}

\begin{theorem} \label{thm:tStructure} Suppose $A$ is  cohomologically connected. Then $(\sD^{\le0}_A,\sD^{\ge0}_A)$ is a non-degenerate $t$-structure on $\sD^{+w}_A$.
\end{theorem}

\begin{proof} Replacing $A$ with its minimal model, we may assume that $A$ is connected. The property (1) of definition~\ref{def:tStructure} is obvious; properties (2) and (3)  follow from lemmata~\ref{lem:Perp} and \ref{lem:Truncation}, respectively.

For $A\in \cap_{n\le0}\sD_A^{\le n}$, it follows that $H^n(qA)=0$ for all $n$, i.e., $qA\cong0$ in $\sD^{+w}_\Q$. Since $q$ is conservative, $A\cong0$ in $\sD^{+w}_A$. The case of $B\in \cap_{n\ge0}\sD_A^{\ge n}$ is similar, hence the $t$-structure is non-degenerate.
\end{proof}

\begin{definition}
Let $\sD^{f,\le0}_A:=\sD_A^f\cap \sD_A^{\le0}$, \ $\sD^{f,\ge0}_A:=\sD_A^f\cap \sD_A^{\ge0}$, \
$\sH^f_A:=\sH_A\cap\sD_A^f=\sD^{f,\le0}_A\cap \sD^{f,\ge0}_A$.
\end{definition}

\begin{corollary}\label{cor:FiniteTStructure}
If  $A$ is cohomologically connected, then $(\sD^{f,\le0}_A,\sD^{f,\ge0}_A)$ is a non-degenerate $t$-structure on $\sD^f_A$ with heart $\sH^f_A$.
\end{corollary}

\begin{proof} Since $\sD^f_A$ is a full triangulated subcategory of $\sD^{+w}_A$, closed under isomorphisms in $\sD^{+w}_A$, all the properties of a non-degenerate $t$-structure are inherited from the non-degenerate $t$-structure on $(\sD^{\le0}_A,\sD^{\ge0}_A)$  on $\sD^{+w}_A$ given by theorem~\ref{thm:tStructure}, except perhaps for the condition (3) of definition~\ref{def:tStructure}. So, take $M\in \sD^f_A$. Since $(\sD^{\le0}_A,\sD^{\ge0}_A)$  is a $t$-structure on $\sD^{+w}_A$, we have a distinguished triangle
\[
M^{\le 0}\to M\to M^{>0}\to M^{\le 0}[1]
\]
with $M^{\le0}$ in $\sD^{\le0}_A$, $M^{>0}$ in $\sD^{\ge0}_A[-1]$. Applying the exact functor $\gr^W_n$ (see remark~\ref{rem:weight})
gives the distinguished triangle
\[
\gr^W_nM^{\le 0}\to \gr^W_nM\to  \gr^W_nM^{>0}\to \gr^W_nM^{\le 0}[1]
\]
in the derived category of $\Q$-vector spaces $D(\Q)$, such that  $\gr^W_nM^{\le0}$ is in $D(\Q)^{\le0}$ and $\gr^W_nM^{>0}$ is in $D(\Q)^{\ge1}$, i.e., $H^n(\gr^W_nM^{\le0})=0$ for $n>0$, $H^n(\gr^W_nM^{>0})=0$ for $n\le 0$. However, since $M$ is in $\sD^f_A$, it follows that $\gr^W_nM$ is in $D^b(\Q)$ for all $n$ and is isomorphic to 0 for all but finitely many $n$  (lemma~\ref{lem:finite}). The long exact cohomology sequence for a distinguished triangle in $D(\Q)$ thus shows that $\gr^W_nM^{\le0}$ and $\gr^W_nM^{>0}$ are  in $D^b(\Q)$ for all $n$ and are isomorphic to zero for all but finitely many $n$. Applying lemma~\ref{lem:finite} again shows $M^{\le 0}$ and $M^{>0}$ are in $\sD^f_A$.
\end{proof}
 
\begin{lemma} \label{lem:HeartProps} (1) The restriction of $\otimes^L$ to $\sH_A$ and $\sH^f_A$ makes these into abelian tensor categories.\\
\\
(2) The weight filtrations on $\sD^{+w}_A$ and $\sD^f_A$ restrict to define exact functorial filtrations on $\sH_A$ and $\sH^f_A$.\\
\\
(3) $\sH^f_A$ is the smallest abelian subcategory of $\sH^f_A$ containing the Tate objects $\Q(n)$, $n\in\Z$ and closed under extensions in $\sH^f_A$.
\end{lemma}

\begin{proof} (1) is more or less obvious: for cell modules $M$ and $N$, we have
$q(M\otimes_AN)\cong qM\otimes_\Q qN$; the K\"unneth formula  for $H^n(qM\otimes_{\Q} qN)$ thus shows that $\sD_A^{\le 0}$ and $\sD^{\ge0}_A$ are closed under $\otimes^L_A$.

For (2), note that the augmentation $\epsilon:A\to \Q$ is a homomorphism of Adams graded cdgas,  and that $q=\epsilon_*$. Thus $q$ is compatible with the weight filtrations on $\sD_A$ and $\sD_\Q$ (and also on the finite categories). In particular, we have
\[
q(\gr^W_nM)\cong \gr^W_nqM.
\]
On the other hand, for $C$ in $\sD^{+w}_\Q$ we have 
\[
C\cong \oplus_m H^m(C)[-m]
\]
Furthermore $H^m(C)$ is isomorphic to its associated weight graded $\oplus_n\gr^W_nH^m(C)$. All this implies that
\[
M \text{ is in }\sD_A^{\le0}\Longleftrightarrow \gr^W_nM \text{ is in }\sD_A^{\le 0}\text{ for all } n
\]
and similarly for $\sD_A^{\ge0}$. Thus, the $t$-structure $(\sD_A^{\le0},\sD_A^{\ge0})$  on $\sD^+_A$ induces a $t$-structure $(W_n\sD_A^{\le0},W_n\sD_A^{\ge0})$  on the full subcategory $W_n\sD^+_A$ with objects the $W_nM$, $M\in\sD^+_A$. The same holds for $\sD_A^{\ge0}$, from which it follows that the truncation functors $\tau_{\le 0}$, $\tau_{\ge 0}$ associated with the $t$-structure $(\sD_A^{\le0},\sD_A^{\ge0})$ commute with the functors $W_n$. This proves (2).

For (3), we argue by induction on the weight filtration. Let $\sH^T_A\subset \sH^f_A$ be any full abelian  subcategory containing all the $\Q(n)$ and closed under extension in $\sH^f_A$. Since $A(0)=\Q\cdot\id$, the full subcategory $\sD^f_A(-n)$ of $\sD^f_A$ consisting of $M$ with $M\cong \gr^W_nM$ is equivalent to the bounded derived category of (ungraded) finite dimensional $\Q$-vector spaces, $D^b(\Q)$, with the equivalence sending a complex $C$ to $\Q(-n)\otimes_\Q C$. The $t$-structure on $\sD^f_A$ restricts to a $t$-structure on $\sD^f_A(-n)$ which is equivalent to the standard $t$-structure on $D^b(\Q)$. 

Thus, if we have $M\in \sH^f_A$, then $\gr^W_nM\cong \Q(-n)^{r_n}$ for some $r_n\ge0$. If $N$ is the minimal $n$ such that $W_nM\neq0$, then we have the exact sequence
\[
0\to \gr^W_NM\to M\to W^{>N}M\to0
\]
By induction on the length of the weight filtration, $W^{>N}M$ is in $\sH^T_A$, hence $M$ is in $\sH^T_A$ and thus $\sH^T_A=\sH^f_A$.
\end{proof}

\begin{lemma} \label{lem:Van}For $N, M\in \sH^f_A$, $n\le m\in\Z$, we have 
\[
\Hom_{\sH^f_A}(W^{>m}M, W_nN)=0
\]
\end{lemma}

\begin{proof} If $M=\Q(-a)$, $N=\Q(-b)$ with $a> b$, then
\[
\Hom_{\sH^f_A}(M,N)=H^0(A(a-b))=0
\]
since $A$ is connected. The result in general follows by induction on the weight filtration.
\end{proof}

\begin{proposition}\label{prop:Tannaka} $\sH^f_A$ is a neutral Tannakian category over $\Q$.
\end{proposition}

\begin{proof} Since $\Q(n)^\vee=\Q(-n)$,  it follows from lemma~\ref{lem:HeartProps} that $M\mapsto M^\vee$ restricts from $\sD^f_A$ to an exact involution on $\sH^f_A$. Since $\sD^f_A$ is rigid, it follows  that $\sH^f_A$ is rigid as well. Also 
\[
\Hom_{\sH^f_A}(\Q(-a),\Q(-b))=\begin{cases}H^0(A(a-b))=0&\text{ if }a\neq b\\
H^0(A(0))=\Q\cdot\id&\text{ if }a=b.\end{cases}
\]
By induction on the weight filtration, this implies that 
$\Hom_{\sH^f_A}(M,N)$ is a finite dimensional $\Q$-vector space for all $M,N$ in $\sH^f_A$. Since the identity for the tensor product is $\Q(0)$, it follows that $\sH^f_A$ is $\Q$ linear.

We have the rigid tensor functor $q:\sH^f_A\to \sH^f_\Q$. Noting that $\sH^f_\Q$ is equivalent to the category of finite dimensional graded $\Q$-vector spaces, composing $q$ with the functor ``forget the grading" from $\sH^f_\Q$ to $\Vect_\Q$ defines the rigid tensor functor
\[
\omega:\sH^f_A\to \Vect_\Q.
\]
The forgetful functor $\sH^f_\Q\to\Vect_\Q$ is faithful, so we need only see that $q:\sH^f_A\to \sH^f_\Q$ is faithful. Sending $M\in\Vect_\Q$ to $\Q(-n)\otimes M$ defines an equivalence of $\Vect_\Q$ with the full subcategory $\gr^W_n\sH^f$ of $\sH^f_A$ consisting of $M$ which are isomorphic to $\gr^W_nM$. Via this identification, we can further identify $q$ with the functor 
\[
M\mapsto \gr^W_*M:=\oplus_n\gr^W_nM.
\]

Let $f:M\to N$ be a map in $\sH^f_A$ such that $\gr^W_nf=0$ for all $n$; we claim that $f=0$. By induction on the length of the weight filtration, it follows that $W^{>n}f=0$, where $n$ is the mininal integer such that $W_nM\oplus W_nN\neq0$. Thus $f$ is given by a map
\[
\tilde f:W^{>n}M\to \gr^W_nN.
\]
But  $\tilde f=0$ by lemma~\ref{lem:Van}, hence $f=0$ as desired.
\end{proof}

\begin{nota}
We denote the truncation to the heart,
\[
\tau_{\le0}\tau^{\ge0}: \sD^{+w}_A\to \sH_A,
\]
by $H^0_A$.
\end{nota}

\subsection{Connection matrices}\label{subsec:connections}  A convenient way to define an  $A$-cell module is by a connection matrix (called a {\em twisting matrix} in \cite{KrizMay}).

Let $(M,d_M)$ be a complex of Adams graded $\Q$-vector spaces. An {\em $A$-connection} for $M$ is a map (of bi-graded $\Q$-vector spaces)
\[
\Gamma:M\to A^+\otimes_\Q M
\]
of Adams degree 0 and cohomological degree 1. One says that $\Gamma$ is {\em flat} if
\[
d\Gamma+\Gamma^2=0.
\]
This means the following:
 $A\otimes_\Q M$ has the standard tensor product differential, so  $d\Gamma:=d_{A^+\otimes_\Q M}\circ \Gamma+\Gamma\circ d_M$ using the usual differential in the complex of maps $M$ to  $A^+\otimes_\Q M$. Also, we extend $\Gamma$ to 
\[
\Gamma:A^+\otimes M\to A^+\otimes M
\]
using the Leibniz rule,  so that $\Gamma^2$ is defined. 

\begin{remark}\label{rem:Dif}
Given a connection $\Gamma: M\to A^+\otimes_\Q M$, define 
\[
d_0: M\to A\otimes_{\Q}M =M\oplus A^{+}\otimes_{\Q}M, \ m\mapsto  d_M m \oplus \Gamma m
\]
 and extend $d_0$ to $ d_\Gamma:A\otimes_{\Q} M\to A\otimes_{\Q}M$ by the Leibniz rule. Then $\Gamma$ is flat  if and only if  $d_\Gamma$   endows $A\otimes_{\Q}M$ with the structure of  a dg $A$-module, i.e.  $d^2_\Gamma=0$.
 \end{remark}
 
If $\Gamma:M\to A^+\otimes_\Q M$ is a connection,  call $\Gamma$ {\em nilpotent} if $M$ admits a filtration by bi-graded $\Q$ subspaces
\[
0=M_{-1}\subset M_0\subset\ldots\subset M_n\subset \ldots\subset M
\]
such that $M=\cup_nM_n$ and such that 
\[
d_M(M_n)\subset M_{n-1};\  \Gamma(M_n)\subset A^+\otimes  M_{n-1}
\]
for every $n\ge0$.

The following result is obvious:

\begin{lemma}\label{lem:ConnCell} Let $\Gamma:M\to A^+\otimes_\Q M$ be a flat nilpotent connection. Then 
 the dg $A$-module $(A\otimes_{\Q}M, d_\Gamma)$ is a cell module. 
 \end{lemma}
 
 Indeed, choosing a $\Q$ basis $\sB$ for $M$ such that $\sB_n:=M_n\cap \sB$ is a $\Q$ basis for $M_n$ for each $n$ gives the necessary filtered $A$ basis for $A\otimes_\Q M$. In addition, we have

\begin{lemma}\label{lem:Nilpotent} Let $\Gamma:M\to A^+\otimes_\Q M$ be a flat connection. Suppose there is an integer $r_0$ such that $|m|\ge r_0$ for all $m\in M$. Then $\Gamma$ is nilpotent.
\end{lemma}

\begin{proof} The proof is essentially the same as that of lemma~\ref{lem:FinGen}(2): If $M$ is concentrated in a single Adams degree $r_0$, then $\Gamma$ is forced to be the zero-map. Thus, taking $M_0=\ker(d_M)\subset M$ and $M_1=M$ shows that $\Gamma$ is nilpotent. In general, one shows by induction on the length of the weight filtration that the restriction of $\Gamma$ to $W_nM:=\oplus_{r\le n}M(r)$ is nilpotent for every $n$, and then a limit argument completes the proof.
\end{proof}

A morphism $f:(M,d_M,\Gamma)\to (M',d_{M'},\Gamma')$ is a map of bi-graded vector spaces
\[
f :=f_0+f^+:M\to A\otimes M'=M'\oplus A^+\otimes M'
\]
such that 
\[
 d_{\Gamma'}f=fd_\Gamma.
\]
 In particular, we may identify the category of complexes of $\Q$-vector spaces with the subcategory consisting of complexes with flat connection 0 and morphisms $f=f^0+f^+$ with $f^+=0$.

 We denote the category of flat nilpotent connections over $A$ by $\Conn_A$. We let $\Conn^{+w}_A$ be the full subcategory consisting of flat nilpotent connections on  $M$ with $M(r)=0$  for $r<<0$, and $\Conn^f_A$ the full subcategory of flat nilpotent connections on $M$ with $M$ finite dimensional over $\Q$. It follows from lemma~\ref{lem:Nilpotent}  that a flat connection on $M$ with $M(r)=0$  for $r<<0$ (or with $M$ finite dimensional over $\Q$) is automatically nilpotent.

Define a tensor operation on $\Conn_A$ by 
\[
(M,\Gamma)\otimes(M',\Gamma'):=(M\otimes M', \Gamma\otimes\id+\id\otimes\Gamma')
\]
with $ \Gamma\otimes\id+\id\otimes\Gamma'$ suitably interpreted as a connection by using the necessary symmetry isomorphisms.

Let  $I$ be the complex
\[
\Q \xrightarrow{\delta} \Q\oplus\Q
\]
with $\Q$ in degree -1, and with connection 0. We have the two inclusions $i_0, i_1:\Q\to I$. Two maps $f,g:(M,\Gamma)\to(M',\Gamma')$ are said to be homotopic if there is a map $h:(M,\Gamma)\otimes I\to (M',\Gamma')$ with $f=h\circ(\id\otimes i_0)$, $g=h\circ(\id\otimes i_1)$. 

\begin{definition} Let    $\sH\Conn_A$ denote the homotopy category of $\Conn_A$, i.e., the objects are the same as $\Conn_A$ and morphisms are homotopy classes of maps in 
$\Conn_A$. Similarly, we have the full subcategories 
\[
\sH(\Conn^f_A)\subset \sH(\Conn^{+w}_A)\subset \sH\Conn_A
\]
with objects $\Conn^f_A$, resp. $\Conn^{+w}_A$.
\end{definition}

If $M$ is an $A$-cell module, then let $M_0$ be the complex of $\Q$-vector spaces $M\otimes_A\Q$. Using the identity splitting $\Q\to A$ to the augmentation $A\to \Q$, we have the canonical isomorphism of 
$A$-modules
\[
A\otimes_\Q M_0\cong M.
\]
Using the decomposition $A=\Q\oplus A^+$, we can decompose the differential on $A\otimes_\Q M_0$ induced by the above isomorphism as
\[
d=d^0+d^+
\]
where $d^0$ maps $\Q\otimes M_0$ to $\Q\otimes M_0$ and $d^+$ maps $\Q\otimes M_0$ to
$A^+\otimes M_0$.

We can thus make $M_0$ into a complex of Adams graded $\Q$-vector spaces by using the differential $d^0$. The map 
\[
d^+:M_0\to A^+\otimes M_0
\]
gives a connection and the flatness condition follows from $d^2=0$. Nilpotence follows from the filtration condition (definition~\ref{def:DgMod}(3b)) for an $A$-basis of $M$.

Conversely, if $(M_0,d^0)$ is a complex of Adams graded $\Q$-vector spaces, and
\[
\Gamma:M_0\to A^+\otimes M_0
\]
is a flat nilpotent connection, make the free Adams graded $A$-module $A\otimes_\Q M_0$ a cell module by taking  $d_\Gamma$ to be the differential (see remark~\ref{rem:Dif} and lemma~\ref{lem:ConnCell}).

It is easy to see that these operations define an equivalence of the category of $A$-cell modules with the category of  flat nilpotent $A$-connections, and that the homotopy relations and tensor products correspond.   Indeed the functor which assigns to a flat nilpotent connection $(M_0,d_{M_0}, \Gamma)$ the cell module $(A\otimes_{\Q} M_0,  d_\Gamma)$ is essentially surjective by the discussion, and the map on Hom groups is an isomorphism.

Define the shift operator by $(M,\Gamma)[1]:=(M[1],-\Gamma[1])$.  Given a morphism $f:(M,\Gamma)\to (M',\Gamma')$ of flat nilpotent connections, decompose $f:M\to A\otimes M'$ as $f:=f^0+f^+$. Define the {\em cone} of $f$ as having underlying complex Cone$(f^0)$, with connection $(-\Gamma[1]\oplus \Gamma')+f^+$. This gives us the cone sequence
 \[
 (M,\Gamma)\to (M',\Gamma')\to \Cone(f)\to (M,\Gamma)[1].
 \]
Using the cone sequences as distinguished triangles makes  $\sH\Conn_A$ into a triangulated tensor category, equivalent to the homotopy category of $A$-cell modules. Via proposition~\ref{prop:Whitehead} we have thus defined an equivalence of  $\sH\Conn_A$ with $\sD_A$ as triangulated tensor categories.  This restricts to equivalences of $\sH(\Conn_A^{+w})$ with $\sD^{+w}_A$ and $\sH(\Conn^f_A)$ with $\sD^f_A$.

The weight filtration in $\sD_A$ can be described in this language: Let $M$ be an Adams graded complex of $\Q$-vector spaces, which we decompose into Adams graded pieces as $M=\oplus_rM(r)$.
Set
\[
W_nM:= \oplus_{r\le n}M(r)
\]
giving us the subcomplex $W_nM$ of $M$. If $\Gamma:M\to A^+\otimes M$ is a flat connection, then as $\Gamma$ has Adams degree 0, it follows that $\Gamma$ restricts to a flat nilpotent connection
\[
W_n\Gamma:W_nM\to A^+\otimes W_nM.
\]
It is easy to see that this filtration corresponds to the weight filtration on $\sD_A$.

 Let $\sH\Conn_A^{+w}\subset \sH\Conn_A$ be the full subcategory of objects $M$ such that $W_nM\cong 0$ for some $n$, and let $\sH\Conn^f_A\subset \sH\Conn_A^{+w}$ be the full subcategory of objects $M$ such that $\oplus_nH^n(M)$ is finite dimensional. It is easy to see that the inclusions $\sH(\Conn^f_A)\subset \sH\Conn^f_A$ and $\sH(\Conn_A^{+w})\subset \sH\Conn_A^{+w}$ are equivalences, giving us the equivalences
\[
\sH\Conn^f_A\sim \sD^f_A,\ \sH\Conn_A^{+w}\sim \sD^{+w}_A.
\]

Now suppose that $A$ is connected. It is easy to see that the standard $t$-structure on  the derived category  $D(\Q)$ of complexes over $\Q$ induces a $t$-structure on the homotopy category  $\sH\Conn^{+w}_A$.  Under the equivalence $\sH\Conn^{+w}_A\sim \sD^{+w}_A$, the $t$-structure  on $\sD^{+w}_A$ defined in section \ref{subsec:tStructure} corresponds to the pair of subcategories $(\sH\Conn_A^{\le 0}, \sH\Conn_A^{\ge 0})$, hence these give the corresponding $t$-structure on $\sH\Conn^{+w}_A$. In particular, the heart $\sH_A$ is equivalent to the category of flat   $A$-connections on Adams graded $\Q$-vector spaces $V$ with $V(r)=0$ for $r<<0$. Denote this latter category by $\Conn^0_A$.

 As we have seen,
 $\sD_A^f$ is equivalent to
 the full subcategory $\sH\Conn_A^f$ of $\sH\Conn_A$ with objects the flat nilpotent connections on complexes $M$ such that $\oplus_nH^n(M)$ is finite dimensional. In case $A$ is connected, we have a similar description of $\sH^f_A$ as the abelian category of flat connections on finite dimensional Adams graded $\Q$-vector spaces, or equivalently, the full subcategory of $\sH\Conn^f_A$ consisting of complexes $M$ with $H^*(M)=H^0(M)$. 

\begin{remark} By lemma~\ref{lem:Nilpotent}, the flat connection $\Gamma$ for an object $(M,\Gamma)$ in $\Conn^0_A$ is automatically nilpotent.  
\end{remark}

We can also give an explicit description of the truncation functors for this $t$-structure in the language of flat nilpotent connections.  Let $(M,d)$ is a  complex of Adams graded $\Q$-vector spaces with a flat nilpotent connection
\[
\Gamma:M\to A^+\otimes M
\]
such that $(M,d,\Gamma)$ is in $\Conn^{+w}_A$. Then we can decompose $\Gamma$  as
\[
\Gamma:=\sum_{i\ge1}\Gamma^{(i)}
\]
by writing 
\[
[A^+\otimes M]^{n+1}=\oplus_{i\ge1}A^i\otimes M^{n-i+1}
\]
and letting $\Gamma^{(i),n}:M^n\to A^i\otimes M^{n-i+1}$ be the composition
\[
M^n\xrightarrow{\Gamma^n}[A^+\otimes M]^{n+1}\to A^i\otimes M^{n-i+1}.
\]
The flatness condition for $\Gamma$ when restricted to the component which maps $M^n$ to  $A^1\otimes M^{n}$ yields the commutative diagram
\[
\xymatrix{\ar[d]_{d^n} M^n \ar[r]^{\Gamma^{(1),n}}&  A^1\otimes_{\Q} M^n \ar[d]^{1\otimes d^{n+1}}\\
M^{n+1} \ar[r]_{\Gamma^{(1),n+1}} & A^1\otimes_{\Q} M^{n+1}. }
\]
This implies that  $\Gamma$  restricts to a flat  connection  $\tau_{\le n}\Gamma$ on the subcomplex $\tau_{\le n}M$:
\[
\tau_{\le n}\Gamma:\tau_{\le n}M\to A^+\otimes\tau_{\le n}M;
\]
$\tau_{\le n}\Gamma$ is nilpotent by lemma~\ref{lem:Nilpotent}.

This in turn implies that $\Gamma$ descends to a connection on the quotient complex $\tau^{>n}M:=M/\tau_{\le n}M$:
\[
\tau^{>n}\Gamma:\tau^{> n}M\to  A^+\otimes\tau^{> n}M
\]
which is in fact a flat nilpotent connection.  Indeed, the only question for flatness is for the terms in $\Gamma^2+d\Gamma$ which factor via $\Gamma$ or $d$ through $A^+\otimes M^{*\le n}$, but which have non-zero image in   $A^+\otimes\tau^{> n}M$. There are three such terms:
\[
\Gamma^{(1),n}\circ \Gamma^{(i+1-n),i}, (1\otimes d^n)\circ \Gamma^{(i+1-n),i}, (1\otimes d^{n-1})\circ \Gamma^{(i+2-n),i}
\]
where we use the convention that $\Gamma^{(0),i}=d^i$. For a term of the first type, the fact that $\Gamma^{(1)}$ commutes with $d$ implies that the composition factors  through $A^{i+1-n}\otimes (M^n/\ker d^n)$. The second term similarly factors through $A^{i+1-n}\otimes (M^n/\ker d^n)$, while the third term goes to zero in $A^{i+2-n}\otimes (M^n/\ker d^n)$.  

As before, the nilpotence of $\tau^{>n}\Gamma$ follows from lemma~\ref{lem:Nilpotent}.

Thus for each $(M,d,\Gamma)$ in $\Conn^{+w}_A$ we have the  sequence of complexes with flat nilpotent connection
\[
0\to(\tau_{\le n}M,d,\tau_{\le n}\Gamma)\to (M,d,\Gamma)\to (\tau^{>n}M,d,\tau^{>n}\Gamma)\to 0
\]
which is exact as a sequence of bi-graded $\Q$-vector spaces.
When we take the associated cell modules, this gives us the canonical distinguished triangle for the $t$-structure we have described for $\sD^{+w}_A$. 

In particular, the truncation functor $H^n_A:=\tau^{\ge n}\tau_{\le n}$ can be explicitly described in the language of flat nilpotent connections. Namely,  the restricted connection
\[
\Gamma^{(1), n}:M^n\to A^1\otimes M^n
\]
 defines a  connection (not necessarily flat) on the  Adams graded $\Q$-vector space $M^n$ for each $n$, and the differential $d$ gives a map  in the category of  connections
 \[
 d^n:(M^n,\Gamma^{(1), n})\to (M^{n+1},\Gamma^{(1), n+1}).
 \]
 In short, $(M,d,\Gamma^{(1)})$ is a complex in the category of connections. Thus $\Gamma^{(1)}$ induces a connection on $H^n(M)$:
 \[
 H^n(\Gamma):=H^n(\Gamma^{(1)}):H^n(M)\to A^1\otimes H^n(M).
 \]
 On $M^n$, the flatness condition for $\Gamma$, when restricted to the component which maps $M^n$ to  $A^2\otimes M^n$, gives the identity:
 \[
(\id\otimes d^{n+1})\circ \Gamma^{(2),n}-\Gamma^{(1),n+1}\circ\Gamma^{(1),n}+\Gamma^{(2),n+1}\circ d^n=0
 \]
 and thus $H^n(\Gamma^{(1)})$ is flat.  $H^n(\Gamma^{(1)})$ is nilpotent by lemma~\ref{lem:Nilpotent}.  
 
 The canonical quasi-isomorphism of complexes
 \[
 \tau^{\ge n}\tau_{\le n}(M,d_M)\to H^n(M,d_M)
 \]
thus  gives rise to a quasi-isomorphism of complexes with flat nilpotent connection
 \[
 \tau^{\ge n}\tau_{\le n}(M,d_M,\Gamma)\to (H^n(M,d_M), H^n(\Gamma^{(1)})).
 \] 
 
 \begin{definition}\label{defn:coLie}
Let $A$ be a cohomologically connected cdga with 1-minimal model $A\{1\}$. We let $QA:=A\{1\}^1$ and let $\del:QA\to\Lambda^2QA$ denote the differential $d:A\{1\}^1\to \Lambda^2A\{1\}^1=A\{1\}^2$.  Then $(QA,\del)$ is  co-Lie algebra over $\Q$. If $A$ is an Adams graded cdga, then $QA$ becomes an Adams graded co-Lie algebra.

In the Adams graded case, we let $\coRep(QA)$ denote the category of co-modules $M$ over $QA$, where $M$ is a bi-graded $\Q$-vector space such that the Adams degrees in $M$ are bounded below.  
\end{definition}

\begin{remark}\label{rem:ConnEquiv}
Let us suppose that $A$ is a generalized nilpotent  Adams graded cdga. Then the co-Lie algebra $QA$ is given by the restriction of $d$ to $A^1$, noting that $d$ factors as
\[
d:A^1\to \Lambda^2A^1\subset A^2.
\]
If now $M$ is an Adams graded $\Q$-vector space (concentrated in cohomological degree 0)
 and $\Gamma:M\to A^+\otimes M$ is a flat connection, then $\Gamma$ is actually a map
 \[
 \Gamma:M\to A^1\otimes M
 \]
 and the flatness condition is just saying the $\Gamma$ makes $M$ into an Adams graded co-module for the co-Lie algebra $QA$.  If in addition the Adams degrees occuring in $M$ have a lower bound, then $\Gamma$ is automatically nilpotent (lemma~\ref{lem:Nilpotent}).
 
 Thus, we have an equivalence of the category $\Conn^0_A$ with  $\coRep(QA)$, which restricts to an equivalence of  $\Conn^0_A\cap\Conn^f_A$ with the category $ \coRep^f(QA)$ of finite dimensional co-modules over $QA$.
 
 Putting this together with the above discussion, we have equivalences
 \[
 \sH_A\sim \Conn^0_A\sim  \coRep(QA)
 \]
 which restrict to equivalences 
 \[
 \sH^f_A\sim \Conn^0_A\cap \Conn^f_A\sim  \coRep^f(QA).
 \]
 \end{remark}

\subsection{Summary} In \cite{KrizMay} the relations between the various constructions we have presented above are discussed. We summarize the main points here.

\begin{definition} 1. Let $H=\Q\cdot \id\oplus \oplus_{r\ge1}$ be an Adams Hopf algebra over $\Q$. We let $\coRep(H)$ denote the abelian tensor category of co-modules $M$ over $H$, where $M$ is a bi-graded $\Q$ vector space such that the Adams degrees in $M$ are bounded below. Let $\coRep^f(H)\subset \coRep(H)$ be the full subcategory of co-modules $M$ such that $M$ is finite dimensional over $\Q$.\\
\\
2. Let $\gamma=\oplus_{r\ge1}\gamma(r)$ be an Adams graded co-Lie algebra over $\Q$. We let $\coRep(\gamma)$ denote the abelian tensor category of co-modules $M$ over $\gamma$, where $M$ is a bi-graded $\Q$ vector space such that the Adams degrees in $M$ are bounded below. Let $\coRep^f(\gamma)\subset \coRep(\gamma)$ be the full subcategory of co-modules $M$ such that $M$ is finite dimensional over $\Q$.
\end{definition}

The Adams grading induces a functorial exact weight filtration on $\coRep(H)$ and $\coRep(\gamma)$ by setting
\[
W_nM:=\oplus_{r\le n}M(r).
\]
The subcategories $\coRep^f(H)$ and $\coRep^f(\gamma)$ are Tannakian categories over $\Q$, with neutral fiber functor the associated graded for the weight filtration $\gr^W_*$.

Let $H_+=\oplus_{r\ge1}H(r)\subset H$ be the augmentation ideal, $\gamma_H:=H_+/H_+^2$ the co-Lie algebra of $H$. For an $H$ co-module $\delta:M\to H\otimes M$ we have the associated $\gamma_H$ co-module $\bar{M}$ with the same underlying bi-graded $\Q$ vector space, and with co-action $\bar\delta:\bar{M}\to \bar{M}\otimes \gamma_H$ given by the composition
\[
M\xrightarrow{\delta}M\otimes H=M\oplus M\otimes H_+\to M\otimes H_+\to M\otimes \gamma_H.
\]
Then the association $M\mapsto \bar{M}$ induces equvalences of filtered abelian tensor categories
\[
\coRep(H)\sim \coRep(\gamma_H),\ \coRep^f(H)\sim\coRep^f(\gamma_H).
\]

For an Adams graded cdga $A$, we have the Adams graded Hopf algebra $\chi_A:=H^0(\bar{B}(A))$ and the 
Adams graded co-Lie algebra $\gamma_A:=\gamma_{\chi_A}$.  We have as well the co-Lie algebra $QA$  defined using the 1-minimal model of $A$ (definition~\ref{defn:coLie}).

\begin{theorem}\label{thm:MayKriz} Let $A$ be an Adams graded cdga. Suppose that $A$ is cohomologically connected.

\medskip
\noindent
(1) There is a functor $\rho:D^b(\coRep^f(\chi_A))\to \sD^f_A$. $\rho$ respects the weight filtrations and sends Tate objects to Tate objects. $\rho$ induces a functor on the hearts
\[
\sH(\rho):\coRep^f(\chi_A)\to \sH^f_A
\]
which is an equivalence of filtered Tannakian categories, respecting the fiber functors $\gr_*^W$.

\medskip
\noindent
(2) Let   $A\{1\}$ be the 1-minimal model of $A$. Then $A\{1\}\to A$ induces an isomorphism of graded Hopf algebras $\chi_{A\{1\}}\to \chi_A$ and graded co-Lie algebras
\[
QA\cong \gamma_{A\{1\}}\cong \gamma_A.
\]

\medskip
\noindent
(3) The functor $\rho$ is an equivalence of triangulated categories if and only if $A$ is 1-minimal.

\medskip
\noindent
(4) Sending a co-module $M\in \coRep(\chi_A)$ to the $\gamma_A$ co-module $\bar{M}$ defines equivalences of neutral Tannakian categories
\[
\coRep(\chi_A)\sim \coRep(\gamma_A);\ \coRep^f(\chi_A)\sim \coRep^f(\gamma_A).
\]
\end{theorem}

Putting this together with our discussion on connections in section~\ref{subsec:connections} gives

\begin{corollary}\label{cor:CatEquiv}  Let $A$ be a cohomologically connected Adams graded cdga. We have equivalences of  filtered abelian tensor categories
\[
\coRep(\chi_A)\sim \coRep(\gamma_A)\sim \coRep(QA)\sim \Conn^0_A
\]
and equivalences of filtered neutral Tannakian categories
\[
\coRep^f(\chi_A)\sim \coRep^f(\gamma_A)\sim \coRep^f(QA)\sim \Conn^0_A\cap \Conn^f_A.
\]
\end{corollary}

\section{Relative theory of cdgas}\label{sec:RelDga} The theory of cdgas over $\Q$ generalizes to a large extent to cdgas over a cdga $\sN$. In this section, we give the main constructions in this direction that we will need. 

\subsection{Relative minimal models} We fix a base cdga  $\sN$.  A {\em cdga over $\sN$} is a cdga $\sA$  together with a homomorphism of cdgas $\phi:\sN\to \sA$. An {\em augmented cdga over $\sN$} has in addition a splitting $\pi:\sA\to \sN$ to $\phi$. The same notions apply for an Adams graded cdga $\sA$ over an Adams graded cdga $\sN$. From now on, we assume we are in the Adams graded setting.

The notions of generalized nilpotent algebras and minimal models (over $\Q$) extend without difficulty to augmented cdgas over $\sN$. Specifically:

\begin{definition}\label{defn:GenNilp}
 An Adams graded cdga $\sA$ over $\sN$  is said to be {\em generalized nilpotent over $\sN$} if
\begin{enumerate}
\item As a bi-graded $\sN$-algebra, $\sA=\Sym^*E\otimes\sN$ for some Adams graded $\Z$-graded $\Q$-vector space $E$, i.e.,
${\sA}=\Lambda^*E^\odd\otimes\Sym^*E^\ev\otimes\sN$, where the parity refers to the cohomological degree. In addition,  $E(r)^n=0$ if $n\le0$ or if $r\le 0$.
\item For $n\ge0$, let $\sA_{(n)}\subset \sA$ be the $\sN$-subalgebra generated  by the subspace $E^{\le n}$ of $E$ consisting of elements of cohomological degree $\le n$. Set $\sA_{(n+1,0)}=\sA_{(n)}$ and for $q\ge0$  define $\sA_{(n+1,q+1)}$ inductively as the $\sN$-subalgebra generated by $\sA_{(n)}$ and 
\[\sA_{(n+1,q+1)}^{n+1}:=
\{x\in \sA_{(n+1)}^{n+1}\ | dx\in {\sA}_{(n+1,q)}.\}
\]
Then for all $n\ge0$,
\[
\sA_{(n+1)}=\cup_{q\ge0}\sA_{(n+1,q)}.
\]
\end{enumerate}
\end{definition}

\begin{remark} \label{rem:Filt}  We can phrase the condition (2) above differently: For each $n\ge0$, $E^{\le n+1}$ has an increasing exhaustive bi-graded filtration
\[
E^{\le n}=F_0E^{\le n+1}\subset F_1E^{\le n+1}\subset\ldots\subset F_mE^{\le n+1}\subset\ldots\subset  E^{\le n+1}
\]
such that
\[
d(F_mE^{\le n+1}\otimes\sN)\subset \Sym^*(F_{m-1}E^{\le n+1})\otimes\sN
\]
Indeed, if $\sA=\Sym^*E\otimes\sN$ satisfies  (2), define $F_mE^{\le n+1}$  by 
\[
F_mE^{\le n+1}\otimes 1=(E^{\le n+1}\otimes1)\cap \sA_{(n+1,m)}^*.
\]
Conversely, it is easy to see that the existence of such a filtration $F_*E^{\le n+1}$ for all $n$ implies (2).
\end{remark}

\begin{lemma}\label{lem:CohConn} Let $\sA$ be a generalized nilpotent cdga over a cdga $\sN$. If $\sN$ is cohomologically connected, then so is $\sA$.
\end{lemma}

\begin{proof} Write $\sA=\Sym^*E\otimes\sN$ as an $\sN$ algebra, with $E$ a bi-graded $\Q$-vector space, so that 
the conditions (1) and (2) of definition~\ref{defn:GenNilp} are satisfied. Let $\sA^p:=\Sym^*E^{*\le p}\otimes\sN$. Then $\sA^p\subset \sA$ is a  generalized nilpotent sub-cdga of $\sA$; as $\sA$ is the limit of the $\sA^p$, it suffices to show that $\sA^p$ is cohomologically connected. We may therefore assume that $E=E^{*\le p}$ and that the result holds for $\sA^{p-1}$.

By remark~\ref{rem:Filt}  there is an exhaustive increasing bi-graded filtration 
\[
E^{*\le p-1}= F_0E\subset\ldots\subset F_nE\subset\ldots\subset E
\]
 on  $E$ so that  $d(F_nE\otimes \sN)\subset \Sym^*F_{n-1}E\otimes \sN$ for all $n>0$, and such that $\sA_n:=\Sym^*F_nE\otimes \sN$ (with differential induced from $\sA$) is a generalized nilpotent cdga over $\sN$. It thus suffices to show that $\sA_n$   is cohomologically connected for each $n>0$. By induction on $n$, it suffices to show that for all $n>0$, the quotient complex
 \[
 \bar{\sA}_n:=\sA_n/\sA_{n-1}
 \]
 has vanishing cohomology $H^i$  for $i\le0$.
 
 Writing $E_n:=F_nE/F_{n-1}E$, we have the filtration $G_*$ on 
 $\Sym^*F_nE/\Sym^*F_{n-1}E$ with 
 \[
\gr^G_m(\Sym^*F_nE/\Sym^*F_{n-1}E)\cong  \Sym^mE_n\otimes\Sym^*F_{n-1}
 \]
for all $m\ge1$ and $G_0(\Sym^*F_nE/\Sym^*F_{n-1}E)=0$.

 By the Leibniz rule, the subspace
 \[
 \bar{\sA}_{n,m}:=G_m(\Sym^*F_nE/\Sym^*F_{n-1}E)\otimes \sN\subset \bar{\sA}_n
 \]
is a subcomplex. Thus it suffices to show that  $H^i(\bar{\sA}_{n,m}/\bar{\sA}_{n,m-1})=0$  for $i\le0$.

For all $n>0$  and $m\ge1$, we have
\[
\bar{\sA}_{n,m}/\bar{\sA}_{n,m-1}\cong \Sym^mE_n\otimes\Sym^*F_{n-1}E\otimes\sN
\]
with differential  $\id\otimes d_{\sN}$.  Since $E_n$ has cohomological degree $\ge1$ and $\sN$ has vanishing cohomology in degrees $<0$, it follows that  $H^i(\bar{\sA}_{n,m}/\bar{\sA}_{n,m-1})=0$  for $i\le0$.  This completes the proof.
\end{proof}

\begin{definition}
Let $\sA$ be an augmented Adams graded cdga over $\sN$.  An {\em $n$-minimal model} over $\sN$ of $\sA$ is a map of augmented Adams graded cdgas over $\sN$
\[
s: \sA\{n\}_\sN\to  \sA,
\]
with  $\sA\{n\}_\sN$ generalized nilpotent over $\sN$, $\sA\{n\}_\sN=\Sym^*E\otimes\sN$, with $E$ satisfying the conditions of definition~\ref{defn:GenNilp}, such that  $\deg e\le n$ for all $e\in E$, and such that $s$ induces an isomorphism on $H^m$ for $1\le m\le n$ and an injection on $H^{n+1}$. 

If the base-cdga $\sN$ is understood, we call an $n$-minimal model over $\sN$ a {\em relative $n$-minimal model}.
\end{definition}

\begin{proposition} \label{prop:RelMin} 
Let $\sN$ be a cohomologically connected Adams graded cdga, $\sA$ an augmented Adams graded cdga over $\sN$.   Then\\
\\
1. For each $n$, there is an $n$-minimal model over $\sN$: $\sA\{n\}_\sN\to \sA$.\\
\\
2. If $\psi:\sA\to \sB$ is a quasi-isomorphism of augmented cdgas over $\sN$, and $s:\sA\{n\}_\sN\to \sA$, $t:\sB\{n\}_\sN\to \sB$ are  relative $n$-minimal models, then there is an isomorphism of augmented cdgas over $\sN$, $\phi:\sA\{n\}_\sN\to \sB\{n\}_\sN$ such that $\psi\circ s$ is homotopic to $t\circ \phi$.\\
\\
 Suppose that   $\sA$ is also cohomologically connected. Then $\sA\{n\}_\sN\to \sA$ induces an isomorphism on $H^i$ for all $i\le n$. In particular, the map   $\sA\{\infty\}_\sN\to \sA$ is a quasi-isomorphism.
\end{proposition}

\begin{proof}  This result is the relative analog of theorem~\ref{thm:MinMod} and the proof is essentially the same (see \cite{BousfieldGuggenheim, Quillen} for the details in the absolute case).

 The construction of the $n$-minimal model over $\sN$ is essentially the same as for cdgas over $\Q$ except that we use both the cohomological degree and the Adams degree for induction: The augmentation gives a canonical decomposition of $\sA$ as
\[
\sA=\sN\oplus \sI
\]
with $\sI$ an Adams graded dg $\sN$-ideal in $\sA$. Let  
$E_{10}(1)\subset \sI^1(1)$ be a $\Q$-subspace of representatives for $H^1(\sI(1))$, in cohomological degree 1, with Adams degree $1$.  We have the evident mapping
\[
 E_{10}(1)\otimes_\Q\sN\to \sA
\]
using the $\sN$-module structure, which extends to 
\[
\Sym^* E_{10}(1)\otimes_\Q\sN\to \sA
\]
using the algebra structure. Clearly this is a map of augmented cdgas over $\sN$, and induces an isomorphism on $H^1(-)(1)$, because $\sN(r)=0$ for $r<0$ and $\sN(0)=\Q\cdot\id$.

One then proceeds as in the case $\sN=\Q$ to adjoin elements in degree 1 and Adams degree 1 to successively kill elements in the kernel of the map on $H^2(-)(1)$. Since $\sN(r)=0$ for $r<0$ and $\sN(0)=\Q\cdot\id$,  this does not affect $H^1$ in Adams degree $\le 1$. Thus we have constructed a bi-graded $\Q$-vector space $E_1(1)$, of Adams degree 1 and cohomological degree 1, a generalized nilpotent cdga over $\sN$, $\sA_{1,1}:=\Sym^*E_1(1)\otimes \sN$ and a map of cdgas over $\sN$, $\sA_{1,1}\to \sA$, that induces an isomorphism on $H^1(-)(1)$  and an injection on  $H^2(-)(1)$  . 

This completes the Adams degree $\le1$ part for the construction of the 1-minimal model. So far, we have not used the cohomological connectivity of $\sN$, this comes in now: Use the canonical augmentation of $\sA_{1,1}$ to write $\sA_{1,1}=\sN\oplus \sI_{1,1}$.

\begin{claim} $H^p(\sI_{1,1}(r))=0$ for $r>1$, $p\le 1$.
\end{claim}

To prove the claim, we use the same filtration that we used in the proof of lemma~\ref{lem:CohConn}. The same induction argument as in lemma~\ref{lem:CohConn}, using of course the cohomological connnectedness of $\sN$,  shows that the lowest degree cohomology of $\sI_{1,1}(r)$ comes from $\oplus_{i=1}^{r-1}\Sym^iE_1(1)\otimes H^1(\sN(r-i))$ plus $\Sym^rE_1(1)\otimes H^0(\sN(0))$. Since all the elements of $E_1(1)$ have cohomological degree 1, this proves the claim.

For the $n$-minimal model with $n>1$,  we continue the construction, first adjoining elements of Adams degree 1 and cohomological degree 2 to generate all of $H^2(\sA)(1)$, and then adjoining elements of Adams degree 1 and cohomological degree 2 to kill the kernel on $H^3(-)(1)$.
 Continuing in this manner   gives the generalized nilpotent cdga over $\sN$, 
 \[
 \sA_{1,n}:=\Sym^*E_n(1)\otimes \sN, 
 \]
with   $E_n(1)$ in Adams degree $1$ and cohomological degree $1,\ldots,  n$, together with a map over $\sN$, $\sA_{1,n}\to \sA$, that induces an isomorphism on $H^i(-)(1)$ for $1\le i\le n$ and an injection for $i=n+1$. If we are in the case $n=\infty$, we just take the colimit of the $\sA_{1,n}$. In addition, writing $\sA_{1,n}=\sN\oplus\sI_{1,n}$, we have
\[
H^p(\sI_{1,n}(r))=0\text{ for }r>1, p\le1.
\]

Now suppose we have constructed bi-graded $\Q$-vector spaces 
\[
E_n(1)\subset E_n(2)\subset\ldots\subset E_n(m)
\]
(for fixed $n$ with $1\le n\le \infty$) with $E_n(j)$ having Adams degrees $1,\ldots, j$ and cohomological degrees $1,\ldots, n$, a differential on $\sA_{m,n}:=\Sym^*E_n(m)\otimes \sN$ making  $\sA_{m,n}$ a generalized nilpotent cdga over $\sN$, and a map $\sA_{m,n}\to \sA$ of cdgas over $\sN$ that is an isomorphism on $H^i(-)(j)$ for $1\le i\le n$, $j\le m$, and an injection for $i=n+1$, $j\le m$. In addition,  writing $\sA_{n,m}=\sN\oplus\sI_{n,m}$, we have
\begin{equation}\label{eqn:InductiveVan}
H^p(\sI_{m,n}(r))=0\text{ for }r>m, p\le1.
\end{equation}
We extend $E_n(m)$ to $E_n(m+1)$ by simply repeating the construct for $E_n(1)$ described above, but working in Adams degree $m+1$ rather than 1; using \eqref{eqn:InductiveVan} allows us to start the construction by adjoining generators for $H^1(\sI(m+1))$, just as in the case of Adams weight 1.
Again, as $\sN(r)=0$ for $r<0$ and $\sN(0)=\Q\cdot\id$,  the inclusion $\sA_{m,n}\to \sA_{m+1,n}$ is an isomorphism  in Adams degree $\le m$. In addition,  the argument used to prove the claim shows that  \eqref{eqn:InductiveVan} extends from $m$ to $m+1$ and the induction goes through.

Taking $E_n:=\cup_mE_n(m)$, we thus have a differential on $\sA\{n\}_\sN:=\Sym^*E_n\otimes \sN$ making  $\sA\{n\}_\sN$ a generalized nilpotent cdga over $\sN$, and a map $\sA\{n\}_\sN\to \sA$ of cdgas over $\sN$ that is an isomorphism on $H^i(-)$ for $1\le i\le n$ and an injection for $i=n+1$, completing the proof of (1).

For (2), the construction of an isomorphism $\phi$ between two $n$-minimal models over $\sN$ is also the same   as for $\sN=\Q$,  using again  induction on the Adams degree. 

To prove (3), the assumption that $\sN$ is cohomologically connected passes to all generalized nilpotent cdgas over $\sN$ (lemma~\ref{lem:CohConn}), in particular, $\sA\{n\}_\sN$ is cohomologically connected. If  in addition $\sA$ is cohomologically connected, then  $\sA\{n\}_\sN\to \sA$  automatically induces an isomorphism on $H^i$ for $i\le 0$. Combining this with (1) proves (3).
\end{proof}

\begin{remark} A generalized nilpotent cdga over $\sN$ is automatically a cell-module over $\sN$.  Indeed, for $\sA=\Sym^*E\otimes\sN$ satisfying the conditions of definition~\ref{defn:GenNilp}, one has the filtration on $E^{\le n}$ given by remark~\ref{rem:Filt}. Combining this filtration with the filtration by degree on $\Sym^*E$  gives a  filtration on $\Sym^*E$ which exihibits $\sA$ as an $\sN$-cell module.
\end{remark}

\subsection{Relative bar construction}\label{subsec:RelativeBar}  One forms the bar construction for a cdga $\sA$ over $\sN$ just as for cdgas over $\Q$, replacing $\otimes_\Q$ with $\otimes_\sN$. However, for this construction to have good cohomological properties, one should replace $\sA$ with a quasi-isomorphic cdga $\sA'$ which is a cell module over $\sN$, so that $\otimes_\sN=\otimes^L_\sN$. This is accomplished by using the minimal model $\sA\{\infty\}$. In any case, we give the ``pre-derived" definition for an arbitrary cdga $\sA$ over $\sN$.

\begin{definition} Let $\sA$ be an augmented Adams graded cdga over $\sN$. Define the simplicial cdga $B^\nai_\bullet(\sA/\sN)$ by
\[
B^\nai_\bullet(\sA/\sN):=\sA^{\otimes_\sN[0,1]}
\]
The inclusion $\{0,1\}\to [0,1]$ makes $B^\nai_\bullet(\sA/\sN)$ a simplicial cdga over $A\otimes A$. 
Given two (possibly equal) augmentations $\epsilon_1.\epsilon_2:\sA\to \sN$, set
\[
B^\nai_\bullet(\sA/\sN,\epsilon_1,\epsilon_2):=B^\nai_\bullet(\sA/\sN)\otimes_{\sA\otimes\sA}\sN.
\]
and let $\bar{B}^\nai_\sN(\sA,\epsilon_1,\epsilon_2)$ be the total complex associated to 
$B^\nai_\bullet(\sA/\sN,\epsilon_1,\epsilon_2)$.
\end{definition}

\begin{remark} 
If $\sA$ is a generalized nilpotent algebra over $\sN$, then $\bar{B}^\nai_\sN(\sA,\epsilon_1,\epsilon_2)$ has a natural structure of an  $\sN$-cell module. In addition, since $\sA$ is Adams graded, $\sA(r)=0$ for $r<0$, hence $\sA^{\otimes_\sN n}$ is in $\CM^{+w}_\sN$ for each $n\ge0$ and $\bar{B}^\nai_\sN(\sA,\epsilon_1,\epsilon_2)$  is in $\CM_\sN^{+w}$.  Finally, if $\epsilon_1=\epsilon_2=\epsilon$, then $\bar{B}^\nai_\sN(\sA,\epsilon)$ has the natural structure of a dg Hopf algebra in $\CM_\sN^{+w}$, and thus a Hopf algebra in $\sD^{+w}_\sN$.
\end{remark}

\begin{definition}  Let $\sA$ be an augmented Adams graded cdga over $\sN$ with augmentation $\epsilon$. Suppose that $\sN$ is cohomologically connnected and let $\sA\{\infty\}_\sN\to\sA$ be the relative minimal model of $\sA$ over $\sN$. Define
\[
B_\bullet(\sA/\sN):=B^\nai_\bullet(\sA\{\infty\}_\sN/\sN),\ \bar{B}_\sN(\sA,\epsilon):=\bar{B}^\nai_{\sN}(\sA\{\infty\}_\sN,\epsilon\{\infty\}).
\]
\end{definition}

\begin{remark}\label{rem:RelCoLieAlg} Still supposing $\sN$ to be cohomologically connected, we may apply the truncation functor
\[
H^0_\sN:\sD^{+w}_\sN\to \sH_\sN
\]
to the  dg Hopf algebra  $\bar{B}_\sN(\sA,\epsilon)$ in $\sD^{+w}_\sN$, giving us the Hopf algebra $H^0_\sN(\bar{B}_\sN(\sA,\epsilon))$ in $\sH_\sN$. We may therefore also form the co-Lie algebra object $\gamma_{\sA/\sN}$ in $\sH_\sN=\Conn^0_\sN$:
\[
\gamma_{\sA/\sN}:=H^0_\sN(\bar{B}_\sN(\sA,\epsilon))_+/H^0_\sN(\bar{B}_\sN(\sA,\epsilon))_+^2
\]
with $H^0_\sN(\bar{B}_\sN(\sA,\epsilon))_+\subset H^0_\sN(\bar{B}_\sN(\sA,\epsilon))$ the augmentation ideal.

We let $\bar{B}_{\bullet\le m}(\sA/\sN,\epsilon)$ denote the restriction of the simplicial object 
 $\bar{B}_{\bullet}(\sA/\sN,\epsilon)$ to the full subcategory $\{[0],\ldots,[m]\}$ of $\Ord$, and 
 $\bar{B}^{\le m}_\sN(\sA,\epsilon)\subset \bar{B}_\sN(\sA,\epsilon)$ the associated total complex of 
 $\bar{B}_{\bullet\le m}(\sA/\sN,\epsilon)$.

If we suppose that $\sA$ is in $\sD^f_\sN$, then $H^0_\sN(\bar{B}^{\le m}_\sN(\sA,\epsilon))$ is in 
 $\sH_\sN^f$ for each $m$, hence $H^0_\sN(\bar{B}_\sN(\sA,\epsilon))$ has the structure of an ind-Hopf algebra in  $\sH_\sN^f$ with
 \[
 H^0_\sN(\bar{B}_\sN(\sA,\epsilon))=\colim_{m\to\infty}H^0_\sN(\bar{B}^{\le m}_\sN(\sA,\epsilon))
 \]
 in $\sH_\sN$.
 \end{remark}

\subsection{Base-change} We consider a quasi-isomorphism $\phi:\sN'\to \sN$ of cohomologically connected cdgas. Given an augmented cdga $\sA$ over $\sN$ with augmentation $\epsilon:\sA\to \sN$, we have $\sA=\sI\oplus \sN$, with $\sI$ the kernel of $\epsilon$. In particular, $\sI$ is a (non-unital) $\sN$-algebra. Via $\phi$, we make $\sI$ a (non-unital) $\sN'$-algebra, and thus give $\sA':=\sI\oplus\sN'$ the structure of a cdga over $\sN'$, with augmentation $\epsilon':\sA'\to \sN'$ the projection on $\sN'$ with kernel $\sI$. 

This construction yields the commutative diagram of cdgas
\begin{equation}\label{eqn:BaseChange}
\xymatrix{
\sA'\ar[r]^{\phi'}\ar@<3pt>[d]^{\epsilon'}&\sA\ar@<3pt>[d]^\epsilon\\
\sN'\ar@<3pt>[u]^{p'}\ar[r]_\phi&\sN\ar@<3pt>[u]^p}
\end{equation}
with $\phi$ and $\phi'$ quasi-isomorphisms.

Now let $f':\sA'\{n\}_{\sN'}\to \sA'$ be a relative  $n$-minimal model over $\sA'$ over $\sN'$.  Since the composition $\phi'f':\sA'\{n\}_{\sN'}\to \sA$ is an $\sN'$-module map, $\phi'f'$ factors through a unique map 
\[
f:\sA'\{n\}_{\sN'}\otimes_{\sN'}\sN\to \sA
\]
of cdgas over $\sN$. Similarly, the $\sN'$-augmentation  of $\sA'\{n\}_{\sN'}$ induces an $\sN$-augmentation of $\sA'\{n\}_{\sN'}\otimes_{\sN'}\sN$, making $f$ a map of augmented cdgas over $\sN$. 

\begin{lemma} $f:\sA'\{n\}_{\sN'}\otimes_{\sN'}\sN\to \sA$ is a relative $n$-minimal model of $\sA$ over $\sN$.
\end{lemma} 

\begin{proof} As $\sA'\{n\}_{\sN'}$ is a generalized nilpotent algebra over $\sN'$, with generators in degree $\le n$, the same follows for  $\sA'\{n\}_{\sN'}\otimes_{\sN'}\sN$ as an algebra over $\sN$. $\phi$ is a quasi-isomorphism, so $\phi_*:\sD_{\sN'}\to \sD_\sN$ is an equivalence of triangulated categories. We can compute cohomology of a dg module via maps in the derived category; as
$\sA'\{n\}_{\sN'}$ is an $\sN'$-cell module, we have $\phi_*(\sA'\{n\}_{\sN'})=\sA'\{n\}_{\sN'}\otimes_{\sN'}\sN$, hence the canonical map
\[
\sA'\{n\}_{\sN'}\to \sA'\{n\}_{\sN'}\otimes_{\sN'}\sN
\]
is a quasi-isomorphism of cdgas. Since $\phi':\sA'\to \sA$ is a quasi-isomorphism, and $\sA'\{n\}_{\sN'}\to \sA'$ is a relative $n$-minimal model, the map on $H^i$ induced by $f$ is an isomorphism for $1\le i\le n$ and an injection for $i=n+1$, i.e., $f:\sA'\{n\}_{\sN'}\otimes_{\sN'}\sN\to \sA$ is a relative $n$-minimal model.
\end{proof}

\begin{remark}\label{rem:QuasiIsoReplacement} Still assuming $\sN$ and $\sN'$ cohomologically connected, write $\sA\{n\}_\sN$ for the $n$-minimal model $\sA'\{n\}_{\sN'}\otimes_{\sN'}\sN$. We have the change of rings isomorphism
\[
\sA'\{n\}_{\sN'}^{\otimes_{\sN'}m}\otimes_{\sN'}\sN\to \sA\{n\}_{\sN}^{\otimes_{\sN}m}
\]
and the quasi-isomorphism
\[
\sA'\{n\}_{\sN'}^{\otimes_{\sN'}m}\to \sA'\{n\}_{\sN'}^{\otimes_{\sN'}m}\otimes_{\sN'}\sN
\]
Thus  on the bar construction
\[
\bar{B}_{\sN'}^\nai(\sA'\{n\}_{\sN'},\epsilon')\xrightarrow{\alpha} \bar{B}_{\sN'}^\nai(\sA'\{n\}_{\sN'},\epsilon')\otimes_{\sN'}\sN\xrightarrow{\beta}
\bar{B}_{\sN}^\nai(\sA\{n\}_{\sN},\epsilon)
\]
the map $\alpha$ is a quasi-isomorphism and the map $\beta$ is an isomorphism. 

In particular, taking $n=\infty$, we have the canonical isomorphism
\[
\phi_*(H^0_{\sN'}(\bar{B}_{\sN'}(\sA',\epsilon')))\cong H^0_{\sN}(\bar{B}_{\sN}(\sA,\epsilon))
\]
of Hopf algebra objects in $\sH_{\sN}$. Since $\phi_*:\sH_{\sN'}\to \sH_\sN$ is an equivalence, we are thus free to replace $\sN$ with a quasi-isomorphic $\sN'$ in a study of $H^0_{\sN}(\bar{B}_{\sN}(\sA,\epsilon))$. For instance, we may use the minimal model $\sN\{\infty\}\to \sN$ as a replacement for $\sN$.
\end{remark}

\subsection{Connection matrices}\label{subsec:RelConnMat} Generalized nilpotent algebras over $\sN$ fit well into the connection matrix point of view described in section~\ref{subsec:connections}. Indeed, suppose that $\sA=\Sym^*E\otimes\sN$ is generalized nilpotent over  $\sN$, with augmentation  $\epsilon:\sA\to\sN$ induced by writing $\Sym^*E=\Q\oplus\Sym^{*\ge1}E$.

Using the augmentation of $\sN$, we write $\sN=\Q\cdot\id\oplus \sN^+$, which writes $\sA$ as 
\[
\sA=\Sym^*E\otimes\id\oplus \Sym^*E\otimes\sN^+.
\]
Thus the differential on $\sA$ is completely determined by its restriction to $\Sym^*E\otimes\id$, giving the decomposition
\[
d=d^0+\Gamma
\]
with $d^0$ a differential on $\Sym^*E$ and $\Gamma:\Sym^*E\to  \Sym^*E\otimes\sN^+$ a flat  connection. In addition, $(\Sym^*E,d^0)$ an Adams graded cdga over $\Q$ with augmentation  $\epsilon^0$ induced by the projection to  $\Sym^0E=\Q$. Finally, the connection $\Gamma$ is nilpotent since $\Sym^*E$ has all Adams degrees $\ge0$ (lemma~\ref{lem:Nilpotent}).

Using the tensor structure in the category of flat nilpotent connections,  the flat nilpotent connection $\Gamma:\Sym^*E\to \Sym^*E\otimes\sN^+$ gives rise to a flat nilpotent connection on $(\Sym^*E)^{\otimes n}$ for all $n$. These fit together  to  give a flat nilpotent connection on the bar construction:
\[
\bar{B}(\Gamma):\bar{B}((\Sym^*E,d^0),\epsilon^0)\to
\bar{B}((\Sym^*E,d^0),\epsilon^0)\otimes\sN^+.
\]
This defines a Hopf algebra object in $\Conn_\sN$.

\begin{proposition} \label{prop:co-LieConnection} Let $\sN$ be cohomologically connected. The  $\sN$-cell module corresponding to  $\bar{B}((\Sym^*E,d^0),\epsilon^0)$ with flat nilpotent connection $\bar{B}(\Gamma)$ is isomorphic to  $\bar{B}^\nai_\sN(\sA,\epsilon)$, as dg Hopf algebra objects in $\CM^{+w}_\sN$. 
\end{proposition}

\begin{proof} We check instead the equivalent statement that the dg Hopf algebra in $\Conn_\sN$ corresponding to $\bar{B}^\nai_\sN(\sA,\epsilon)$ is $(\bar{B}((\Sym^*E,d^0),\epsilon^0),\bar{B}(\Gamma))$. 

We note that we have canonical isomorphisms
\[
\sA^{\otimes_\sN n}\cong (\Sym^*E)^{\otimes_\Q n}\otimes_\Q \sN=
(\Sym^*E)^{\otimes_\Q n}\otimes\id\oplus (\Sym^*E)^{\otimes_\Q n}\otimes_\Q \sN^+ 
\]
respecting differentials and multiplications. Tracing this isomorphism through the definition we have given of the flat nilpotent connection on $\bar{B}((\Sym^*E,d^0),\epsilon^0)$ completes the proof. 
\end{proof}

\subsection{Semi-direct products}\label{subsec:SemiDirect} Let $\epsilon:\sA\to \sN$ be an augmented Adams graded cdga over $\sN$. We suppose that $\sN$ is generalized nilpotent and that $\sA$ is generalized nilpotent over $\sN$. We let $G_\sA:=\Spec H^0(\bar{B}(\sA))$, $G_\sN:=\Spec H^0(\bar{B}(\sN))$
 be the $\Q$-algebraic group schemes defined with respect to the canonical augmentations $\sA\to\Q$, $\sN\to \Q$. The $\sN$-algebra structure $\pi^*:\sN\to \sA$ induces the map of algebraic groups $\pi:G_\sA\to G_\sN$; the augmentation $\epsilon$ gives a splitting $s:G_\sN\to G_\sA$ to $\pi$. 

\begin{lemma} The map $\pi$ is flat.
\end{lemma}

\begin{proof} As $\Q$-algebras, $H^0(\bar{B}(\sA))$ and $H^0(\bar{B}(\sN))$ are polynomial algebras on $\sA^1$, $\sN^1$ respectively, and the map
\[
H^0(\bar{B}(\pi^*)):H^0(\bar{B}(\sN))\to H^0(\bar{B}(\sA))
\]
is just the polynomial extension of the linear injection 
\[
\pi^*:\sN^1\to \sA^1, 
\]
i.e., $H^0(\bar{B}(\pi^*))$ identifies $H^0(\bar{B}(\sA))$ with a polynomial extension of $H^0(\bar{B}(\sN))$.
\end{proof}

\begin{lemma}\label{lem:Kernel} Let $e$ denote the identity in $G_\sN$. The fiber $\pi^{-1}(e)$ is canonically isomorphic to $\Spec H^0(\bar{B}(\sA\otimes_\sN\Q))$ as  group schemes over $\Q$.
\end{lemma}

\begin{proof} We have the natural map of Hopf algebras
\[
 H^0(\bar{B}(\sA))\otimes_{H^0(\bar{B}(\sN))}\Q\to H^0(\bar{B}(\sA\otimes_\sN\Q)).
\]
Writing $\sA=\Sym^*E\otimes\sN$ as an $\sN$-algebra, $H^0(\bar{B}(\sA\otimes_\sN\Q))$ is a polynomial algebra on $(\Sym^*E)^1$, while $H^0(\bar{B}(\sA))$ is the polynomial algebra on $\sA^1=(\Sym^*E)^1\oplus \sN^1$, and $H^0(\bar{B}(\sN))$ is the polynomial algebra on $\sN^1$. This shows that the above map is an algebra isomorphism.
\end{proof}

Set $K:=\Spec  H^0(\bar{B}(\sA\otimes_\sN\Q))=\Spec H^0(\bar{B}(\Sym^*E))$. The splitting $s$ gives an action of $G_\sN$ on $K$ and an isomorphism of $G_\sA$ with the semi-direct product
\[
G_\sA\cong K\ltimes G_\sN.
\]
Let $K_s$ denote the $\Q$-group scheme $K$ with this $G_\sN$-action.

On the other hand, we have seen  (proposition~\ref{prop:co-LieConnection}) that writing $\sA=\Sym^*E\otimes \sN$ gives $\Sym^*E$ a flat nilpotent connection
\[
\Gamma: \Sym^*E \to \sN^+\otimes \Sym^*E 
\]
and an isomorphism of 
$H^0_\sN(\bar{B}_\sN(\sA))$ with $H^0(\bar{B}(\Sym^*E))$ as Hopf algebras in $\Conn^0_\sN$.

Replacing $\sN$ with its 1-minimal model, and noting that $\Conn^0_\sN\sim \Conn^0_{\sN\{1\}}$ we have the canonical structure of $H^0(\bar{B}(\Sym^*E))$ as a Hopf algebra in the category of co-modules over the co-Lie algebra $Q\sN=\gamma_\sN$  (remark~\ref{rem:ConnEquiv}). But this category is equivalent to the category of representations of $G_\sN$, giving us another action of $G_\sN$ on $K$.

\begin{theorem}\label{thm:Kernel}  The action of $G_\sN$ on $K=\Spec H^0(\bar{B}(\Sym^*E))$ induced by the splitting $s$ is the same as the action given by the  flat nilpotent $\sN$-connection $\Gamma$ on $\Sym^*E$. In other words, there is an isomorphism 
\[
K_s\cong  \Spec H^0_\sN(\bar{B}_\sN(\sA))
\]
as $\Q$-group schemes with $G_\sN$-action.
\end{theorem}

\begin{proof} It suffices to check that the two co-actions of the co-Lie algebra $\gamma_\sN$ are the same, in fact, it suffices to check that the  two co-actions of $\gamma_\sN$ on the co-Lie algebra $\gamma_{\Sym^*E}$ of $K$ are the same.  

By Quillen's theorem (theorem~\ref{thm:MayKriz}(2)), we can identify the co-Lie algebras  $\gamma_\sA$, $\gamma_\sN$ and $\gamma_{\Sym^*E}$ with $Q\sA$, $Q\sN$ and $Q\Sym^*E$, respectively.  Since we are assuming $\sA$ and $\sN$ are both generalized nilpotent, $Q\sA$, $Q\sN$ and $Q\Sym^*E$  are the respective co-Lie algebras
\[
d_\sA:\sA^1\to \Lambda^2\sA^1,\ d_\sN:\sN^1\to \Lambda^2\sN^1,\ d_E:E^1\to \Lambda^2E^1.
\]
On the level of co-Lie algebras, the splitting $s$ is just the decomposition of $\sA^1=(\Sym^*E\otimes\sN)^1$ as
\[
\sA^1=(\Sym^*E)^1\oplus \sN^1.
\]
The co-action of $\sN^1$ on $\sA^1$ determined by the splitting $s$ is thus given by $d_\sA$ followed by the projection of $\Lambda^2\sA^1$ on $\sN^1\otimes\sA^1$ via the isomorphism
\[
\Lambda^2\sA^1=\Lambda^2((\Sym^*E)^1\oplus \sN^1)\cong 
\Lambda^2(\Sym^*E)^1\oplus \sN^1\otimes(\Sym^*E)^1\oplus \Lambda^2\sN^1.
\]
This induces the co-action of $\sN^1$ on $(\Sym^*E)^1$ by taking the composition
\[
(\Sym^*E)^1\to \sA^1\xrightarrow{d_\sA} \Lambda^2\sA^1\to  \sN^1\otimes(\Sym^*E)^1.
\]
Via our identifications, this gives us the co-action of $\gamma_\sN$ on $\gamma_{\Sym^*E}$ determined by the section $s$.

On the other hand, the flat nilpotent connection $\Gamma$ on $\Sym^*E$ giving the  isomorphism of 
$H^0_\sN(\bar{B}_\sN(\sA))$ with $H^0(\bar{B}(\Sym^*E))$ in $\Conn^0_\sN$  is just the restriction of $d_\sA$ to $\Sym^*E$ followed by the projection of $\sA=\sN\otimes \Sym^*E$ to $\sN^+\otimes\Sym^*E$. However, by reasons of degree, the restriction of $d_\sA$ to $(\Sym^*E)^1=E^1$ decomposes as
\[
d_\sA:E^1\to \Lambda^2E^1\oplus \sN^1\otimes E^1
\]
from which it follows that $\Gamma:E^1\to \sN^1\otimes E^1$ is the same as the co-action defined by $s$. 
\end{proof}

\section{Motives over a base} \label{sec:Motives}
This section summarizes  the material we need from  the work of 
Cisinski-D\'eglise \cite{CisinskiDeglise} and {\O}stvar-R\"ondigs \cite{RoendigsBigMot}.
Together with the content of  sections \ref{sec:CycAlg} and \ref{sec:NModMot},  
 this will be developed in a forthcoming article \cite{LevineTateMotive} by the second author. In this section, we always assume that $k$ admits resolution of singularities.

\subsection{The construction} We summarize the main points of the construction of the category $\DM^\eff(S)$ of effective motives over $S$, and the category $\DM(S)$ of motives over $S$, from \cite{CisinskiDeglise}. Although $S$ is allowed to be a quite general scheme in \cite{CisinskiDeglise}, we restrict ourselves to the case of  a  base-scheme $S$ that is  separated, smooth and essentially of finite type over a field. We let  $\Sch_S$  denote the category of finite type separated $S$-schemes and let $\Sm/S$ denote the full subcategory of $\Sch_S$  consisting of smooth $S$-schemes.

For $X\in \Sm/S$, $Y\in \Sch_S$, define the group of finite $S$-correspondences $c_S(X,Y)$ as the free abelian group on the integral closed subschemes $W\subset X\times_SY$ with $W\to X$ finite and surjective over an irreducible component of $X$.

For $X, Y$ in $\Sm/S$, $Z\in \Sch_S$,  let $p_{XY}$, $p_{YZ}$ and $p_{XZ}$ be the evident projections from $X\times_SY\times_SZ$. One checks that the formula
\begin{equation}\label{eqn:CorComp}
W\circ W':=p_{XZ*}(p_{XY}^*(W)\cdot p_{YZ}^*(W'))\in c_S(X,Z)
\end{equation}
where $\cdot$ is the intersection product, is well-defined for all $W\in c_S(X,Y)$, $W'\in c_S(Y,Z)$\footnote{Even though $Z$ may be singular, one can locally embed $Z$ in an $\A^N_S$ and compute the intersection multiplicites there}; this follows from the fact that $\supp(W)\times_SZ\cap X\times_S\supp(W')$ is finite over $X$ and each irreducible component of this intersection dominates a component of $X$. This is called the {\em composition of correspondences}.

We start with the category $\SmCor(S)$. Objects are the same as $\Sm/S$, morphisms are
\[
\Hom_{\SmCor(S)}(X,Y):=c_S(X,Y)
\]
with composition law given by the formula \eqref{eqn:CorComp}. Define the abelian category of {\em presheaves with transfer} on $\Sm/S$, $\PST(S)$, as the category of presheaves of abelian groups on $\SmCor(S)$. We have the representable presheaves $\Z_S^{tr}(Z)$ for $Z\in \Sm/S$ by 
$\Z_S^{tr}(Z)(X):=c_S(X,Z)$ and pull-back maps given by the composition of correspondences. In fact, the same formula defines $\Z_S^{tr}(Z)$ for $Z\in\Sch_S$.

One gives the category of complexes $C(\PST(S))$  the  {\em Nisnevich local} model structure (which we won't need to specify). The homotopy category is equivalent to the (unbounded) derived category $D(\Sh^{tr}_\Nis(S))$, where $\Sh^{tr}_\Nis(S)$ is the full subcategory of $\PST(S)$  consisting of the presheaves with transfer which restrict to Nisnevich sheaves on $\Sm/S$.  

The  operation 
\[
\Z_S^{tr}(X) \otimes_S^{tr}\Z_S^{tr}(X') :=\Z_S^{tr}(X\times_SX') 
\]
extends to  a tensor structure $\otimes^{tr}_S$ making $\PST(S)$ a tensor category: one forms the {\em canonical left resolution} $\sL(\sF)$ of a presheaf $\sF$ by taking the canonical surjection \[
\sL_0(\sF):=\bigoplus_{ X\in \Sm/S, s\in \sF(X)}\Z_S^{tr}(X)\xrightarrow{\phi_0} \sF
\]
setting $\sF_1:=\ker\phi_0$ and iterating. One then defines
\[
\sF\otimes_S^{tr}\sG:=H_0(\sL(\sF)\otimes_S^{tr} \sL(\sG))
\]
noting that $\sL(\sF)\otimes_S^{tr} \sL(\sG)$ is defined since both complexes are degreewise direct sums of representable presheaves.

The restriction of $\otimes^{tr}_S$ to the subcategory of cofibrant objects in $C(\Sh^{tr}_\Nis(S))$ induces a tensor operation  $\otimes^L_S$ on $D(\Sh^{tr}_\Nis(S))$ which makes $D(\Sh^{tr}_\Nis(S))$ a triangulated tensor category.

\begin{definition}[\cite{CisinskiDeglise}]   $\DM^\eff(S)$ is  the localization of the triangulated category $D(\Sh^{tr}_\Nis(S))$ with respect to the localizing category generated by the complexes $\Z_S^{tr}(X\times\A^1)\to \Z_S^{tr}(X)$.  Denote  by $m_S(X)$ the image of $\Z_S^{tr}(X)$ in $\DM^\eff(S)$.
\end{definition}

\begin{remark} 1.  $\DM^\eff(S)$ is a triangulated tensor category  with tensor product $\otimes_S$  induced from the tensor product   $\otimes^L_S$  via the localization map
\[
Q_S:D(\Sh^{tr}_\Nis(S))\to \DM^\eff(S),
\]
   and satisfying $m_S(X)\otimes_S  m_S(Y)=m_S(X\times_SY)$.\\
\\
2. One has the model category $C(\PST_{\A^1}(S))$ with underlying category $C(\PST(S))$ defined as the Bousfield localization of $C(\PST(S))$ with respect to the  
complexes
\begin{enumerate}
\item For each {\em elementary Nisnevich square}  with $X\in\Sm/$:
\[
\xymatrix{
W\ar@{^{(}->}[r]\ar@{=}[d]&X'\ar[d]^f\\
W\ar@{^{(}->}[r]&X}
\]
one has the complex 
\[
\Z^{tr}_S(X'\setminus W)\to \Z^{tr}_S(X\setminus W)\oplus \Z^{tr}_S(X')\to \Z^{tr}_S(X)
\]
Recall that the square above is an elementary Nisnevich square if $f$ is \'etale, the horizontal arrows are closed immersions of reduced schemes and the square is cartesian.
\item For $X\in \Sm/S$, one has the complex $\Z_S^{tr}(X\times\A^1)\to \Z_S^{tr}(X)$.
\end{enumerate}

 The homotopy category of $C(\PST_{\A^1}(S))$ is equivalent to $\DM^\eff(S)$.
\end{remark}

\begin{definition}\label{def:TateObj} Let $T^{tr}$ be the presheaf with transfers
\[
T^{tr}:=\coker(\Z_S^{tr}(S)\xrightarrow{i_{\infty*}}\Z_S^{tr}(\P_S^1))
\]
and let $\Z_S(1)$ be the image in $\DM^\eff(S)$ of $T^{tr}[-2]$. Let 
\[
\otimes T^{tr}: C(\PST(S))\to C(\PST(S))
\]
be the functor $C\mapsto C\otimes_S^{tr}T^{tr}$.
\end{definition}

Let $\Spt_{T^{tr}}(S)$ be the model category of $\otimes T^{tr}$ spectra in $C(\PST_{\A^1}(S))$, i.e., objects are sequence $E:=(E_0, E_1, \ldots)$, $E_n\in C(\PST(S))$, with bonding maps
\[
\epsilon_n:E_n\otimes^{tr}_S T^{tr}\to E_{n+1}.
\]
Morphisms are given by sequences of maps in $C(\PST(S))$ which strictly commute with the respective bonding maps.

 The model structure on the category of $T^{tr}$-spectra is defined by following the construction of Hovey \cite{HoveySymSpec}. The weak equivalences are the {\em stable weak equivalences}: for each $E\in \Spt_{T^{tr}}(S)$ there is a canonical fibrant model $E\to E^f$, where $E^f:=(E^f_0, E^f_1,\ldots)$ with each $E^f_n$ fibrant in 
$C(\PST_{\A^1}(S))$ and the map
\[
E^f_n\to \sHom(T^{tr}, E^f_{n+1})
\]
adjoint to the bonding map $E^f_n\otimes^{tr}_S T^{tr}\to E^f_{n+1}$ is a weak equivalence in the model category  $C(\PST_{\A^1}(S))$.

\begin{definition} The ``big" category of triangulated motives over $S$, $\DM(S)$, is the homotopy category of $\Spt_{T^{tr}}(S)$.
\end{definition}

\begin{remark} Concretely, a stable weak equivalence $f:(E_0, E_1,\ldots)\to (F_0, F_1,\ldots)$ is a map such that, for each $X\in C(\PST(S))$, the map
\[
\colim_n\Hom_{\DM^\eff(S)}(X\otimes^{tr}_S(T^{tr})^{\otimes n}, E_n)\xrightarrow{f_{n*}}
\colim_n\Hom_{\DM^\eff(S)}(X\otimes^{tr}_S(T^{tr})^{\otimes n}, F_n)
\]
is an isomorphism.
\end{remark}

We will use the following result from \cite{CisinskiDeglise}.
\begin{theorem}[\hbox{\cite[section~10.4]{CisinskiDeglise}}]\label{thm:MotivesMain} Suppose that $S$ is in $\Sm/k$ for a field $k$, take $X$ in $\Sm/S$, and let $m_k(X)$, $m_S(X)$ denote the motives of $X$ in $\DM(k)$, $\DM(S)$, respectively. Then there is a natural isomorphism
\[
\Hom_{\DM(S)}(m_S(X), \Z(n)[m])\cong \Hom_{\DM(k)}(m_k(X), \Z(n)[m])
\]
\end{theorem}

\begin{remark} \label{rem:MotCohIso} By Voevodsky's embedding theorem \cite[chapter~V, theorem~3.2.6]{FSV} the functor $\DM^\eff_\gm(k)\to \DM^\eff_-(k)$ induces a full embedding
\[
\DM_\gm(k)\to \DM(k),
\]
hence $\Hom_{\DM(k)}(m_k(X), \Z(n)[m])$ is motivic cohomology in the sense of Voevodsky \cite[chapter~V]{FSV},  that is  
\[
\Hom_{\DM(k)}(m_k(X), \Z(n)[m])=H^m(X,\Z(m)).
\]
\end{remark}

\subsection{Tensor structure}
The tensor structure on $C(\PST(S))$ induces a ``tensor  operation" on the spectrum category by the usual device of choosing a cofinal subset $\N\subset \N\times\N$, $i\mapsto (n_i,m_i)$, with $n_{i+1}+m_{i+1}=n_i+m_i+1$ for each $i$: each pair of $T^{tr}$ spectra $E:=(E_0, E_1,\ldots)$ and $F:=(F_0, F_1,\ldots)$ gives rise to a $T^{tr}$ bispectrum
\[
E\boxtimes^{tr}_S F:=\begin{pmatrix}&\vdots&\\\ldots&E_i\otimes^{tr}_S F_j&\ldots\\&\vdots&\end{pmatrix}
\]
with vertical and horizontal bonding maps induced by the bonding maps for $E$ and $F$, respectively. The vertical bonding maps use in addition the symmetry isomorphism in $C(\PST_{\A^1}(S))$. Finally, the choice of the cofinal $\N\subset \N\times\N$ converts a bispectrum to a spectrum.

Of course, this is not even associative, so one does not achieve a tensor operation on $\Spt_{T^{tr}}(S)$, but $\boxtimes^{tr}_S$ (on cofibrant objects) does pass to the localization $\DM(k)$, and gives rise there to a tensor structure, making $\DM(S)$ a tensor triangulated category. We write this tensor operation as $\otimes_S$, as before.

\subsection{Tate motives}  In $\DM(X)_\Q$ we have the full subcategory of {\em Tate motives over $X$}, $\DTM(X)$, this being the full triangulated subcategory of $\DM(X)_\Q$ closed under isomorphism and generated by the Tate motives $\Q_X(n)$, $n\in \Z$. Since $\Q_X(n)\otimes\Q_X(m)\cong \Q_X(n+m)$, $\DTM(X)$ is a tensor triangulated subcategory of $\DM(X)_\Q$.

Just as for the case of motives over a field, the category $\DTM(X)$ admits a canonical weight filtration, and, in case $X$ satisfies the Beilinson-Soul\'e vanishing conjectures, a $t$-structure with heart generated by the Tate objects $\Q_X(n)$. In fact, the results of \cite{LevineTate} apply directly, so we will content ourselves here with giving the relevant definitions.

\begin{definition} Let $W_n\DTM(X)$ denote the full triangulated subcategory of $\DTM(X)$ generated by the Tate motives $\Q_X(-a)$ with $a\le n$. Let $W_{[n,m]}\DTM(X)$ be the full triangulated subcategory of $\DTM(X)$ generated by the Tate motives $\Q_X(-a)$ with $n\le a\le m$, and let $W^{>n}\DTM(X)$ be the full triangulated subcategory of $\DTM(X)$ generated by the Tate motives $\Q_X(-a)$ with $a> n$.
\end{definition}

\begin{lemma}\label{lem:TateMotCoh} For $X\in\Sm/k$ there is a natural isomorphism
\[
\Hom_{\DMT(X)}(\Q_X(a),\Q_X(b)[m])\cong H^m(X,\Q(b-a))
\]
\end{lemma}

\begin{proof} This follows directly from theorem~\ref{thm:MotivesMain} and remark~\ref{rem:MotCohIso}, noting that $\otimes\Q_X(a)$ is invertible in $\DMT(X)\subset\DM(X)_\Q$.
\end{proof}

\begin{lemma}  $\DMT(X)$ is a rigid tensor triangulated category.
\end{lemma}

\begin{proof} 
The unit $\1$  for the tensor operation is $\Q_X(0)$. It suffices to check that the generators $\Q_X(n)$ of  $\DMT(X)$ admit a dual  (see e.g. \cite[part~I, IV.1.2]{MixMot}). Setting $\Q_X(n)^\vee=\Q_X(-n)$, with maps $\delta:\1\to\Q_X(n)^\vee\otimes\Q_X(n)$, $\epsilon:\Q_X(n)\otimes \Q_X(n)^\vee\to\1$ being the canonical isomorphisms shows that $\Q_X(n)$ has a dual.
\end{proof}

\begin{theorem} \label{thm:MixTateWeight} 1. $(W_n\DTM(X),W^{>n}\DTM(X))$ is a $t$-structure on $\DTM(X)$ with heart consisting of 0-objects.\\
\\
2. Denote the truncation functors for the $t$-structure  $(W_n\DTM(X),W^{>n}\DTM(X))$ by
\begin{align*}
&W_n:\DTM(X)\to W_n\DTM(X)\subset \DTM(X)\\
&W^{>n}:\DTM(X)\to  W^{>n}\DTM(X)\subset\DTM(X).
\end{align*}
Then
\begin{enumerate}
\item[(a)] $W_n$ and   $W^{>n}$ are exact
\item[(b)] $W_n$ is right adjoint to the inclusion $W_n\DTM(X)\to \DTM(X)$ and   $W^{>n}$ is left adjoint to the inclusion   $W^{>n}\DTM(X)\to \DTM(X)$.
\item[(c)] For each $n<m$ there is an exact functor 
\[
W_{[n+1,m]}:\DTM(X)\to W_{[n+1,m]}\DTM(X)\subset\DTM(X)
\]
and a natural distinguished triangle
\[
W_n\to W_m\to W_{[n+1,m]}\to W_n[1].
\]
\item[(d)] $\DTM(X)=\cup_{n\in\Z}W_n\DTM(X)=\cup_{n\in\Z}W^{>n}\DTM(X)$.
\end{enumerate}
\end{theorem}

\begin{proof} By lemma~\ref{lem:TateMotCoh}, we have an isomorphism
\begin{align*}
\Hom_{\DM(X)_\Q}(\Q_X(a),\Q_X(b)[m])&\cong H^m(X,\Q(b-a))\\
&= \begin{cases}0&\text{ for }b<a\\
0&\text{ for }b=a, m\neq0\\
\Q\cdot\id&\text{ for }b=a, m=0.\end{cases}
\end{align*}
Thus,  \cite[lemma~1.2]{LevineTate} 
applies to prove the theorem.
\end{proof}

We denote the exact functor $W_{[n,n]}:\DTM(X)\to W_{[n,n]}\DTM(X)$ by $\gr^W_n$ and the category 
$W_{[n,n]}\DTM(X)$ by $\gr^W_n\DTM(X)$.
\begin{remark} Since 
\[
\Hom_{\DTM(X)}(\Q_X(-n),\Q_X(-n)[m])=\begin{cases} 0&\text{ for }m\neq0\\
\Q\cdot\id&\text{ for }m=0,\end{cases}
\]
the category $\gr^W_n\DTM(X)$ is equivalent to $D^b(\Q)$. Thus, we can define the $\Q$-vector space  $H^n(\gr^W_nM)$ for $M$ in $\DTM(X)$.
\end{remark}

\begin{definition} 1. We say that $X$ {\em satisfies the Beilinson-Soul\'e vanishing conjectures} if $H^m(X,\Q(n))=0$ for $m\le 0$ and $n\neq0$.\\
\\
2. Let $\DTM(X)^{\le0}$ be the full subcategory of $\DTM(X)$ with objects $M$ such that $H^m(\gr^W_nM)=0$ for all $m>0$ and all $n\in\Z$. Let $\DTM(X)^{\ge0}$ be the full subcategory of $\DTM(X)$ with objects $M$ such that $H^m(\gr^W_nM)=0$ for all $m<0$ and all $n\in\Z$. Let $\MT(X):=\DTM(X)^{\le0}\cap \DTM(X)^{\ge0}$.
\end{definition}

\begin{theorem}\label{thm:MixTateTStruct} Suppose $X$ satisfies the Beilinson-Soul\'e vanishing conjectures. Then \\
\\
1. $(\DTM(X)^{\le0}, \DTM(X)^{\ge 0})$ is a non-degenerate $t$-structure on $\DTM(X)$ with heart $\MT(X)$ containing the Tate motives $\Q_X(n)$, $n\in\Z$. \\
\\
2. $\MT(X)$ is equal to its smallest abelian subcategory containing the $\Q_X(n)$, $n\in\Z$, and closed under extensions in $\MT(X)$.\\
\\
3. The tensor operation in $\DTM(X)$ restricted to $\MT(X)$ makes $\MT(X)$ a rigid $\Q$-linear abelian tensor category.\\
\\
4. The functor $\oplus_n\gr^W_n:\MT(X)\to\Vect_\Q$ is a fiber functor, making $\MT(X)$ a neutral Tannakian category.
\end{theorem}

\begin{proof} By lemma~\ref{lem:TateMotCoh}, the assumption that $X$  satisfies the Beilinson-Soul\'e vanishing conjectures implies that
\[
\Hom_{\DTM(X)_\Q}(\Q_X(a),\Q_X(b)[m])=  \begin{cases}0&\text{ for }b>a, m\le0\\
0&\text{ for }b=a, m\neq0\end{cases}
\]
With this, the result follows from \cite[theorem~1.4, proposition~2.1]{LevineTate}. 
\end{proof}

\section{Cycle algebras} \label{sec:CycAlg}
Bloch's cycle complex $z^p(X,*)$ is defined using cycles on $X\times\Delta^n$, where $\Delta^n$ is the algebraic $n$-simplex
\[
\Delta^n:=\Spec k[t_0,\ldots, t_n]/(\sum_it_i-1).
\]
One can also use cubes instead of simplices to define the various versions of the cycle complexes. The major advantage is that the product structure for the cubical complexes is easier to define and, with $\Q$-coefficients, one can construct cycle complexes which have a strictly commutative and associative product. This approach is used by Hanamura in his construction of a category of mixed motives, as well as in the construction of categories of Tate motives by  Bloch \cite{BlochLieAlg}, Bloch-Kriz \cite{BlochKriz} and Kriz-May \cite{KrizMay}. 

We combine the cubical variation with the strictly functorial constructions of Friedlander-Suslin-Voevodsky to give a functorial version of the cycle complex. This allows us to extend the representation theorem of Spitzweck to give a description of mixed Tate motives over a smooth base in terms of cell modules over a cycle algebra.
 
We give a rather sketchy treatment in this section. The second author will make a complete treatment in his forthcoming article \cite{LevineTateMotive}. In this section, we assume that $k$ admits resolution of singularities.

\begin{SecRemark} Joshua \cite{Joshua} has given a definition of a cycle algebra $\sA_X$ over a smooth quasi-projective variety $X$, and has defined the triangulated category of mixed Tate motives over $X$ as $\sD^f_{\sA_X}$, along the lines outlined above. However, there is apparantly an error in his construction, in that he uses the cubical Suslin complex construction on the sheaves $\Z^{tr}(\A^n)$ in his definition of the cycle complex $\sA_X$. Since the projection $\A^n\to \Spec k$ induces a weak equivalence $C_*^\Sus(\Z^{tr}(\A^n))\to C_*^\Sus(\Z)$, the degree $n$ portion of Joshua's cycle algebra does not represent weight $n$ motivic cohomology, as is claimed in \cite{Joshua}.
\end{SecRemark}

\subsection{Cubical complexes}\label{subsection:cubes} 
We  recall the definition of the cubical version of the Suslin-complex $C_*$. 

Let $(\sq^1,\del\sq^1)$ denote the pair $(\A^1,\{0,1\})$, and $(\sq^n,\del\sq^n)$  the $n$-fold product of $(\sq^1,\del\sq^1)$. Explicitly, $\sq^n=\A^n$, and $\del\sq^n$ is the divisor $\sum_{i=1}^n(x_i=0)+\sum_{i=1}^n(x_i=1)$, where $x_1,\ldots, x_n$ are the standard coordinates on $\A^n$.
 A {\em face} of $\sq^n$ is a face of the normal crossing divisor $\del\sq^n$, i.e., a subscheme defined by equations of the form $t_{i_1}=\epsilon_1,\ldots, t_{i_s}=\epsilon_s$, with the $\epsilon_j$ in $\{0,1\}$. If a face $F$ has codimension $m$ in $\sq^n$, we write 
 $\dim F=n-m$.

For $\epsilon\in\{0,1\}$ and $j\in\{1,\ldots, n\}$ we let $\iota_{j,\epsilon}:\sq^{n-1}\to\sq^n$ be the closed embedding defined by inserting an $\epsilon$ in the $j$th coordinate. We let $\pi_j:\sq^n\to\sq^{n-1}$ be the projection which omits the $j$th factor.

\begin{definition}\label{defn:cubical} Let $X$ be a noetherian scheme and let $\sF$ be presheaf on $\Sm/X$. Let $C_n^\cb(\sF)$ be the presheaf
\[
C_n^\cb(\sF)(X):=\sF(X\times\sq^n)/\sum_{j=1}^n\pi_j^*(\sF(X\times\sq^{n-1})),
\]
and let $C_*^\cb(\sF)$ be the complex with differential 
\[
d_n=\sum_{j=1}^n(-1)^{j-1}F(\iota_{j,1})-\sum_{j=1}^n(-1)^{j-1}F(\iota_{j,0}).
\]
\end{definition}
If  $\sF$ is a Nisnevich sheaf, then $C_*^\cb(\sF)$ is a complex of Nisnevich sheaves, and if $\sF$ is a presheaf (resp. Nisnevich sheaf) with transfers, then $C_*^\cb(\sF)$ is a complex of presheaves (resp. Nisnevich sheaves) with transfers. We extend the construction to 
complexes of sheaves (with transfers) by taking the total complex of the evident double complex.

For a presheaf $\sF$ on $\Sm/X$ and $Y\in \Sm/X$, let 
\[
C_n^\Alt(\sF)(Y)\subset C_n^\cb(\sF)(Y)_\Q=\sF(Y\times\sq^n)_\Q
\]
be the $\Q$-subspace consisting of the alternating elements of $\sF(Y\times\sq^n)_\Q$ with respect to the action of the symmetric group $\Sigma_n$ on $\sq^n$, i.e., the elements $x$ satisfying
\[
(\id\times\sigma)^*(x)=\sgn(\sigma)\cdot x.
\]
Here $\Sigma_n$ acts on $\sq^n=\A^n$ by permuting the coordinates.
$Y\mapsto C_n^\Alt(\sF)(Y)$ evidently forms a sub-presheaf of $ C_n^\cb(\sF)_\Q$, which we denote by $C_n^\Alt(\sF)$; in fact the $C_n^\Alt(\sF)$ form a  subcomplex $C_*^\Alt(\sF)\subset C_*^\cb(\sF)_\Q$. We extend this to  
complexes of presheaves by taking the total complex of the evident double complex.

The arguments of  e.g. \cite[section~2.5]{LevineHandbook}  show
\begin{lemma}\label{lem:OmnibusCube} Let $\sF$ be a 
complex of presheaves on  $\Sm/X$.
\begin{enumerate}
\item There is a natural isomorphism $C_*^\Sus(\sF)\cong C_*^\cb(\sF)$ in the derived category of presheaves on $\Sm/X$. If $\sF$ is a complex of presheaves with transfer, we have an isomorphism $C_*^\Sus(\sF)\cong C_*^\cb(\sF)$ in the derived category  $D(\PST(X))$.
\item The inclusion $C_*^\Alt(\sF)(Y)\subset C_*^\cb(\sF)_\Q(Y)$ is a quasi-iso\-morph\-ism for all $Y\in\Sm/X$.
\end{enumerate}
\end{lemma}

\begin{remark}\label{rem:BlochCubeAlt} One can define a cubical version of Bloch's cycle complex, following the pattern of definition~\ref{defn:cubical}. That is, define $\un{z}^q(X,n)^\cb$ to be the free abelian group on the codimension $q$ subvarieties $W\subset  X\times\sq^n$ such that $W\cap X\times F$ has codimension $q$ for every face $F\subset \sq^n$, and let $z^q(X,n)^\cb$ be the quotient of  $\un{z}^q(X,n)^\cb$ by the  ``degenerate" cycles coming from  $\un{z}^q(X,n-1)^\cb$ by pull-back. This gives us the complex $z^q(X,*)^\cb$, which is quasi-isomorphic to the simplicial version $z^q(X,*)$ defined in \cite{AlgCyc}.

Taking the subgroups of alternating cycles gives us the subcomplex $z^q(X,*)^\Alt\subset z^q(X,*)^\cb_\Q$, quasi-isomorphic to $z^q(X,*)^\cb_\Q$.
\end{remark}

\begin{lemma} Let $\sF$ be in $C(\PST(X))$. Suppose that $C_*^\cb(\sF)$ satisfies Nisnevich excision. Then $C_*^\cb(\sF)$ is quasi-fibrant for the Nisnevich local model structure on 
$C(\PST(X))$ and is  $\A^1$-local.
\end{lemma}

\begin{proof} Let $C_*^\cb(\sF)\to C_*^\cb(\sF)^f$ be a fibrant model for $C_*^\cb(\sF)$, with respect to the Nisnevich local model structure on $C(\PST(X))$. Since $C_*^\cb(\sF)$ satisfies Nisnevich excision, the map of complexes
\[
C_*^\cb(\sF)(Y)\to C_*^\cb(\sF)^f(Y)
\]
is a quasi-isomorphism for every $Y\in \Sm/X$. Thus, since the homotopy category of  $C(\PST(X))$ for the
Nisnevich local model structure is equivalent to the derived category of $\Sh^{tr}_\Nis(X)$,  we have isomorphisms for  every $Y\in \Sm/X$ and $n\in\Z$:
\begin{align*}
\Hom_{D(\Sh^{tr}_\Nis(X))}&(\Z_X^{tr}(Y), C_*^\cb(\sF)[n])\\
&\cong 
\Hom_{D(\Sh^{tr}_\Nis(X))}(\Z_X^{tr}(Y), C_*^\cb(\sF)^f[n])\\
&\cong \Hom_{K(\Sh^{tr}_\Nis(X))}(\Z_X^{tr}(Y), C_*^\cb(\sF)^f[n])\\
&\cong H^n(C_*^\cb(\sF)^f(Y))\\
&\cong  H^n(C_*^\cb(\sF)(Y)).
\end{align*}
Since the presheaves  $\Z_X^{tr}(Y)$ form a set of cofibrant generators for $C(\PST(X))$,  it follows that $C_*^\cb(\sF)$ is quasi-fibrant for the Nisnevich local model structure on 
$C(\PST(X))$.

On the other hand, for every $\sF$,  the cubical Suslin complex construction $C_*^\cb(\sF)$ is homotopy invariant as a complex of presheaves, i.e.,
\[
C_*^\cb(\sF)(Y)\to C_*^\cb(\sF)(Y\times\A^1)
\]
is a quasi-isomorphism for each $Y\in \Sm/X$. Thus
\begin{multline*}
\Hom_{D(\Sh^{tr}_\Nis(X))}(\Z_X^{tr}(Y), C_*^\cb(\sF)[n])\\
\to
\Hom_{D(\Sh^{tr}_\Nis(X))}(\Z_X^{tr}(Y\times \A^1), C_*^\cb(\sF)[n])
\end{multline*}
is an isomorphism for all $Y\in \Sm/X$, i.e., $C_*^\cb(\sF)$ is $\A^1$-local.
\end{proof}

\begin{example} Let $W$ be a finite type $k$-scheme. We recall the presheaf with transfers $z_\qfin(W)$ (also denoted $z_\equi(W,0)$ in \cite{FSV}) on $\Sm/k$.   For  $Y\in \Sm/k$,  $z_\qfin(W)(Y)$ is defined to be the free abelian group on integral closed subschemes $Z\subset Y\times_kW$ such that $Z\to Y$ is quasi-finite, and dominant over a component of $Y$.

It follows from \cite[chapter~V, theorem~4.2.2(4)]{FSV} and lemma~\ref{lem:OmnibusCube} that one has a natural isomorphism for $Y\in\Sm/k$
\[
H_n(C^\cb_*(z_\qfin(W))(Y))\cong \H^{-n}_\Nis(Y,C^\cb_*(z_\qfin(W))),
\]
and hence   $C^\cb_*(z_\qfin(W))$ satisfies Nisnevich excision as a complex of presheaves on $\Sm/X$. Thus $C^\cb_*(z_\qfin(W))$ is $\A^1$-local in $C(\PST_{\A^1}(X))$.
\end{example}

Denote by $\Z^{tr}_X(\P^1/\infty)$ the sheaf defined by the exactness of the split exact sequence
\[
0\to\Z^{tr}_X\xrightarrow{i_{\infty*}}\Z^{tr}_X(\P^1)\to \Z^{tr}_X(\P^1/\infty)\to0
\]
Of course, $\Z^{tr}_X(\P^1/\infty)=\Z^{tr}_X(1)[2]$. Similarly, let 
$\Z^{tr}_X((\P^1/\infty)^r)$ be defined by the exactness of
\[
\oplus_{j=1}^r\Z^{tr}_X((\P^1)^{r-1})\xrightarrow{\sum_j i_{j,\infty*}}\Z^{tr}((\P^1)^r)\to
\Z^{tr}_X((\P^1/\infty)^r)\to0
\]
where $i_{j,\infty}:(\P^1)^{r-1}\to(\P^1)^r$ inserts $\infty$ in the $j$th spot. Thus  $\Z^{tr}_X((\P^1/\infty)^r)$ is isomorphic to $\Z^{tr}_X(r)[2r]$.

\begin{remark} We used the notation $T^{tr}$ for $\Z^{tr}_X(\P^1/\infty)$ in the context of ``Tate spectra" (definition~\ref{def:TateObj}); we introduce this new notation to make clear the relation with the sheaf $z_\qfin(\A^1)$.
\end{remark}

\subsection{The cycle cdga in $\DM^\eff(X)_\Q$} 
For $Y\in\Sm/k$, we denote $\Z^{tr}_{\Spec k}(Y)$ by $\Z^{tr}(Y)$.

The symmetric group $\Sigma_q$ acts on  $\Z^{tr}((\P^1/\infty)^q)$ by permuting the coordinates in $(\P^1)^q$. We let 
$\sN^\efc(q)\subset C_*^\Alt(\Z^{tr}((\P^1/\infty)^q)$ be the subsheaf of {\em symmetric} sections with respect to this action. 

\begin{lemma}\label{lem:symComp} The inclusion $\sN^\efc(q)\subset C_*^\Alt(\Z^{tr}((\P^1/\infty)^q)$  is an isomorphism in $\DM^\eff_-(k)$, in fact a quasi-isomorphism of complexes of presheaves on $\Sm/k$.
\end{lemma}
\begin{proof} Fix $Y\in \Sm/k$. We have the sequence of maps
\[
C_*(\Z^{tr}((\P^1/\infty)^q))(Y)\to C_*(z_\qfin(\A^q))(Y)\to z^q(Y\times\A^q,*),
\]
the first map induced by the inclusion $\A^q\subset (\P^1)^q$, the second by the inclusion of the quasi-finite cycles on $Y\times\A^q\times\Delta^n$ to the cycles in good position on $Y\times\A^q\times\Delta^n$. Both maps are quasi-isomorphisms: for the first, use the localization sequence of
 \cite[chapter~IV, corollary~5.12]{FSV} together with \cite[chapter~IV, theorem~8.1]{FSV}; for the second, 
 use the duality theorem \cite[chapter~IV, theorem~7.4]{FSV} and Suslin's comparison theorem \cite[chapter~VI, theorem~3.1]{FSV}. 

Passing to the cubical versions, tensoring with $\Q$ and taking the alternating subcomplexes, it follows from lemma~\ref{lem:OmnibusCube}  and remark~\ref{rem:BlochCubeAlt} that we have the sequence of quasi-isomorphisms
\[
C^\Alt_*(\Z^{tr}((\P^1/\infty)^q))(Y)\to C^\Alt_*(z_\qfin(\A^q))(Y)\to z^q(Y\times\A^q,*)^\Alt.
\]
As the pull-back  by the projection $p:Y\times\A^q\to Y$
\[
z^q(Y,*)^\Alt\to z^q(Y\times\A^q,*)^\Alt
\]
is also a quasi-isomorphism by the homotopy property, $\Sigma_q$ acts trivially on $z^q(Y\times\A^q,*)^\Alt$, in $D^-(\Ab)$, and thus  $\Sigma_q$ acts trivially on the cohomology of the complex  $C^\Alt_*(\Z^{tr}((\P^1/\infty)^q))(Y)$. Since $C^\Alt_*(\Z^{tr}((\P^1/\infty)^q))(Y)$ is a complex of $\Q$-vector spaces, it follows that 
$\sN^\efc(q)(Y)\to C_*^\Alt(\Z^{tr}((\P^1/\infty)^q)(Y)$  is a quasi-isomorphism, as claimed.
\end{proof}

For $X,Y\in\Sm/k$, the external product of correspondences gives the associative external product
\begin{multline*}
C_n^\cb(\Z^{tr}((\P^1/\infty)^q)(X)\otimes C_m^\cb(\Z^{tr}((\P^1/\infty)^p))(Y)\\
\to
C_{n+m}^\cb(\Z^{tr}((\P^1/\infty)^{p+q}))(X\times_kY).
\end{multline*}
Taking $X=Y$ and pulling back by the diagonal $X\to X\times_kX$ gives the cup product of complexes of sheaves
\[
\cup:C_*^\cb(\Z^{tr}((\P^1/\infty)^p))\otimes C_*^\cb(\Z^{tr}((\P^1/\infty)^q))\to
C_*^\cb(\Z^{tr}((\P^1/\infty)^{p+q})).
\]
Taking the alternating projection with respect to  the $\sq^*$ and symmetric projection with respect to the $\A^*$ yields the associative, commutative product
\[
\cdot:\sN^\efc(p)\otimes \sN^\efc(q)\to \sN^\efc(p+q),
\]
which makes $\sN^\efc:=\Q\oplus \oplus_{r\ge1}\sN^\efc(r)$ into an Adams graded cdga object in $C(\Sh_\Nis^{tr}(k))$. By abuse of notation, we write $\sN^\efc(0)$ for the constant presheaf $\Q$.

\begin{definition}\label{def:CycleComplex}
For $X\in\Sm/k$, we let $\sN_X^\efc(q)$ denote the restriction of $\sN^\efc(q)$ to $\SmCor(X)$; similarly define the  Adams graded cdga object in $C(\Sh^{tr}_\Nis(X))$:
\[
\sN^\efc_X=\Q\oplus \oplus_{q\ge1}\sN^\efc_X(q).
\]
\end{definition}

Taking sections of $\sN^\efc$ on $X$ gives us  the Adams graded cdga $\sN^\efc(X)$. In fact, $\sN_X^\efc$ is a presheaf of Adams graded cdgas over $\sN^\efc(X)$, where for $f:Y\to X$ in $\Sm/X$, the algebra structure comes from the pull-back map
\[
f^*:\sN^\efc(X)\to \sN_X^\efc(Y)=\sN^\efc(Y).
\]
\begin{lemma}\label{lem:ProductIso} The multiplication map
\[
\sN^\efc_X(r)\otimes^{tr}_X \sN^\efc_X(s)\to \sN^\efc_X(r+s)
\]
is an isomorphism in $\DM^\eff(X)_\Q$.
\end{lemma}

\begin{proof} First note that the restriction of $\Z^{tr}((\P^1/\infty)^n)$ to $\SmCor(X)$ is canonically isomorphic to $\Z^{tr}_X((\P^1/\infty)^n)$, where the isomorphism is induced by the natural isomorphisms
\[
Y\times_X(\P^1_X\times_X\ldots\times_X\P^1_X)\cong Y\times_X((\P^1)^n\times_kX)\cong Y\times_k(\P^1)^n
\]
for $Y\in \Sm/X$.   In $\DM^\eff(X)$ we have the canonical isomorphism
\[
m_X(X)(n)[2n]\to C^\cb_*(\Z^{tr}_X((\P^1/\infty)^n))
\]
which identifies the multiplication
\[
C^\cb_*(\Z^{tr}_X((\P^1/\infty)^r))\otimes^{tr}_XC^\cb_*(\Z^{tr}_X((\P^1/\infty)^s))\to C^\cb_*(\Z^{tr}_X((\P^1/\infty)^{r+s}))
\]
with the canonical isomorphism
\[
m_X(X)(r)[2r]\otimes^{tr}_X m_X(X)(s)[2s]\to m_X(X)(r+s)[2r+2s].
\]
\end{proof}

\section{$\sN^\efc(X)$-modules and motives} \label{sec:NModMot} We relate the category of Tate motives to the dg modules over the cycle algebra $\sN^\efc(X)$. We give a rather sketchy treatment in this section, a detailed account will be written in the forthcoming article 
 \cite{LevineTateMotive} by the second author. In this section, we assume $k$ admits resolution of singularities.

\subsection{The contravariant motive}

We define the functor
\[
h_X:\Sm/X^\op\to \DM(X)
\]
as follows: For $Y\to X$ in $\Sm/X$ we have the internal Hom presheaf on $\SmCor(X)$ defined by
\[
\sHom(\Z^{tr}_X(Y),C_*(\Z^{tr}_X(n)[2n]))(W):= C_*(\Z^{tr}_X(n)[2n])(Y\times_XW).
\]
The multiplication
\[
\Z_X^{tr}(n)[2n]\otimes^{tr}_X\Z_X^{tr}(1)[2]\to \Z^{tr}_X(n+1)[2n+2]
\]
gives rise to the bonding maps
\[
\sHom(\Z^{tr}_X(Y),C_*(\Z^{tr}_X(n)[2n]))\otimes^{tr}_XT^{tr}\to 
\sHom(\Z^{tr}_X(Y),C_*(\Z^{tr}_X(n+1)[2n+2]))
\]
defining the $T^{tr}$ spectrum $h_X(Y)\in\Spt_{T^{tr}}(X)$:
\[
h_X(Y):=(\sHom(\Z^{tr}_X(Y),C_*(\Z^{tr}_X)),\ldots, \sHom(\Z^{tr}_X(Y),C_*(\Z^{tr}_X(n)[2n])),\ldots).
\]

Using the action of correspondences on $\Z^{tr}_X(Y)$, one sees immediately that $h_X$ extends to a  functor
\[
h_X:\SmCor(X)^\op\to \Spt_{T^{tr}}(X),
\]
which in turn extends to 
\[
C(h_X):C(\SmCor(X)^\op)\to \Spt_{T^{tr}}(X).
\]

\begin{lemma} The composition of $h_k$ with the canonical localization functor
\[
\Spt_{T^{tr}}(k)\to \DM(k)
\]
gives a dual to the composition
\[
\SmCor(k)\xrightarrow{M_\gm}\DM^\eff(k)\xrightarrow{\Sigma^\infty_{T^{tr}}}\DM(k).
\]
\end{lemma}

\begin{proof} For $Y\in\Sm/k$, we denote $C_*(z_\qfin(Y))$ by $C_*^c(Y)$ and let $M_\gm^c(Y)$ denote the image of $C_*^c(Y)$ in $\DM^\eff_-(k)$. We recall the presheaf with transfers $z_\equi(Y,r)$, with $z_\equi(Y,r)(X)$ the free abelian group on the subvarieties $W\subset X\times Y$ such that $W\to X$ is dominant and equi-dimensional of relative dimension $r$ over some component of $X$.

For $Y\in\Sm/k$ of dimension $d$, one has the dual motive $M_\gm(Y)^\vee$ in $\DM_\gm(k)$, since $k$ admits resolution of singularities.
Also, $M_\gm(Y)^\vee(d)[2d]$ is in $\DM^\eff_\gm(k)$ and the image of  
$M_\gm(Y)^\vee(d)[2d]$ in $\DM^\eff_-(k)$ is canonically isomorphic to $M_\gm^c(Y)$ (see \cite[chapter~V, section~4.3]{FSV}). Letting $\Sigma^\infty_{T^{tr}}M_\gm^c(Y)(-d)[-2d]$ denote the $T^{tr}$ spectrum
\[
(0,\ldots,0, C_*^c(Y), C_*^c(Y)(1)[2],\ldots)
\]
with $C_*^c(Y)$ in degree $d$, we see that in $\DM(k)$, $\Sigma^\infty_{T^{tr}}M_\gm(Y)$ has a dual, namely, the object represented by  $\Sigma^\infty_{T^{tr}}M_\gm^c(Y)(-d)[-2d]$. 

The restriction by the open immersion $\A^n\to (\P^1)^n$ induces a quasi-isomorphism of presheaves
\[
\sHom(\Z^{tr}_k(Y),C_*(\Z^{tr}_k(n)[2n]))\to \sHom(\Z^{tr}_X(Y),C^c_*(\A^n)).
\]
By the duality theorem \cite[chapter~IV, theorem~7.1]{FSV}, the inclusion of complexes of presheaves
\[
\sHom(\Z^{tr}(Y),C^c_*(\A^n))\to C_*(z_\equi(Y\times\A^n,d))
\]
is a quasi-isomorphism of complexes of presheaves, as is each morphism in the sequence
\[
C^c_*(Y\times \A^{n-d})\to \sHom(\Z^{tr}(\A^d),C^c_*(Y\times \A^{n-d}))\to C_*(z_\equi(Y\times\A^n,d))
\]
for all $n\ge d$. 

By \cite[chapter~V, corollary~4.1.8]{FSV} we have $M^c_\gm(Y\times\A^{n})\cong M^c_\gm(Y)(n)[2n]$ for all $n\ge0$ Thus we have the canonical isomorphisms in $\DM^\eff_-(k)$:
\[
C_*^c(Y)(n-d)[2n-2d]\cong C^c(Y\times\A^{n-d}) 
\cong \sHom(\Z^{tr}_k(Y),C_*(\Z^{tr}_k(n)[2n]))
\]
for all $n\ge d$. One checks that this isomorphism is compatible with the bonding morphisms for $\Sigma_{T^{tr}}^\infty M_\gm^c(Y)(-d)[-2d]$ and $h_k(Y)$, giving the desired isomorophism
$M_\gm(Y)^\vee\cong h_k(Y)$ in $\DM(k)$. 
\end{proof}

We let
\[
h_X:K(\SmCor(X)^\op)\to \DM(X)
\]
be the exact functor induced by the composition
\[
C(\SmCor(X)^\op)\xrightarrow{C(h_X)}\Spt_{T^{tr}}(X)\to \DM(X).
\]

We can use the cycle complex construction $\sN^\efc_X$ (definition~\ref{def:CycleComplex}) to define a $\Q$ version of $h_X$. Indeed, for $Y\in\Sm/X$, set
\[
\nh_X(Y)(n):=\sHom(\Q^{tr}_X(Y), \sN^\efc_X(n)).
\]
The composition
\[
\Z^{tr}_X(1)[2]\to C^\cb_*(\Z_X^{tr}(1)[2])\to \sN^\efc_X(1)
\]
together with the multiplication in $\sN^\efc_X$ induces bonding maps
\[
\epsilon_n:\sHom(\Q^{tr}_X(Y), \sN^\efc_X(n))\otimes^{tr}_X T^{tr}_X\to \sHom(\Q^{tr}_X(Y), \sN^\efc_X(n+1)),
\]
giving us the $T^{tr}$ spectrum
\[
\nh_X(Y):=(\sHom(\Q^{tr}_X(Y), \sN^\efc_X(0)),\sHom(\Q^{tr}_X(Y), \sN^\efc_X(1)),\ldots).
\]
Sending $Y$ to $\nh_X(Y)$ gives an exact functor
\[
\nh_X:K(\SmCor(X))^\op\to \DM(X)_\Q.
\]

We have the canonical isomorphism in $D(\Q)$
\[
\sN^\efc(n)(Y)\cong C_*(\Z^{tr}_X(n)[2n])(Y)_\Q.
\]
This gives an isomorphism (in $D(\PST(X))_\Q$) 
\[
\sHom(\Q^{tr}_X(Y), \sN^\efc_X(n))\cong\sHom(\Z^{tr}_X(Y),C_*(\Z^{tr}_X(n)[2n]))_\Q,
\]
which induces a canonical isomorphism
\[
\nh_X(Y)\cong h_X(Y)_\Q
\]
natural in $Y$, in fact an isomorphism of functors
\[
\nh_X\cong h_{X\Q}:K(\SmCor(X))^\op\to \DM(X)_\Q.
\]

\subsection{Cell modules as Tate motives}
Recall the Adams graded cdga $\sN^\efc(X)$ gotten by evaluating the presheaf $\sN^\efc_k$ of 
Adams graded cdgas at $X\in\Sm/k$. The following result extends Spitzweck's representation theorem  (see \cite[section~5]{LevineHandbook}) from fields to $X\in \Sm/k$. 

\begin{theorem} \label{thm:CellModMain}  Let $X$ be in $\Sm/k$. There is an exact tensor functor
\[
\sM_X:\sD_{\sN^\efc(X)}\to \DM(X)_\Q
\]
with $\sM_X(\Q(n))\cong \Q_X(n)$. In addition\\
\\
1. The restriction of $\sM_X$ to
\[
\sM_X^f:\sD^f_{\sN^\efc(X)}\to \DM(X)_\Q
\]
defines an equivalence of $\sD^f_{\sN^\efc(X)}$ with $\DTM(X)$, as tensor triangulated categories, natural in $X$.\\
\\
2. $\sM^f$ transforms the weight filtration in $\sD^f_{\sN^\efc(X)}$ to that in $\DTM(X)$.\\
\\
3. Suppose that $X$ satisfies the Beilinson-Soul\'e vanishing conjectures. Then $\sM^f$ is a functor of triangulated categories with $t$-structure. In particular, $\sM^f$ intertwines the respective truncation functors and induces an equivalence of Tannakian categories
\[
H^0(\sM^f):\sH^f_{\sN^\efc(X)}\to \MT(X),
\]
which identifies $\sD^f_{\sN^\efc(X)}$ with $\DTM(X)$. 
\end{theorem}

\begin{proof}[Sketch of proof.] We just describe the construction of the exact functor $\sM_X$.

Let $M=\oplus_rM(r)$ be an $\sN^\efc(X)$-cell module. This gives us the presheaf of Adams graded dg $\sN^\efc_X$-modules $M\otimes_{\sN^\efc(X)}\sN^\efc_X$. Decompose 
$\sM:=M\otimes_{\sN^\efc(X)}\sN^\efc_X$ as the sum of its Adams graded pieces
\[
\sM=\oplus_r \sM(r).
\]
We have the canonical map $T^{tr}\to \sN^\efc_X(1)$; combining with the multiplication
\[
\sM\otimes^{tr}_X\sN^\efc_X\to \sM
\]
gives us the bonding maps $\epsilon_r:\sM(r)\otimes^{tr}_XT^{tr}\to \sM(r+1)$. We set
\[
\sM_X(M):=((\sM(0), \sM(1),\ldots), \epsilon_r)
\]
giving the functor of dg categories
\[
\sM_X:\CM_{\sN^\efc(X)}\to \Spt_{T^{tr}}(X)_\Q.
\]
As a homotopy equivalence of $\sN^\efc(X)$-cell modules clearly gives rise to a weak equivalence of the associated motives, we have the exact functor
\[
\sM_X:\KCM_{\sN^\efc(X)}\to \DM(X)_\Q.
\]
Since $\KCM_{\sN^\efc(X)}\to \sD_{\sN^\efc(X)}$ is an equivalence, we have as well the exact functor
\[
\sM_X:\sD_{\sN^\efc(X)}\to \DM(X)_\Q.
\]
\end{proof}

\subsection{From cycle algebras to motives} Let $p:Y\to X$ be in $\Sm/X$, giving us the map of cycle algebras
\[
p^*:\sN^\efc(X)\to \sN^\efc(Y);
\]
in particular, we may consider $\sN^\efc(Y)$ as a dg module over $\sN^\efc(X)$. We let
\[
\rho_Y:\cycmod_X(Y)\to \sN^\efc(Y)
\]
be a quasi-isomorphism with $\cycmod_X(Y)$ an $\sN^\efc(X)$-cell module. 

We proceed to define a natural transformation
\[
\psi_Y:\sM_X(\cycmod_X(Y))\to \nh_X(Y).
\]
Indeed, recall that $\nh_X(Y)$ is the $T^{tr}$-spectrum
\[
(\sHom(\Q^{tr}_X(Y),\sN^\efc_X(0)),\ldots,\sHom(\Q^{tr}(Y),\sN^\efc_X(r)),\ldots),
\]
with bonding maps induced by the multiplication in $\sN^\efc_X$ and the structure map $T^{tr}\to \sN^\efc_X(1)$, while $\sM_X(\cycmod_X(Y))$ is the $T^{tr}$ spectrum
\[
([\cycmod_X(Y)\otimes_{\sN^\efc(X)}\sN^\efc_X](0),\ldots, [\cycmod_X(Y)\otimes_{\sN^\efc(X)}\sN^\efc_X](r),\ldots)
\]
with bonding maps also given by the multiplication with $T^{tr}\to \sN^\efc_X(1)$.

Now take $W\in\Sm/X$. Then
\[
\sHom(\Q^{tr}(Y),\sN^\efc_X(r))(W):=\sN^\efc(Y\times_XW)(r)
\]
Using the external products in $\sN^\efc$, we thus have the canonical map of Adams graded complexes 
\[
\tilde{\psi}_X(W):\sN^\efc(Y)\otimes_{\sN^\efc(X)}\sN^\efc(W)\to \sN^\efc(Y\times_XW).
\]
The maps $\tilde{\psi}_X(W)$ clearly define a map of Adams graded complexes of presheaves with transfer
\[
\tilde{\psi}_Y:\sN^\efc(Y)\otimes_{\sN^\efc(X)}\sN^\efc\to \oplus_{r\ge0}\sHom(\Q^{tr}_X(Y),\sN^\efc_X(r));
\]
restricting to the component of Adams weight $r$ gives the map of complexes of presheaves with 
transfer
\[
\tilde{\psi}_Y(r):[\sN^\efc(Y)\otimes_{\sN^\efc(X)}\sN^\efc](r)\to \sHom(\Q^{tr}_X(Y),\sN^\efc_X(r)).
\]
It is easy to see that $\tilde{\psi}_Y$ respects the action (on the right) by $\sN^\efc$.

Composing $\tilde{\psi}_Y(r)$ with the structure map
\[
\cycmod_X(Y)\otimes_{\sN^\efc(X)}\sN^\efc\xrightarrow{\rho_Y\otimes\id}
\sN^\efc(Y)\otimes_{\sN^\efc(X)}\sN^\efc
\]
gives us the map
\[
\psi_Y(r):[\cycmod_X(Y)\otimes_{\sN^\efc(X)}\sN^\efc](r)\to \sHom(\Q^{tr}_X(Y),\sN^\efc_X(r)).
\]
also respecting the right  $\sN^\efc$ action. Thus   the maps $\psi_Y(r)$ define a map of $T^{tr}$ spectra
\[
\psi_Y:\sM_X(\cycmod_X(Y)) \to \nh_X(Y)
\]
as desired.

In general, $\psi_Y$ does not define an isomorphism in $\DM(X)_\Q$. In this direction our main result is

\begin{proposition}\label{prop:TateKuenneth} Suppose that $h_X(Y)_\Q$ is in $\DTM(X)$. Then 
\[
\psi_Y:\sM_X(\cycmod_X(Y))\to \nh_X(Y)
\] is an isomorphism.
\end{proposition}

\begin{proof}[Sketch of proof]  
We introduce the intermediate category $\sM_{\sN^\efc_X}$ of presheaves with transfers of Adams graded dg modules over the presheaf (on $\Sm/X$) of cdgas $\sN^\efc_X$, and the subcategory of finite cell modules  $\sC\sM^f_{\sN^\efc_X}$.  The transfers and module structure on $M\in \sM_{\sN^\efc_X}$  are required to be compatible in the following manner: For $Y\to X$ in $\Sm/X$, $W,T\in \Sm/Y$,  $M(W)$ and $M(T)$ are dg $\sN^\efc(Y)$ modules via the structure morphisms $W\to Y$, $T\to Y$. Then for $f\in c_Y(W,T)$,  one has
 \[
 a\cdot f^*(m)=f^*(a\cdot m)
 \]
 for $a\in \sN^\efc(Y)$, $m\in M(T)$.

 A modification of the functor $\sM_X$ defines the exact functor
 \[
 \widetilde\sM_X:\sD_{\sN^\efc_X}\to \DM(X)_\Q.
 \]
 
 For $Y\in \Sm/X$,  we have the object $\sHom(\Q^{tr}_X(Y),\sN_X^\efc)$ in $\sD_{\sN^\efc_X}$ and a natural isomorphism
 \[
\phi_Y:  \widetilde\sM_X(\sHom(\Q^{tr}_X(Y),\sN_X^\efc))\to \nh_X(Y).
 \]
 In addition, we have the base-extension functor
 \[
 -\otimes_{\sN^\efc(X)}^L\sN_X^\efc:\sD_{\sN^\efc(X)}\to \sD_{\sN^\efc_X}
 \]
 and natural isomorphism
 \[
 \sM_X\cong  \widetilde\sM_X\circ( -\otimes_{\sN^\efc(X)}^L\sN_X^\efc).
 \]
Finally, the external products for $\sN_X^\efc$ define a natural map
\[
\psi'_Y:\cycmod_X(Y)\otimes_{\sN^\efc(X)}\sN_X^\efc\to \sHom(\Q^{tr}_X(Y),\sN_X^\efc)
\]
making the diagram
\[
\xymatrix{
\sM_X(\cycmod_X(Y))\ar[r]^{\psi_Y}\ar[d]_{\widetilde\sM_X(\psi_Y')}&\nh(Y)\\
\widetilde\sM_X(\sHom(\Q^{tr}_X(Y),\sN_X^\efc))\ar[ur]_{\phi_Y}}
 \]
 commute. Thus we need only show that $\psi'_Y$ is an isomorphism in $\sD_{\sN^\efc_X}$. 
 
For this, let $\sM:=\oplus_{r\ge0}\sM(r)$ be an Adams graded dg $\sN^\efc_X$-module. Define the Adams graded dg $\sN^\efc(X)$ module $\Gamma(\sM)$ by taking sections on $X$:
\[
\Gamma(\sM)(r):=\sM(r)(X)
\]
The $\sN^\efc_X$-module struction gives a canonical map
\[
\psi'_{\sM}:\Gamma(\sM)\otimes^L_{\sN^\efc(X)}\sN_X^\efc\to \sM
\]
which is just the map $\psi_Y'$ in case $\sM=\sHom(\Q^{tr}_X(Y),\sN_X^\efc)$. 

Using the weight filtration in  $\DMT(X)$, one shows that the fact that $\nh(Y)$ is in $\DMT(X)$ implies that the dg $\sN^\efc_X$-module $\sHom(\Q^{tr}_X(Y),\sN_X^\efc)$ is quasi-isomorphic to a finite $\sN^\efc_X$-cell module.  Thus, we need only show that, if  $\sM$ is a finite $\sN^\efc_X$-cell module, then $\psi'_\sM$ is an isomorphism in $\sD_{\sN^\efc_X}$.

For this, one argues again by induction on the weight filtration, reducing to the case of $\sM$ a twist of $\sN^\efc_X$ (i.e. a pure Tate motive), for which the result is obvious.
\end{proof}

\subsection{The cell  algebra of an $X$-scheme} We now assume that $\sN(X)$ is cohomologically connected.

Let $p:Y\to X$ be in $\Sm/X$ with a section $s:X\to Y$. We thus have the map of cycle algebras
$p^*: \sN^\efc(X)\to \sN^\efc(Y)$ making $\sN^\efc(Y)$ a cdga over $\sN^\efc(X)$ 
with augmentation $s^*: \sN^\efc(Y)\to  \sN^\efc(X)$. Let $\sN^\efc(Y)_X\<\infty\>\to \sN^\efc(Y)$ be the relative minimal model of $\sN^\efc(Y)$ over $\sN^\efc(X)$. In particular, $\sN^\efc(Y)_X\<\infty\>$ is  a cell module over $\sN^\efc(X)$. In addition, the multiplication 
\[
\sN^\efc(Y)_X\<\infty\>\otimes \sN^\efc(Y)_X\<\infty\>\to \sN^\efc(Y)_X\<\infty\>
\]
given by the cdga structure on $\sN^\efc(Y)_X\<\infty\>$ descends to
\[
\mu_Y:\sN^\efc(Y)_X\<\infty\>\otimes_{\sN^\efc(X)} \sN^\efc(Y)_X\<\infty\>\to \sN^\efc(Y)_X\<\infty\>.
\]
\begin{definition} The {\em motivic cell algebra of $Y$}, 
\[
\calg_X(Y)\in\CM_{\sN^\efc(X)}
\]
 is  $\sN^\efc(Y)_X\<\infty\>$, considered as a cell module over $\sN^\efc(X)$.
\end{definition}

The same construction we used to define the map $\sM_X(\cycmod_X(Y))\to \nh_X(Y)$ gives us the map in $\Spt_{T^{tr}}(X)_\Q$
\[
\psi_Y:\sM_X(\calg_X(Y))\to \nh_X(Y).
\]

\begin{lemma}\label{lem:TateKuenneth} Suppose that $h_X(Y)_\Q$ is in $\DTM(X)$ and that $Y$ satisfies the Beilinson-Soul\'e vanishing conjectures. Then 
\[
\psi_Y:\sM_X(\calg_X(Y))\to \nh_X(Y)
\]
 is an isomorphism.
\end{lemma}

\begin{proof} This follows from proposition~\ref{prop:TateKuenneth}, once we know that $\calg_X(Y)\to \sN^\efc(Y)$ is a quasi-isomorphism. 

Recall that the Beilinson-Soul\'e vanishing conjectures for $Y$ are just saying that 
$\sN^\efc(Y)$ is cohomologically connected. Using the section $s:X\to Y$, we see that $X$ also satisfies the Beilinson-Soul\'e vanishing conjectures, hence  $\sN^\efc(X)$ is cohomologically connected. The structure map  $\sN^\efc(Y)_X\<\infty\>\to \sN^\efc(Y)$ is thus a quasi-isomorphism by 
proposition~\ref{prop:RelMin}(3).
\end{proof}

\section{Motivic $\pi_1$}

We can now put  all our constructions together  to give a description of the Deligne-Goncharov motivic $\pi_1$ in terms of a relative bar construction.  In this section, we assume $k$ admits resolution of singularities.

\subsection{Cosimplicial constructions} Fix a base-field $k$. We have the action of finite sets on $\Sch_k$ by
\[
X^S:=\prod_{s\in S}X
\]
for  $X\in\Sch_k$ and $S$ a finite set, where $\prod$ means product over $k$. As this defines a functor
\[
X^?:\Sets_{fin}^\op\to \Sch_k
\]
we have an induced functor (also denoted $X^?$) from simplicial objects in finite sets to cosimplicial schemes.

\begin{examples} 1. We have the simplicial finite set $[0,1]$:
\[
[0,1]([n]):=\Hom_\Ord([n],[1])
\]
giving us the {\em cosimplicial path space of $X$}, $X^{[0,1]}$. The two inclusions $i_0, i_1:[0]\to [1]$ induce the projection
\[
\pi:X^{[0,1]}\to X^{\{0,1\}}.
\]
Explicitly, $X^{\{0,1\}}$ is the constant cosimplicial scheme $X^2$. $X^{[0,1]}$ has $n$-cosimplices $X^{n+2}$ with the $i$th coface map given by 
\[
(t_0,\ldots, t_n)\mapsto (t_0,\ldots, t_i,t_i, t_{i+1},\ldots,t_n)
\]
and the codegeneracies given by projections. The structure morphism $\pi$ is given by the projection $X^{n+2}\to X^2$ on the first and last factor.

2. Fix $k$-points $a,b\in X(k)$, giving the map $i_{b,a}:\Spec k\to X^2$ corresponding to $(b,a)$. The {\em pointed path space} $\sP_{b,a}(X)$ is
\[
\sP_{b,a}(X):=\Spec k\times_{i_{b,a},\pi}X^{[0,1]}.
\]
We write $\sP_a(X)$ for $\sP_{a,a}(X)$.
\end{examples}

\subsection{The ind-motive of a cosimplicial scheme}\label{subsec:IndMot} Let $X^\bullet:\Ord\to \Sm/k$ be a smooth cosimplicial $k$-scheme, $[i]\mapsto X[i]\in \Sm/k$. Following Deligne-Goncharov, we define $h_\gm(X^\bullet)$ as an ind-object in $\DM_\gm(k)$.

Let $\Z\Sm/k$ be the additive category generated by $\Sm/k$,  i.e., objects are denoted $\Z(X)$ for $X\in\Sm/k$, for $X$ irreducible, $\Hom_{\Z\Sm/k}(\Z(X),\Z(Y))$ is the free abelian group on $\Hom_{\Sm/k}(X,Y)$ and disjoint union is direct sum.  Sending $X$ to $M_\gm(X)$ extends to
\[
M_\gm:K^b(\Z\Sm/k)\to \DM_\gm(k);
\]
composing with the duality involution ${}^\vee$ on $\DM_\gm(k)$ gives 
\[
h_\gm:={}^\vee\circ M_\gm:K^b(\Z\Sm/k^\op)\to \DM_\gm(k).
\]
Since $DM_\gm(k)$ is pseudo-abelian, $h_\gm$ extends canonically to 
\[
h_\gm:K^b(\Z\Sm/k^\op)^\natural \to \DM_\gm(k).
\]
 where $(-)^\natural$ means the pseudo-abelian hull.

For each $n$, one has the complex $C^*(\Delta_n,X^\bullet)$ with
\[
C^i(\Delta^n,X^\bullet):=\oplus_{g:[i]\hookrightarrow[n]}\Z(X([i])),
\]
where
 the sum is over all injective maps $g:[i]\to[n]$ in $\Ord$.  The boundary
\[
d^i:C^i(\Delta^n,X^\bullet)\to C^{i+1}(\Delta^n,X^\bullet)
\]
is defined as follows: For $0\le j\le i+1$, we have the coface map $\delta^i_j:[i]\to [i+1]$ (see section~\ref{subsection:bar}).   Fix an  injection $g:[i+1]\to [n]$. Define 
\[
\delta^{i,g}_{j*}:C^i(\Delta^n,X^\bullet)\to C^{i+1}(\Delta^n,X^\bullet)
\]
by projecting $C^i(\Delta^n,X^\bullet)$ to the component $\Z(X[i])$ indexed by $g\circ\delta^i_j$ followed by the map
\[
X(\delta^i_j):X[i]\to X[i+1]
\]
and then the inclusion $\Z(X[i+1])\to C^{i+1}(\Delta^n,X^\bullet)$ indexed by $g$. Set
\[
d^i:=\sum_{j,g}\sgn(j,g)\cdot \delta^{i,g}_{j*}
\] where $\sgn(j,g)$ is the sign of the shuffle permutation of $[n]$ given by $(g\circ\delta^i_j([i])^c,g\circ\delta^i_j([i]))$.

Projecting on the factors $g$ with $0$ in the image of $g$ defines a map of complexes
\[
\pi_{n+1,n}:C^*(\Delta_{n+1},X^\bullet)\to C^*(\Delta_n,X^\bullet)
\]
giving us a projective system in $C^b(\Z\Sm/k)$. Reindexing so that $C^n$ is now in degree $-n$ gives an inductive system in $C^b(\Z\Sm/k^\op)$
\[
\ldots\to C_*(\Delta_n,X^\bullet)\to C_*(\Delta_{n+1},X^\bullet)\to \ldots
\]

\begin{definition} $h_\gm(X^\bullet)$ is the ind-object of $\DM_\gm(k)$ defined by the ind-system
\[
n\mapsto h_\gm(C_*(\Delta_n,X^\bullet))
\]
\end{definition}

Let $\Z^*(X^\bullet)$ be the complex which is $X[n]$ in degree $n$, and differential the alternating sum of the coface maps, and let $\Z_*(X^\bullet)$ be the associated complex in $\Z\Sm/k^\op$. Taking the sum of the identity maps defines a map
\[
q_n:C_*(\Delta_n,X^\bullet)\to \Z_*(X^\bullet) 
\]
in $C^-(\Z\Sm/k^\op)$, giving a map of the above ind-system to $\Z_*(X^\bullet)$.

\begin{lemma}[\cite{LevineHandbook}]\label{lem:Ind} Let $F:\Z\Sm/k^\op\to \sA$ be an additive functor to a pseudo-abelian category, closed under filtered inductive limits. Then 
\[
\colim_nF(C_*(\Delta_n,X^\bullet))\to F(\Z_*(X^\bullet))
\]
is a homotopy equivalence in $C^-(\sA)$.
\end{lemma}

The category $\DM(k)$ is large enough to define the object $h(X^\bullet)$ directly.

\begin{definition} For a cosimplicial scheme $X^\bullet$, define $h_k(X^\bullet)$ by
\[
h_k(X^\bullet):=h_k(\Z_*(X^\bullet)). 
\]
Sending $X^\bullet$ to $h_k(X^\bullet)$ extends to a functor
\[
h_k:[\Sm/k^\Ord]^\op\to \DM(k).
\]
\end{definition}

\begin{proposition} We have a natural isomorphism in $\DM(k)$
\[
\colim_n h_\gm(C_*(\Delta_n,X^\bullet))\cong h_k(X^\bullet)
\]
\end{proposition}
\begin{proof}
This follows directly from lemma~\ref{lem:Ind}.
\end{proof}
 
Finally, we may replace $h_\gm$ and $h_k$ with the functor $\nh_k$. Sending $X^\bullet$ to $\nh_k(X^\bullet):=\nh_k(\Z_*(X^\bullet))$ extends to the functor
\[
\nh_k:[\Sm/k^\Ord]^\op\to \DM(k)_\Q,
\]
the natural isomorphism $h_{k\Q}\to \nh_k$ gives natural isomorphisms
\[
\phi_{X^\bullet}:h_k(X^\bullet)_\Q\to \nh_k(X^\bullet)
\]
and we have natural isomorphisms
\[
h_\gm(C_*(\Delta_n,X^\bullet))_\Q\to \nh_k((C_*(\Delta_n,X^\bullet))
\]
and
\[
\colim_n h_\gm(C_*(\Delta_n,X^\bullet))\cong \nh_k(X^\bullet).
\]

\subsection{Motivic $\pi_1$}  Let $X$ be a smooth $k$-scheme with a $k$-point $x\in X(k)$. This gives us the ind-system in $\DM_\gm(k)_\Q$
\[
n\mapsto h_\gm(C_*(\Delta_n,\sP_x(X)))_\Q
\]
as well as the object $h_k(\sP_x)\in  \DM(k)_\Q$ with isomorphism
\[
\colim_nh_\gm(C_*(\Delta_n,\sP_x(X)))\cong h_k(\sP_x).
\]

Suppose that $h_\gm(X)_\Q$ is in $\DTM(k)$. As $\DTM(k)$ is a tensor subcategory of $\DM(k)_\Q$ and as $h_\gm(X^n)=h_\gm(X)^{\otimes n}$, it follows that $h_\gm(X^n)_\Q$ is in $\DTM(k)$ for all $n\ge0$. Since the individual terms in the complex $h_\gm(C_*(\Delta_n,\sP_x(X)))_\Q$ are all direct sums of self-products of $X$, the motive $h_\gm(C_*(\Delta_n,\sP_x(X)))_\Q$ is in $\DTM(X)$ for all $n$.

If  $k$ satisfies the Beilinson-Soul\'e vanishing conjectures, we have the truncation functor 
\[
H^0_{mot}:\DTM(k)\to \MT(k).
\]
Thus we have the ind-system  in $\MT(k)$
\[
n\mapsto H^0_{mot}(h_\gm(C_*(\Delta_n,\sP_x(X)))_\Q):=\chi(X,x)_n.
\]
Deligne-Goncharov \cite{DelGon} 
note that the standard structures of product, coproduct and antipode in the classical bar construction make the ind-system $\chi(X,x)_*$ into an ind-Hopf algebra object in $\MT(k)$. In case $X$ is the complement of a finite set of $k$-points of $\P^1_k$, and $x\in X(k)$,  Deligne and Goncharov define  $\pi_1^{mot}(X,x)$ to be the dual group scheme object in pro-$\MT(k)$.  They generalize the definition of $\pi_1^{mot}(X,x)$ to the case where $X$ is a smooth uni-rational variety defined over $k$: they show in \cite[the\'eor\`eme~4.13]{DelGon} that a suitable object of Deligne's  realization category comes from the mixed Artin-Tate category $MAT(k)$ (which is a bit larger than $MT(k)$ as it takes into account trivial motives defined over a finite extension of $k$). However, they do not give a direct construction as a motive in $DM_{gm}(k)$. 
We extend their definition in the following direction:

\begin{definition}\label{defn:MotPi1} Suppose that $k$ and $X$ both satisfy the Beilinson-Soul\'e vanishing conjectures, and that $\nh_k(X)$ is in $\DMT(k)$. Let $i_x:\Spec k\to X$ be a $k$-point. Define $\pi_1^\mot(X,x)$ to be the group scheme object in pro-$\MT(k)$ dual to the ind-Hopf algebra  $\chi(X,x)_*$.
\end{definition}

\subsection{Simplicial constructions} Let $\sA\xrightarrow{\epsilon}\sN$ be an  augmented cdga over a cdga $\sN$.  Recall from section~\ref{subsec:RelativeBar} the simplicial version of the relative bar construction
 \[
B^\nai_\bullet(\sA/\sN,\epsilon):=\sA^{\otimes_\sN[0,1]}\otimes_{\sA\otimes\sA}\sN.
\]
The total complex associated to the simplicial object $n\mapsto B^\nai_n(\sA/\sN,\epsilon)$ is the relative bar complex $\bar{B}^\nai_\sN(\sA,\epsilon)$.

Using the opposite of the construction described in    section~\ref{subsec:IndMot}, we have the ind-system of ``finite" complexes $C_*(\Delta_n,B^\nai_\bullet(\sA/\sN,\epsilon))$, and a homotopy equivalence
\[
\colim_nC_*(\Delta_n,B^\nai_\bullet(\sA/\sN,\epsilon))\to \bar{B}^\nai_\sN(\sA,\epsilon).
\]
Replacing $\sA$ with its relative minimal model over $\sN$ (assuming for this that $\sN$ is cohomologically connected), we have the refined version of the simplicial bar construction, $B_\bullet(\sA/\sN,\epsilon)$, the associated complex $\bar{B}_\sN(\sA,\epsilon)$, the approximations $C_*(\Delta_n, B_\bullet(\sA/\sN,\epsilon))$ and the homotopy equivalence
\[
\colim_nC_*(\Delta_n,B_\bullet(\sA/\sN,\epsilon))\to \bar{B}_\sN(\sA,\epsilon).
\]

\subsection{The comparison theorem}  \label{subsec:CompThm}

Take $X\in\Sm/k$. with $k$-point $i_x:\Spec k\to X$.   We apply the construction of the preceeding section to the augmented cdga $\sN^\efc(X)$ over $\sN^\efc(k)$:
\[
\xymatrix{
\sN^\efc(X)\ar@<3pt>[r]^{i_x^*}&\sN^\efc(k).\ar@<3pt>[l]^{p^*}
}
\]
Assuming that $\sN(k)$ is cohomolgoically connected, we have  the relative minimal model $\sN_\infty(X/k):=\sN^\efc_k(X)\<\infty\>_{\sN^\efc(k)}$, which is an augmented $\sN^\efc(k)$ algebra via $i_x^*:\sN_\infty(X/k)\to \sN^\efc(k)$. The multiplication in $\sN_\infty(X/k)$ gives the natural maps
\[
\mu_n:\sN_\infty(X/k)^{\otimes_{\sN^\efc(k)} n}\to \sN^\efc(X^n)
\]
which thus gives natural maps in $\DM(k)$
\[
\phi_n(X,x):\sM_k(C_*(\Delta_n, B_\bullet(\sN^\efc(X)/\sN^\efc(k),i_x^*)))\to
\nh_k(C_*(\Delta_n,\sP_x(X)))
\]
and
\[
\phi(X,x):\sM_k(\bar{B}_{\sN^\efc(k)}(\sN^\efc(X),i_x^*))\to
\nh_k(\sP_x(X))).
\]
The maps $\phi_n(X,x)$ give a map of ind-Hopf algebra objects in $\DM(k)$.

\begin{theorem} \label{thm:Main1} Suppose that $\nh(X)$ is in $\DTM(k)$ and $X$ satisfies the Beilinson-Soul\'e vanishing conjectures. Then  both $\phi_n(X,x)$ and $\phi(X,x)$ are isomorphisms in $\DM(k)$.
\end{theorem}

\begin{proof} Note that the Beilinson-Soul\'e vanishing conjectures for $X$ imply the vanishing conjectures for $k$, hence $\sN_(k)$ is cohomologically connected and thus the relative bar complex $\bar{B}_{\sN^\efc(k)}(\sN^\efc(X),i_x^*)$ is defined.

As $\phi(X,x)$ is identified with the filtered homotopy colimit of the maps $\phi_n(X,x)$, it suffices to show that $\phi_n(X,x)$ is an isomorphism for each $n$. But on the individual terms in the complexes defining $C_*(\Delta_n, B_\bullet(\sN^\efc(X)/\sN^\efc(k),i_x^*))$ and $C_*(\Delta_n,\sP_x(X))$, $\phi_n(X,x)$ is the map
\[
\phi_n(X,x)_n:\sM_k(\sN_\infty(X/k)^{\otimes_{\sN^\efc(k)} n})\to\sHom(\Q^{tr}(X^n),\sN^\efc_k) =\nh_k(X^n)
\]
induced by $\psi_k(X^n)\circ\mu_n$.

Since $X$ is a Tate motive and satisfies the Beilinson-Soul\'e vanishing conjectures, it follows from  lemma~\ref{lem:TateKuenneth} that  $\psi_k(X)$ is an isomorphism. Since $\nh(X^n)=\nh(X)^{\otimes n}$ and since $\sM_k$ is a tensor functor (theorem~\ref{thm:CellModMain}), this implies that $\psi_k(X^n)$ is an isomorphism for all $n$.

 In addition,the structure map $\mu_1$ is a quasi-isomorphism since $\sN^\efc(X)$ is cohomologically connected. Since $X$ is a Tate motive, it follows that the motivic cohomology of $X^n$ satisfies a K\"unneth formula for each $n$. Thus $\mu_n$ is a quasi-isomorphism for each $n$, and hence 
$\phi_n(X,x)_n$ is an isomorphism for each $n$.
\end{proof}

\begin{corollary} \label{cor:Main1} Suppose that $\nh_k(X)$ is in $\DTM(k)$ and  $X$ satisfies the Beilinson-Soul\'e vanishing conjectures. Then we have canonical isomorphisms of ind-Hopf algebras in $\MT(k)$, 
\begin{multline*}n\mapsto[
\sM_k(H^0_{\sN^\efc(k)}(C_*(\Delta_n,  B_\bullet(\sN^\efc(X)/\sN^\efc(k),i_x^*))\\\xrightarrow{H^0(\phi_n(X,x))}
H^0_\mot(h_\gm(C_*(\Delta_n,\sP_x(X)))].
\end{multline*}
\end{corollary}

\begin{proof} This follows from theorem~\ref{thm:Main1} and theorem~\ref{thm:CellModMain}.
\end{proof}

\subsection{The fundamental exact sequence}\label{subsec:FundSeq}

Let $X$ be in $\Sm/k$, with structure morphism $p:X\to\Spec k$. We thus have the exact functor of triangulated tensor categories $p^*:\DM(k)\to \DM(X)$; since $p^*(\Z(n))\cong \Z_X(n)$, $p^*$ induces the exact tensor functor
\[
p^*:\DTM(k)\to \DTM(X).
\]
Similarly, if $x\in X(k)$ is a $k$-rational point, we have
\[
i_x^*:\DTM(X)\to \DTM(k).
\]
Both $p^*$ and $i_x^*$ are compatible with the weight filtrations.

Similarly, the maps $p$ and $i_x$ induce maps of cdgas
\[
p^*:\sN^\efc(k)\to \sN^\efc(X);\ i_x^*:\sN^\efc(X)\to \sN^\efc(k)
\]
and thus exact tensor functors
\[
p^*:D^f_{\sN^\efc(k)}\to D^f_{\sN^\efc(X)},\ i_x^*:D^f_{\sN^\efc(X)}\to D^f_{\sN^\efc(k)}.
\]
Recall that the equivalence $\sM_X^f$ of theorem~\ref{thm:CellModMain} is natural in $X$, so we have natural isomorphisms
\[
\sM^f_X\circ p^*\cong p^*\circ \sM^f_k;\ \sM^f_k\circ i_x^*\cong i_x^*\circ\sM^f_X.
\]

Now suppose that  $X$ satisfies the Beilinson-Soul\'e vanishing conjectures; this property is inherited by $k$ using the splitting $i_x^*$. Thus we have this entire structure for the Tannakian categories $\MT(X)$ and $\MT(k)$, with $p^*$ and $i_x^*$ respecting the fiber functors $\gr^W_*$. Similarly, we have functors $p^*$ and $i_x^*$ for the Tannakian categories $\sH^f_{\sN^\efc(X)}$ and $\sH^f_{\sN^\efc(k)}$, respecting the fiber functors $\gr^W_*$. Finally,  $H^0(\sM^f_X)$ and $H^0(\sM^f_k)$ give an equivalence between these two structures.

Let $G(\MT(X),\gr^W_*)$, $G(\MT(k), \gr^W_*)$ denote the Tannaka groups (i.e. pro-group schemes over $\Q$)  of $(\MT(X), \gr^W_*)$ and $(\MT(k),\gr^W_*)$. We sometimes omit the ``base-point" $\gr^W_*$ from the notation.

 The functors $p^*$ and $i_x^*$ gives maps of  pro-group schemes over $\Q$
\[
p_*:G(\MT(X),\gr^W_*)\to G(\MT(k,\gr^W_*)),\ i_{x*}:G(\MT(k),\gr^W_*)\to G(\MT(X),\gr^W_*).
\]
Letting $K=\ker p_*$, we thus have the split exact sequence
\[
\xymatrix{
1\ar[r]&K\ar[r]& G(\MT(X), \gr^W_*)\ar@<3pt>[r]^{p_*}&G(\MT(k), \gr^W_*)\ar@<3pt>[l]^{i_{x*}}\ar[r]&1
}
\]
of pro-group schemes over $\Q$. Via the splitting $i_{x*}$,  $G(\MT(k))$  acts by conjugation on $K$. Thus  the pro-affine Hopf algebra $\Q[K]$ is a 
$G(\MT(k))$-representation. Tannaka duality yields the  corresponding ind object in $\MT(k)$, and its dual is a pro-group scheme object in $\MT(k)$, which we  denote by $K_x$.
As we have seen above, the Deligne-Goncharov motivic  fundamental group $\pi^\mot_1(X,x)$, is also a pro-group  scheme object in $\MT(k)$.

\begin{theorem} \label{thm:Main}
Let $X$ be in $\Sm/k$ with $k$-point $i_x:\Spec k\to X$. Suppose that $X$ satisfies the Beilinson-Soul\'e vanishing conjectures and that the motive $\nh_k(X)\in \DM_\gm(k)_\Q$ is in $\DMT(k)$. Then there is a natural isomorphism
\[
\pi^\mot_1(X,x)\cong K_x
\]
as pro-group   objects in $\MT(k)$.
\end{theorem}

\begin{proof} As we have seen above, we may identify $G(\MT(X))$ and $G(\MT(k))$ with the Tannaka groups of the categories $\sH^f_{\sN^\efc(X)}$ and $\sH^f_{\sN^\efc(k)}$, respectively. By theorem~\ref{thm:MayKriz}, this gives an isomorphism of $K$ with the kernel of the map of pro-groups schemes over $\Q$:
\[
p_*:\Spec(H^0(\bar{B}(\sN^\efc(X))))\to \Spec(H^0(\bar{B}(\sN^\efc(k))))
\]
induced by 
\[
H^0(\bar{B}(p^*)):H^0(\bar{B}(\sN^\efc(k)))\to H^0(\bar{B}(\sN^\efc(X)))
\]
Similarly, the splitting $i_{x*}$ becomes identified with 
\[
i_{x*}:\Spec(H^0(\bar{B}(\sN^\efc(k))))\to \Spec(H^0(\bar{B}(\sN^\efc(X)))).
\]
By lemma~\ref{lem:Kernel}  and theorem~\ref{thm:Kernel}, we have the identification 
\[
K_x\cong \Spec(H^0_{\sN^\efc(k)}(\bar{B}_{\sN^\efc(k)}(\sN^\efc(X),i_x^*)))
\]
as group schemes in $\sH_{\sN^\efc(k)}$, hence as pro-group schemes  in  $\sH^f_{\sN^\efc(k)}$.

But by theorem~\ref{thm:Main1}, the equivalence
\[
H^0(\sM^f):\sH^f_{\sN^\efc(k)}\to \MT(k)
\]
identifies $\Spec(H^0_{\sN^\efc(k)}(\bar{B}_{\sN^\efc(k)}(\sN^\efc(X),i_x^*)))$ with $\pi^\mot_1(X,x)$, completing the proof.
\end{proof}

\begin{corollary} \label{cor:main} Let $k$ be a number field and $S\subset \P^1(k)$ a finite set of $k$-points of $\P^1$. Set $X:=\P^1_k\setminus S$ and let $a\in X(k)$ be a $k$-point. Then both $k$ and  $X$ satisfy the Beilinson-Soul\'e vanishing conjectures. Furthermore, there is an isomorphism 
 \[
\pi^\mot_1(X,a)\cong K_a
\]
as pro-group objects in $\MT(k)$.
\end{corollary}

\begin{proof} $k$ satisfies the Beilinson-Soul\'e vanishing conjectures by Borel's theorem on the rational $K$-groups of $k$ \cite{Borel}. For $X$, we have the Gysin distinguished triangle
\[
M_\gm(X)\to M_\gm(\P^1)\to \oplus_{x\in S}\Z(1)[2]\to M_\gm(X)[1].
\]
Taking motivic cohomology gives the long exact sequence
\begin{multline*}
\ldots\to \oplus_{x\in S}H^{p-2}(k,\Z(q-1)) \to H^p(k,\Z(q))\oplus H^{p-2}(k,\Z(q-1))\\
\to H^p(X,\Z(q))\xrightarrow{\partial} \oplus_{x\in S}H^{p-1}(k,\Z(q-1))\to\ldots
\end{multline*}
Thus the vanishing conjectures for $k$ imply the vanishing conjectures for $X$. In addition, since $M_\gm(\P^1)=\Z\oplus \Z(1)[2]$, the Gysin exact triangle shows that $M_\gm(X)_\Q$ is in $\DTM(k)$, and thus the dual $\nh_k(X)=M_\gm(X)_\Q^\vee$ is also in $\DTM(k)$. 

We may therefore apply theorem~\ref{thm:Main} to give the isomorphism 
\[
\pi^\mot_1(X,a)\cong K_a.
\]
\end{proof}

\end{document}